\documentclass[12pt]{article}

\usepackage{makeidx}
\usepackage{graphics,amsmath,amssymb,amscd}
\usepackage{amscd}
\usepackage[matrix,arrow,curve]{xy}
\usepackage[dvips]{graphicx}

\usepackage{mathptmx}
\usepackage{helvet}
\usepackage{courier}
\usepackage{type1cm}

\usepackage{multicol}        
\usepackage[bottom]{footmisc}

\makeindex

\makeatletter

\@addtoreset{equation}{section}

\makeatother

\newtheorem{theo}{Theorem}[section]

\newtheorem{df}{Definition}[section]
\newtheorem{prop}{Proposition}[section]
\newtheorem{cor}{Corollary}[section]
\newtheorem{rem}{Remark}[section]
\newtheorem{lem}{Lemma}[section]
\newtheorem{conj}{Conjecture}[section]

\newcommand{\bda}{{\bf a}}
\newcommand{\bdb}{{\bf b}}
\newcommand{\bdk}{{\bf k}}

\newcommand{\bdu}{{\bf u}}

\newcommand{\bdx}{{\bf x}}
\newcommand{\bdy}{{\bf y}}
\newcommand{\bdz}{{\bf z}}

\newcommand{\bdH}{{\bf H}}

\newcommand{\eps}{\varepsilon}

\newcommand{\thet}{\vartheta}
\newcommand{\mR}{{\mathbb R}}
\newcommand{\mZ}{{\mathbb Z}}
\newcommand{\mT}{{\mathbb T}}

\newcommand{\mN}{{\mathbb N}}

\newcommand{\calA}{{\mathcal A}}
\newcommand{\calB}{{\mathcal B}}

\newcommand{\calD}{{\mathcal D}}
\newcommand{\calF}{{\mathcal F}}

\newcommand{\calK}{{\mathcal K}}
\newcommand{\calL}{{\mathcal L}}
\newcommand{\calO}{{\mathcal O}}
\newcommand{\calP}{{\mathcal P}}
\newcommand{\calR}{{\mathcal R}}
\newcommand{\calT}{{\mathcal T}}
\newcommand{\calV}{{\mathcal V}}
\newcommand{\calZ}{{\mathcal Z}}

\newcommand{\loc}{{\mathrm loc}}
\newcommand{\topo}{{\mathrm top}}

\newcommand{\Cr}{\mbox{\rm Cr}}
\newcommand{\Int}{\operatorname{Int}}

\newcommand{\dist}{\operatorname{dist}}

\newcommand{\SM}{\rm {SM}}

\newcommand\qed{{\unskip\nobreak\hfil\penalty50
 \hskip2em\hbox{}\nobreak\hfil\mbox{\rule{1ex}{1ex} \qquad}
  \parfillskip=0pt \finalhyphendemerits=0\par\medskip}}

\renewcommand{\Re}{\mbox{\rm Re}\,}

\begin{document}
\author{Sergey Bolotin\\
University of Wisconsin--Madison and\\ Moscow Steklov Mathematical Institute \and
Dmitry Treschev\\ Moscow Steklov Mathematical Institute }
\title{Anti-integrable limit}

\date{}

\maketitle

\begin{abstract}
  Anti-integrable limit is one of convenient and relatively simple methods for construction
  of  chaotic  hyperbolic invariant sets in Lagrangian, Hamiltonian and other dynamical systems.
  We discuss the most natural context of the method  -- discrete Lagrangian systems.
  Then we present examples and applications.
\end{abstract}

\tableofcontents

\section{Introduction}

Development of mathematics is inhomogeneous and unpredictable.
Unexpected wonderful breakthroughs alternate with strange delays.
Such a delay, difficult to explain, concerns the appearance of the
general concept of anti-integrable (AI) limit. Theory of dynamical
systems was technically ready for such a concept about 100 years
ago partially and 50 years ago completely. The main ideas are
contained in works of Poincar\'e, Birkhoff, Hedlund, Morse, and
other creators of symbolic dynamics. Remarks on hyperbolicity of
the invariant set that appears in AI limit are more or less
straightforward. They need just the general concept of
hyperbolicity and standard technical tools such as the cone
criterium of hyperbolicity. However, anti-integrable limit as a
universal approach appeared only in the end of the last century.

A starting point for the formation of the ideology of AI limit is
the paper by Aubry and Abramovici \cite{aubry0}, where the term
``anti-integrable limit'' was introduced and the main features of
the method were fixed: the context (discrete Lagrangian systems),
the main technical tool (the contraction principle) and the
language for the presentation of the results (symbolic dynamics).

Recall that the standard map is the area-preserving self-map of
the cylinder $\{(x,y) : x\in\mT, y\in\mR\}$, $\mT =
\mR/(2\pi\mZ)$, where
\begin{equation}
\label{standard:ham}
    \Big(\! \begin{array}{c}
               x \\ y
            \end{array}
    \!\Big)
  \mapsto
    \Big(\! \begin{array}{c}
               x_+ \\ y_+
            \end{array}
    \!\Big)
  = \Big(\! \begin{array}{c}
             x + y +  \lambda\sin x \\ y + \lambda\sin x
            \end{array}
    \!\Big).
\end{equation}
If compactness is needed, it is possible to take as the phase
space the torus $\mT^2$.

For $\lambda=0$ the map is integrable: the variable $y$ remains
constant on trajectories. The case of small $|\lambda|$ was
studied extensively in physical and mathematical literature. The
limit $\lambda\to\infty$ is called the AI limit \cite{aubry0}. The
main result of \cite{aubry0} is the construction for large
$|\lambda|$ of an uncountable set of trajectories which are in
one-to-one correspondence with a certain set of quasi-trajectories
(codes).

AI limit appeared as a realization of the simple idea that it is
natural to expect chaotic properties to become more pronounced
when the ``distance to the set of integrable systems'', whatever
this means, increases.

The absence of an analytic first integral in systems close to the AI limit is intuitively obvious. However, a formal proof of this fact needs some work.\footnote
 {In the case of two-dimensional phase space
 nonintegrability follows from the existence
 of transversal homoclinics to periodic solutions.}
On the other hand one should not think that an unbounded growth of chaos in the anti-integrable limit leads to ergodicity of the system. Indeed, the number of elliptic periodic points can be arbitrarily large for large values of $\lambda$ and moreover, for an increasing sequence of the parameter, these points are asymptotically dense on the phase torus, \cite{Dua}. Since in a general situation elliptic periodic trajectories are surrounded by stability islands, the standard map does not become ergodic for large $|\lambda|$. Although
numerically for large values of the parameter $\lambda$ the dynamics looks more and more chaotic, rigorously it is even not known if the measure of the chaotic set is positive. For example, if the metric entropy of the standard map is positive for at least one $\lambda$?

It is easy to show that trajectories constructed in AI limit form a {\it uniformly hyperbolic Cantor set}. Although this set is uncountable, its measure vanishes. A standard way of studying chaos quantitatively is computing (or estimating) the topological entropy which unlike the metric entropy admits relatively simple positive lower estimates. These estimates are based on the
standard fact that if a system has an invariant set with dynamics
conjugated to a topological Markov chain (TMC) then the topological
entropy of the system is not less than that of the TMC. In the case of the standard map (\ref{standard:ham}) this argument gives that for large $\lambda$ the topological entropy is greater than $c\log\lambda$ for some positive constant $c$, compare with \cite{Knill}.
\smallskip

The plan of the paper is as follows. In Section \ref{sec:standard} we give main ideas, methods and results of AI limit in the simplest nontrivial example: the standard map. We define the corresponding set of quasi-trajectories (codes) and prove that for each code there is a unique trajectory of the standard map shadowing this quasi-trajectory. The set of trajectories obtained in this way form a hyperbolic set. We show that topological entropy supported by this set is of order $\log\lambda$. In Section \ref{sec:DLS} we define the class of systems for which methods of the AI limit will be developed. These systems are called
the discrete Lagrangian systems (DLS). We start from globally defined DLS (Section \ref{ssec:DLS}) and then explain how discrete Lagrangian systems can be generated by the so called ambient systems. In section \ref{ssec:DLS_map} the ambient systems are symplectic maps.

Then we give examples of systems which can be studied by methods
of AI limit. In each example we have a large parameter, an analog
of $\lambda$ in the standard map. We tried to arrange the examples
in the order of increasing complexity. The first portion of
examples (Section \ref{sec:exa_discr}) contains systems with
discrete time. We start with the problem on the motion of a light particle in the field of a potential with $\delta$-like dependence on $t$. Then we consider a billiard system in a strip with walls, formed by graphs of periodic functions, the large parameter is the distance between the walls. After that we consider the billiard in the domain with small scatterers, where the large parameter is
the inverse of the scatterers size. Then we use AI limit to prove
several well known (and less well known) shadowing theorems,
starting with Shilnikov description of symbolic dynamics near a
transverse homoclinic orbit and ending with applications to
Celestial Mechanics. In the last section we discuss application of
the AI limit to the problem of Arnold diffusion.

Many examples presented in this survey are classical and go back to the dynamical folklore, so we do not always give references to the original results.

\section{The standard map}
\label{sec:standard}

Let us rewrite the map (\ref{standard:ham}) in the ``Lagrangian
form''. To this end suppose that
$$
  \left(\begin{array}{c} x_- \\ y_- \end{array}\right)
  \mapsto
  \left(\begin{array}{c} x  \\  y  \end{array}\right)
  \mapsto
  \left(\begin{array}{c} x_+ \\ y_+ \end{array}\right) .
$$
Then $x_-,x,x_+$ satisfy the equation
\begin{equation}
\label{eq:sml}
  \lambda^{-1} (x_+ - 2x + x_-) = \sin x.
\end{equation}
The standard map $\SM$ written in this form is defined on the
cylinder ${\calZ} = \mR^2 / \sim$, where the equivalence relation
$\sim$ is as follows:
\begin{equation}
\label{id}
  (x'_1,x'_2) \sim (x_1,x_2) \quad
  \mbox{if and only if} \quad
  x'_1 - x_1 = x'_2 - x_2 \in 2\pi\mZ.
\end{equation}
In the other words, the cylinder $\calZ$ is the quotient space of
the plane $\mR^2$ with respect to the action of the group of
shifts
$$
  (x_1,x_2) \mapsto (x_1 + 2\pi l, x_2 + 2\pi l),\qquad l\in \mZ.
$$
The map $\SM$ acts as follows:
\begin{equation}
\label{SM_lagr}
  (x_-,x) \mapsto \SM(x_-,x) = (x,x_+), \quad
  \mbox{where $x_-,x,x_+$  satisfy~(\ref{eq:sml})}.
\end{equation}
Infinite sequences $\bdx=(x_k)_{k\in\mZ}$ such that the triple
$(x_-,x,x_+) = (x_{k-1},x_k,x_{k+1})$ satisfies (\ref{eq:sml}) for
any integer $l$, are called trajectories of the standard map.

The Lagrangian form of the standard map admits a variational
formulation. Namely, trajectories of the system are critical
points (extremals) of the formal sum
\begin{equation}
\label{action} A(\bdx)=  \sum_{k\in\mZ} L(x_k,x_{k+1}),\qquad
  L(x',x'') = \frac{1}{2\lambda} (x'-x'')^2 - \cos x'' .
\end{equation}
This means that $\bdx^0=(x_k^0)_{k\in\mZ}$ is a trajectory if and
only if for any integer $n$
\begin{equation}
\label{dA}
  \frac{\partial}{\partial x_n}
    \sum_{k = -\infty}^\infty  L(x_k, x_{k+1}) = 0 \quad
  \mbox{at the point} \quad
  \bdx = \bdx^0.
\end{equation}

In the limit $\lambda\to\infty$, the standard map becomes
dynamically meaningless because $x_+$ cannot be found in terms of
$x$ and $x_-$ from equation~(\ref{eq:sml}) when $\lambda^{-1}=0$.
However, the corresponding variational problem remains well-defined. Its solutions are sequences
\begin{equation}
\label{singtra}
  \bda=(a_j)_{j\in\mZ}, \qquad
  a_j \in\pi\mZ.
\end{equation}
For large values of the parameter $\lambda$ the standard map has
many trajectories close to sequences~(\ref{singtra}).

More precisely,  let
$$
    \calA_\Lambda
  = \{ \bda = (a_k)_{k\in\mZ} : a_k \in\pi\mZ, \;
         |a_{k-1} - 2a_k + a_{k+1}| \le \Lambda \;
         \mbox{ for any $k\in\mZ$
         } \}.
$$
For any $\bda\in\calA_\Lambda$ we define the complete metric space
$\Pi = \Pi(\bda)$ of the sequences
$$
  \bdx = \{x_k\}_{k\in\mZ}, \qquad
  \sup\limits_{k\in\mZ} |x_k - a_k| \le \pi/2 .
$$
The metric on $\Pi$ is defined as follows:
$$
     \rho(\bdx',\bdx'')
   = \sup\limits_{k\in\mZ} |x'_k - x''_k|,\qquad
     \bdx',\bdx'' \in \Pi.
$$

\begin{theo}[\cite{aubry0}]
\label{theo:aubry} Given $\Lambda>0$, $\sigma\in (0,\pi/2)$, and
\begin{equation}
\label{lambda0}
    \lambda_0
  = \lambda_0(\Lambda,\sigma)
  = \max\Big\{ \frac{\Lambda+4\sigma}{\sin\sigma}, \frac{8}{\cos\sigma}
        \Big\}
\end{equation}
for any $|\lambda| \ge \lambda_0$ and any $\bda\in\calA_\Lambda$,
the standard map has a unique trajectory $\bdx$ with
$\rho(\bda,\bdx) < \sigma$.
\end{theo}

Sequences from $\calA_\Lambda$ can be regarded as codes of the
corresponding trajectories. This possibility to code trajectories
by elements of a sufficiently large set is typical for chaotic
systems.

Theorem~\ref{theo:aubry} means that for large values of $\lambda$
there is an invariant set $\calK_\Lambda =
\calK_\Lambda(\lambda,\sigma)\subset\calZ$ such that trajectories
from $\calK_\Lambda$ are in a one-to-one correspondence with
elements of $\calA_\Lambda$. Formal definition of $\calK_\Lambda$
is as follows. Any code $\bda\in\calA_\Lambda$ uniquely determines
an orbit $\bdx = \bdx(\bda)$. Consider the map
$$
          \calA_\Lambda\ni\bda
  \mapsto \zeta(\bda)
      =   (x_0,x_1)/\mZ \in \calZ,
$$
where $(x_0,x_1)/\mZ$ means identification (\ref{id}). Then by
definition $\calK_\Lambda = \zeta(\calA_\Lambda)$.

All trajectories in the set $\calK_\Lambda$ are hyperbolic
(Theorem \ref{theo:sm:hyp}). Taking another $\mZ$-quotient, we
make the phase space of $\SM$ compact: $\calZ/\mZ$ is
diffeomorphic to the 2-torus $\mT^2$. Then we prove that
$K_\Lambda = \calK_\Lambda / \mZ \subset\mT^2$ is compact (Theorem
\ref{theo:closed}) and therefore, hyperbolic for the quotient map
$\SM_0$. Finally we estimate topological entropy of the
corresponding dynamical system (Theorem \ref{theo:entropy}).
\medskip

\noindent {\it Proof of Theorem \ref{theo:aubry}}. The proof is
based on the contraction principle in the metric space
$(\Pi(\bda),\rho)$, where $\bda\in\calA_\Lambda$ is fixed. First
we rewrite equations (\ref{eq:sml}) in the form
\begin{equation}
\label{Chi_ai}
    x_k
  = \arcsin_k \Big( \frac{x_{k+1} - 2x_k + x_{k-1}}
                         {\lambda}
              \Big) ,
\end{equation}
where $\arcsin_k$ is the branch of $\sin^{-1}$ such that
$\arcsin_k(0) = a_k$. Thus $\arcsin_k$ maps the interval $(-1,1)$
to the interval $\big(a_k- \frac\pi 2,a_k+\frac\pi 2\big)$. In
other words $\arcsin_k x=a_k\pm\arcsin x$, where $\arcsin =
\arcsin_0:(-1,1)\to(-\pi/2,\pi/2)$ is the standard branch of
$\sin^{-1}$ and  $+$ $(-)$ is taken when $a_k/\pi$ is even (odd).
The sequence $\bdx=\bda$ is a trajectory of equations
(\ref{Chi_ai}) for $\lambda = \infty$.

 Consider the map
$\bdx\mapsto \bdy = \Phi(\bdx)$ such that
$$
    y_k
  = \arcsin_k \Big( \frac{x_{k+1} - 2x_k + x_{k-1}}
                         {\lambda}
              \Big) .
$$
Then any fixed point of $\Phi$ is a trajectory of $\SM$.

Let $B_\sigma(\bda)\subset\Pi$ be the closed ball with center
$\bda$ and radius $\sigma$.

\begin{lem}
\label{lem:ai} Suppose that $|\lambda| > \lambda_0$. Then
\begin{enumerate}
\item $\Phi$ is defined on $B = B_\sigma(\bda)$ and
$\Phi(B)\subset B$; \item $\Phi$ is a contraction on $B$ i.e.,
\begin{equation}
\label{eq:contract}
    \rho\big( \Phi(\bdx'),\Phi(\bdx'') \big)
  < \frac 12 \,\rho(\bdx',\bdx'') \qquad
    \mbox{for all}\quad \bdx',\bdx''\in B.
\end{equation}
\end{enumerate}
\end{lem}

By contraction principle, Theorem~\ref{theo:aubry} follows from
Lemma \ref{lem:ai}. Now we turn to the proof of the lemma. Recall
that $\sigma < \pi / 2$. To check that $\Phi(B)\subset B$, it is
sufficient to show that for any $\bdx\in B$
\begin{equation}
\label{||<sin}
    \Big| \frac{x_{k+1} - 2x_k + x_{k-1}}{\lambda}
    \Big|
  < \sin\sigma .
\end{equation}
Since $\rho(\bdx,\bda) < \sigma$, we have:
$$
      |x_{k+1} - 2x_k + x_{k-1}|
  \le \Lambda + 4\sigma.
$$
Hence inequality (\ref{||<sin}) holds for
$$
  \lambda_0 > \frac{\Lambda + 4\sigma}{\sin\sigma} .
$$

Note that for any pair of real numbers $u',u'' \in
(-\sin\sigma,\sin\sigma)$
$$
      \big| \arcsin u' - \arcsin u''
      \big|
  \le \frac 1{\cos\sigma} |u' - u''| .
$$
Here
$$
\frac 1{\cos\sigma} = \sup\limits_{u\in (-\sin\sigma,\sin\sigma)}
                       \Big| \frac d{du} \arcsin u \Big|.
$$

We put $\bdy' = \Phi(\bdx')$, $\bdy'' = \Phi(\bdx'')$. Then for
any $k\in\mZ$ we have:
\begin{eqnarray}
\nonumber
       |y'_k - y''_k|
 &=&  \bigg|
       \arcsin \Big( \frac{x'_{k+1} - 2x'_k + x'_{k-1}}{\lambda}
                 \Big)
     - \arcsin \Big( \frac{x''_{k+1} - 2x''_k + x''_{k-1}}{\lambda}
                 \Big)
      \bigg| \\
\nonumber
 &\le& \frac 1{\cos\sigma} \,
       \Big| \frac{x'_{k+1} - 2x'_k + x'_{k-1}}{\lambda}
            - \frac{x''_{k+1} - 2x''_k + x''_{k-1}}{\lambda}
       \Big| \\
\nonumber
 &\le& \frac{   |x'_{k+1} - x''_{k+1}|
               + 2 |x'_k - x''_k|
               + |x'_{k-1} - x''_{k-1}|}
            {\lambda\cos\sigma} \\
\label{y-y}
 &\le& \frac 4{\lambda\cos\sigma} \rho(\bdx',\bdx'').
\end{eqnarray}
Hence inequality (\ref{eq:contract}) holds if $\lambda_0 > 8 /
\cos\sigma$. Lemma \ref{lem:ai} is proved. \qed

\begin{theo}
\label{theo:sm:hyp} Given $\Lambda>0$, $\sigma\in (0,\pi/2)$, and
$\lambda$, where $|\lambda| > \lambda_0(\Lambda,\sigma)$. Then any
orbit $\bdx(\bda)$, $\bda\in\calA_\Lambda$ is hyperbolic.
\end{theo}

\noindent {\it Proof}. Recall that the standard definition of a
hyperbolic orbit $\bdx = \bdx(\bda)$ is given in terms of an
invariant decomposition into expanding and contracting subspaces
of tangent spaces at any point \cite{KH}. Here we use instead the
{\it cone criterium} of hyperbolicity of V.M.~Alexeyev, also
contained in \cite{KH}.

In this proof for brevity we use the notation $f = \SM$. By
(\ref{eq:sml}) the Jacobi matrix of the map $f$ is as follows:
$$
     Df(x_-,x)
  =  \Big( \begin{array}{rc} 0   &   1  \\
                            -1   &   2 + \lambda\cos x
           \end{array}
     \Big).
$$
At each point $q=(x_-,x)\in K_\Lambda$ we define the cones
\begin{equation}
\label{HV}
    H_q
  = \big\{(u_-,u)\in T_q\calZ : \|u_-\| \le \alpha_H \|u\| \big\}, \;\;
    V_q
  = \big\{(u_-,u)\in T_q\calZ : \|u\| \le \alpha_V \|u_-\| \big\},
\end{equation}
where $(u_-,u)$ are coordinates on $T_q\calZ$ associated with the
coordinates $(x_-,x)$ on $\calZ$.

By the cone criterium to prove hyperbolicity of trajectories in the
invariant set $\calK_\Lambda$, it is sufficient to check that
there exists $\mu>1$ such that for all $q\in\calK_\Lambda$ the
following holds:
\begin{eqnarray}
\label{cone1} &  Df_q H_q \subset \Int H_{f(q)}, \quad
   Df_q^{-1} V_{f(q)} \subset \Int V_q, & \\
\label{cone2} & \mbox{
 $\|Df_q\xi\| \ge \mu \|\xi\|$ for $\xi\in H_q$, and
 $\|Df_q^{-1}\xi\| \ge \mu \|\xi\|$ for $\xi\in V_{f(q)}$}. &
\end{eqnarray}

We take $\alpha_H = \alpha_V = 1/2$. For any $q\in
K_\Lambda(\lambda,\sigma)$ and $(u_-,u)\in H_q$ we have:
$$
    Df(q) \Big(\begin{array}{c} u_- \\ u \end{array}\Big)
  = \Big(\begin{array}{c} u \\ u_+ \end{array}\Big), \qquad
    u_+ = -u_- + (2 + \lambda\cos x) u.
$$
We see that
$$
  \frac{|u|}{|u_+|} \le \frac{1}{\lambda\cos\sigma - 5/2}.
$$
Therefore
$$
  \Big(\begin{array}{c} u \\ u_+ \end{array}\Big)
  \in \mbox{Int}\, H_{f(q)} \quad
  \mbox{provided}\quad
  \lambda > \frac{3}{\cos\sigma}.
$$
This implies the first inclusion (\ref{cone1}). We also have:
$$
     \Big\|\Big(\begin{array}{c} u \\ u_+ \end{array}\Big)\Big\|^2
  =  u^2 + (\lambda\cos x - 5/2)^2 u^2
 \ge 4 \Big\|\Big(\begin{array}{c} u_- \\ u \end{array}\Big)\Big\|^2
$$
provided $\lambda > 9 / (2\cos\sigma)$. Hence the first inequality
(\ref{cone2}) holds with $\mu = 2$.

The second inclusion (\ref{cone1}) and the second inequality
(\ref{cone2}) can be checked analogously. It remains to note that
$\lambda_0 > 9/(2\cos\sigma)$. \qed

To make chaotic properties of the map $\SM$ quantitative, we
estimate its topological entropy. Since in the case of noncompact
phase space topological entropy is usually infinite, we come to a
quotient-system with a compact phase space. To this end we note
that equations (\ref{standard:ham}) can be considered modulo $2\pi$. This generates the quotient system
$$
  \SM_0 : \mT^2\to\mT^2, \qquad
  \mT^2 = \mR^2/(2\pi\mZ^2)
$$
such that the diagram
$$
  \begin{CD}
  \calZ  @>{\SM}>>  \calZ \\
  @V{\mbox{pr}}VV                   @VV{\mbox{pr}}V    \\
  \mT^2  @>{\SM_0}>>  \mT^2
  \end{CD}\quad , \qquad
  \mbox{pr} : (x_1,x_2) \mapsto (x_1\bmod 2\pi, x_2\bmod 2\pi).
$$
is commutative.

Consider on $\calA_\Lambda$ and on the space $\calO$ of orbits of
$\SM$ the following action of the group $\mZ$: for any
$(l_1,l_2)\in\mZ^2$
\begin{eqnarray*}
   \bda  =    (a_i)_{i\in\mZ}
    &\mapsto& l(\bda)
        =    (a_i + 2\pi l_1 + 2\pi il_2)_{i\in\mZ},  \\
  \calO \ni  \bdx
        =    (x_i)_{i\in\mZ}
    &\mapsto& l(\bdx)
        =    (x_i + 2\pi l_1 + 2\pi il_2)_{i\in\mZ}.
\end{eqnarray*}
We define the quotient spaces
$$
  A_\Lambda = \calA_\Lambda/\mZ^2, \quad
  O = \calO/\mZ^2, \quad
  K_\Lambda = \mbox{pr}\,\calK_\Lambda .
$$
Then $A_\Lambda$ can be regarded as the space of codes while $O$
as the corresponding space of orbits for $\SM_0$. The space
$A_\Lambda$  can be identified with
$$
  B_\Lambda = \{\bdb = (b_i)_{i\in\mZ} : b_i\in\pi\mZ,\; |b_i|\le\Lambda\},
$$
via the bijection
$$
  \bda = (a_i)_{i\in\mZ} \mapsto \bdb = (b_i)_{i\in\mZ}, \qquad
  b_i = a_{i-1} - 2a_i + a_{i+1}
$$
which respects the shift operator $\calT:(a_i)_{i\in\mZ}\mapsto
(a_{i+1})_{i\in\mZ}$.

We recall the definition of the Bernoulli shift. Let $J$ be a
finite set of $q$ symbols and $\Sigma_q=J^\mZ$ the set of
sequences $\bda=(a_i)_{i\in \mZ}$ of elements of $J$. We equip
$\Sigma_q$ with the product (Tikhonov) topology: the base consists
of cylinders
\begin{equation}
\label{top_base}
  U_I(\bda)=\{\bda':a'_i=a_i\;\mbox{ for } i\notin I\},\quad
  \mbox{where $I\subset\mZ$ is a finite set.}
\end{equation}
It is well known \cite{KH} that $\Sigma_q$ is a compact totally
disconnected space with no isolated points and hence homeomorphic
to the standard Cantor set. The Bernoulli shift with $q$ symbols
is the shift $\calT:\Sigma_q\to\Sigma_q$.

Let
\begin{equation}
\label{eq:q}
  q = \#(\pi\mZ\cap[-\Lambda,\Lambda])=1 + 2[\Lambda/\pi]
\end{equation}
where $[\;]$ is  the integer part of a real number.

\begin{theo}
\label{theo:closed} For any $\lambda>\lambda_0$ the set $K_\Lambda$ is
a compact hyperbolic invariant set for the standard map and
$\SM_0:K_\Lambda\to K_\Lambda$ is topologically conjugate to the
Bernoulli shift on the space of $q$ symbols.
\end{theo}

\noindent {\it Proof}. We only need to show that
$\zeta:A_\Lambda\to K_\Lambda$ is continuous. The product topology
on $B_\Lambda$ has the base of open sets (\ref{top_base}), where
$J = B_\Lambda$. In particular, for any $\bda\in B_\Lambda$, the
cylinder
$$
  C_n(\bda)=\{\bda'\in B_\Lambda: a'_i=a_i\;\mbox{for}\;|i|\le n\}
$$
is an open set containing $\bda$. Continuity of $\zeta$ follows
from continuity of the map $\bda\mapsto\bdx(\bda)$ which is given
by the following lemma.

\begin{lem}
\label{lem:a'=a''} Suppose that $|\lambda|\ge \lambda_0$ and
$\bda'\in C_n(\bda)$. Then the orbits $\bdx = \bdx(\bda)$ and
$\bdx' = \bdx(\bda')$ satisfy the estimates
\begin{equation}
\label{x'-x''<5}
  |x'_k - x_k| \le  5^{|k|-n}\cdot 2\sigma \;
  \mbox{ for any } |k| \le n.
\end{equation}
\end{lem}

\noindent {\it Proof of Lemma \ref{lem:a'=a''}}. We obtain the
sequences $\bdx$ and $\bdx'$ as the limits
\begin{eqnarray*}
& \bdx = \lim_{j\to\infty} \bdx^{(j)}, \quad
  \bdx' = \lim_{j\to\infty} \bdx^{\prime (j)}, & \\
& \bdx^{(0)} = \bda, \quad
  \bdx^{\prime (0)} = \bda', \quad
  \bdx^{(j+1)} = \Phi\bdx^{(j)}, \quad
  \bdx^{\prime (j+1)} = \Phi\bdx^{\prime (j)} , &
\end{eqnarray*}
where $\Phi$ is the operator from Lemma \ref{lem:ai}.

We use induction in $j$. For $j = 0$ inequalities (\ref{x'-x''<5})
hold because $x^{(0)}_k = x^{\prime (0)}_k$ for $|k|\le n$. If
$k=n$ then (\ref{x'-x''<5}) for all $j>0$ follows from assertion 1
of Lemma \ref{lem:ai}.

Suppose that (\ref{x'-x''<5}) holds for some $j = s$. To prove it
for $j = s+1$, we will use the estimate
$$
      |\bdx^{(j+1)}_k - \bdx^{\prime (j+1)}_k|
  \le \frac{|\bdx^{(j)}_{k+1} - \bdx^{\prime (j)}_{k+1}|
            + 2 |\bdx^{(j)}_k - \bdx^{\prime (j)}_k|
            +  |\bdx^{(j)}_{k-1} - \bdx^{\prime (j)}_{k-1}|}
           {\lambda\cos\sigma}
$$
which follows from (\ref{y-y}). For any $|k| \le n - 1$ by the
induction assumption we have:
$$
      |\bdx^{(s+1)}_k - \bdx^{\prime (s+1)}_k|
  \le \frac{(5^{- n + |k| + 1}
            + 2\cdot 5^{- n + |k|}
            + 5^{- n + |k| - 1})\cdot 2\sigma}
           {\lambda\cos\sigma}.
$$
The right-hand side of this inequality does not exceed $5^{-n +
|k|}\cdot 2\sigma$ because by (\ref{lambda0}) $\lambda_0\cos\sigma
\ge 8$. \qed

\medskip

In \cite{Gel-Tur} there is another proof of a similar estimate.

\begin{theo}
\label{theo:entropy} There exists a constant $c > 0$ such that for sufficiently large $\lambda$ topological entropy of the map $\SM_0$ satisfies the estimate\footnote {Explicit estimates for $c$ and minimal $\lambda$  can be easily obtained from the proof, see also \cite{Knill}. However we  do not expect that these estimates are close to optimal.}
$$
  h_\topo(\SM_0) \ge c\log\lambda.
$$
\end{theo}

\noindent {\it Proof}. We already proved that for sufficiently
large $\lambda$ the map $\SM_0$ has an invariant set $K_\Lambda$
such that the restriction $\SM_0|_{K_\Lambda}$ is conjugated to
the Bernoulli shift $\calT:\Sigma_q\to\Sigma_q$ on the space of $q$ symbols given by (\ref{eq:q}). Hence $h_\topo(\SM_0)$ is not less than
topological entropy $h_q$ of   the Bernoulli shift \cite{KH}. The
quantity $h_q$ can be computed for example, from the equation,
\cite{Walters} (which holds also for topological Markov chains)
\begin{equation}
\label{top_ent}
  h_q = \lim_{n\to\infty} \frac1n \log\theta_n,
\end{equation}
where $\theta_n$ is the cardinality of all admissible
$n$-sequences. In this case $\theta_n = q^n$. Therefore $h_q =
\log q$.

It remains to note that according to  (\ref{lambda0}) the quantity
$\Lambda$ can be chosen greater than a positive constant
multiplied by $\lambda$. \qed

\section{Discrete Lagrangian systems}
\label{sec:DLS}

Anti-integrable limit is usually discussed in Lagrangian systems.
 We will concentrate on the Lagrangian systems with discrete time
which are called {\it discrete Lagrangian systems} (DLS). The case
of continuous time will be reduced to the case of discrete time,
Section \ref{sec:ail_cls}.

\subsection{The simplest case}
\label{ssec:DLS}

Let $M$ be an $m$-dimensional manifold and $L$ a smooth\footnote
{$C^3$ is more than enough.} function on $M^2=M\times M$. By
definition a sequence $\bdx = (x_i)_{i\in\mZ}$, $x_i\in M$, is a
trajectory of the  discrete Lagrangian system with the Lagrangian
$L$ if $\bdx$ is an extremal of the action functional which is the
formal sum
$$
  A(\bdx) = \sum L(x_i,x_{i+1})
$$
in the same sense as this was explained   for the standard map,
see (\ref{dA}). Equivalently, for any $i\in\mZ$
\begin{equation}
\label{-0+}
  \partial_{x_i}\big( L(x_{i-1},x_i) + L(x_i,x_{i+1})\big) = 0.
\end{equation}
Trajectories remain the same after multiplication of the
Lagrangian by a constant, after addition of a constant to $L$ and
after the  gauge transformation
$$
  L(x,y) \mapsto L(x,y) + f(x) - f(y)
$$
with an arbitrary smooth function $f$ on $M$. The gauge
transformation does not change the action functional.

If the equation
$$
  \partial_x (L(x_-,x)+L(x,x_+)) = 0
$$
can be globally solved for $x_- = x_-(x,x_+)$ as well as for
$x_+=x_+(x_-,x)$, then the map
\begin{equation}
\label{T}
  (x_-,x) \mapsto T(x_-,x) = (x,x_+)
\end{equation}
is a diffeomorphism and determines a discrete dynamical system on
$M\times M$. However usually the map $T$ is only locally defined.

To present conditions for  the local existence and smoothness of
the map $T$, we define\footnote {Sometimes $B$ is defined as
$-\partial_x \partial_y L$, then  in many natural DLS $B$ will be
positive definite.}
$$
  B(x,y) = \partial_x \partial_y L(x,y).
$$
In local coordinates,
\begin{equation}
\label{eq:Bij}
    B(x,y)
  = \left(\frac{\partial^2 L}{\partial y_j\partial x_i}\right).
\end{equation}
In invariant terms, $B(x,y)$ is a linear operator $T_xM\to
T_y^*M$, or a bilinear form  on $T_xM\times T_yM$. We say that $L$
is a twist Lagrangian if it satisfies the following condition.
\medskip

\noindent {\bf Twist condition}: $B(x,y)$ is nondegenerate for all
$x,y\in M$.
\medskip

In this case the map $T$ is locally well defined and smooth. It is easy to
check (see e.g.\ \cite{Veselov}) that $T$ is symplectic with
respect to the symplectic 2-form $\omega=B(x,y)\,dx\wedge dy$,
\begin{equation}
\label{eq:omega} \omega(u,v)=\langle B(x,y)u_1,v_2\rangle-\langle
B(x,y)v_1,u_2\rangle,\qquad u=(u_1,u_2),\; v=(v_1,v_2).
\end{equation}

If the twist condition holds, then the Legendre transform
$S:M^2\to T^*M$,
$$
  (x,y)\mapsto (x,p_x),\qquad  p_x=-\partial_x L(x,y),
$$
is locally invertible and we can represent $T$ by a locally
defined map $F=STS^{-1}:T^*M\to T^*M$. The map $F$  is symplectic
with respect to the standard symplectic form $dp_x\wedge dx$ on
$T^*M$, and $L$ is the generating function of $F$:
\begin{equation}
\label{eq:gen}
  F(x,p_x) = (y,p_y)\quad\Leftrightarrow\quad
  p_x = -\partial_x L(x,y),\quad
  p_y = \partial_y L(x,y).
\end{equation}
Such symplectic map $F$ is  called a twist map, see \cite{Gole}.
Usually it is only locally defined.

The twist condition imposes strong topological restrictions on
$M$, it is rarely satisfied on all $M\times M$. Here are two
canonical examples of DLS.

{\bf 1}. The multidimensional standard map:
\begin{equation}
\label{eq:standard}
    L(x,y)
  = \frac12 \big\langle B(x-y),x-y\big\rangle
    - \frac12 \big(V(x) + W(y)\big),\qquad
     x,y\in\mR^m,
\end{equation}
where $B$ is a symmetric constant nondegenerate matrix. In this case $L$
defines a symplectic twist map $F: \mR^{2m}\to\mR^{2m}$.

Usually $V$ is assumed to be $\mZ^m$-periodic. Then the phase
space $\mR^{2m}$ can be factorized with respect to the
$\mZ^m$-action $(x,y)\mapsto (x+k,y+k)$, $k\in\mZ^m$, and  $F$ is
a symplectic twist self-map of $\mT^m\times\mR^m$.

We discuss such Lagrangians in more detail in section \ref{sec:equi}.
\medskip

{\bf 2}. Consider a  domain   $D\subset \mR^{m}$  bounded by a
smooth convex hypersurface $M$.  The  billiard system in $D$ is a
DLS with the Lagrangian $L(x,y) = |x-y|$ on $M\times M$. Let us
check that $L$ satisfies the twist condition on $(M\times
M)\setminus\Delta$, where $\Delta=\{(x,x):x\in M\}$. Let $\langle
B(x,y)v,w\rangle$ be the bilinear form on $\mR^m\times \mR^m$
corresponding to the operator $B(x,y) =
\partial_x\partial_y L(x,y)$. A direct computation gives
\begin{equation}
\label{eq:B_bill}
    \langle B(x,y)v,w\rangle
  = \frac{- \langle v,w\rangle
           + \langle v,e\rangle \langle w,e\rangle}
         {|x-y|},\qquad
     e
  =  \frac {x-y}{|x-y|}.
\end{equation}
Evidently, $e$ lies in both left and right kernel:
$B(x,y)e=B^*(x,y)e=0$. Hence, if $e\notin T_xM$ and $e\notin
T_yM$, then the restriction of the bilinear form $B(x,y)$ to
$T_xM\times T_yM$ is nondegenerate and hence $L$ is a twist
Lagrangian. This always holds if  the boundary is strictly convex.
This is well known, for a recent reference see \cite{Bia}.

In fact we may identify $T_xM$ and $T_yM$ by an isomorphism
$\Pi(x,y) : T_x M \to T_y M$ which is the parallel projection in
$\mR^{m}$ along the segment $[x,y]$: $\Pi v=v\bmod e$.
 Then
$$
    \langle B(x,y)v,\Pi(x,y)v\rangle
  = \frac{- |v|^2 + \langle v,e\rangle^2}{|x-y|} < 0,\qquad
    v\in T_xM\setminus\{0\}.
$$

We orient $M$ as the boundary. Since $\Pi(x,y)$ changes
orientation, we obtain that $(-1)^m\det B(x,y)>0$. Note that since image
and range of $B(x,y)$ are different, $\det B(x,y)$ is not
invariantly defined (at least as a number), but its sign is.

If $D$ is not convex, the billiard  map $T$ is   defined   on
$$
\{(x,y)\in (M\times M)\setminus\Delta : \mbox{ the segment $(x,y)$
is contained in $D$}\},
$$
and it is  singular when the segment $(x,y)$ is tangent to $M$.
\medskip

In \cite{Veselov} the reader can find other examples of DLS
(mostly integrable) including multivalued ones.

\subsection{Multivalued Lagrangians}
\label{ssec:DLS_multi}

In applications the discrete Lagrangian $L$ of  a DLS with the
configuration space $M$ is usually multivalued. A common situation
is when $L$ is a function on the universal covering of $M\times M$.
For example, for the standard map \ref{eq:standard} we can assume that $M=\mT^m$ and $L$ is a function
on $\mR^{2m}$. One can think of such $L$ as a multivalued
function, i.e.\ a collection of functions  on open sets in $M\times M$.
To cover most possible applications  we  will use the following
definition of a DLS.

A DLS with $m$ degrees of freedom is defined by a finite or
countable collection $\calL = \{L_\kappa\}_{\kappa\in J}$  of
functions (Lagrangians) on  $U_\kappa^-
\times U_\kappa^+$, where $U_\kappa^\pm$ are $m$-dimensional manifolds. Often $U_\kappa^\pm$ are open sets in the
configuration space $M$. We are also given a graph with
the set of vertices $J$ and the set of edges $E\subset J^2$.
The vertices $\kappa,\kappa'\in J$ are joined by an edge
$\gamma=(\kappa,\kappa')$ if $W_\gamma=U_\kappa^-\cap U_{\kappa'}^+\ne
\emptyset$.

Let $\bdk = (\kappa_i)_{i\in\mZ}$, where $\gamma_i=(\kappa_i,\kappa_{i+1})\in E$, be a path in the graph.
A trajectory of the DLS $\calL$ corresponding to  the code $\bdk $ is a
sequence
\begin{equation}
\label{eq:bdx}
  \bdx=(x_i)_{i\in\mZ}, \qquad
  x_i\in
       W_{\gamma_i},
\end{equation}
 which is a   critical point of  the formal discrete action functional
\begin{equation}
\label{eq:Ak}
   A_{\bdk}(\bdx)
 = \sum_{i\in\mZ} L_{\kappa_i}(x_i,x_{i+1}),\qquad \bdk=(\kappa_i)_{i\in\mZ}.
\end{equation}
More precisely,
\begin{equation}
\label{eq:crit}
    \partial_{x_i} (L_{\kappa_{i-1}}(x_{i-1},x_i) + L_{\kappa_{i}}(x_{i},x_{i+1}))
  = 0,\qquad
      i\in\mZ .
\end{equation}
Thus a trajectory is a pair $(\bdk,\bdx)$, where $\bdk\in J^\mZ$.

In this paper we do not use global methods, so choosing local
coordinates, we may assume that $W_\gamma$, $\gamma\in E$, are open sets in
$\mR^m$.

If the Lagrangians satisfy the twist condition
$$
  \det B_\kappa(x,y)\ne 0,\qquad
    B_\kappa(x,y)
  = \partial_x\partial_y L_{\kappa}(x,y),
$$
then  an edge $\gamma=(\kappa,\kappa')$ defines a local
diffeomorphism
$$
  f_{\gamma}:
    \calD_{\gamma}^- \to \calD_{\gamma}^+
$$
of open sets $\calD_{\gamma}^\pm$ in $\mR^{2m}$ by the equation
$$
f_\gamma(x,y)=(y,z)\quad\Leftrightarrow\quad \partial_{y}
(L_{\kappa}(x,y) + L_{\kappa'}(y,z))
  = 0
  $$

Then the  trajectory $(\bdk,\bdx)$ corresponds to a composition
$$
   (x_n,x_{n+1})
 = f_{\gamma_n}\circ f_{\gamma_{n-1}}\circ
            \cdots\circ  f_{\gamma_1}(x_0,x_1),\qquad \gamma_i=(\kappa_i,\kappa_{i-1}).
$$
It is customary to represent such dynamics by a single map $\calF$ by taking the skew product \cite{KH} of the maps $\{f_\gamma\}$.
If the sets $U_\kappa^\pm$ lie in the configuration space $M$, then $\calF$ a map of $(J\times M^2)^\mZ$. But often this is not necessary: we can represent dynamics by a single self-map of a smooth manifold $P$.

We say that a dynamical system $F:P\to P$ is ambient of a DLS $\calL$ if to every trajectory $(\bdk,\bdx)$ of $\calL$ there corresponds a trajectory $\bdz=(z_i)_{\in\in\mZ}$ of $F$. A   formal definitions is as follows:

\begin{df}
The dynamical system $F:P\to P$ is said to be ambient for the DLS $\calL$ if there exist two systems of open sets $\calP_{\gamma}^\pm\subset P$ and the diffeomorphisms
$\beta_{\gamma}^\pm : \calD_{\gamma}^\pm \to \calP_{\gamma}^\pm$ which conjugate $f_{\gamma}$ with $F$. In other words, if the following diagram commutes:
$$
 \begin{CD}
 \calD_{\gamma}^-
    @> f_{\gamma} >> \calD_{\gamma}^+ \\
  @V \beta_{\gamma}^- VV
               @VV \beta_{\gamma}^+ V    \\
 \calP_{\gamma}^-
    @> F|_{\calP_{\gamma}^-} >> \calP_{\gamma}^+
 \end{CD}
$$
\end{df}

Thus $\beta_\gamma^\pm$ are charts in $P$, and $f_\gamma$ is a
coordinate representation of $F$.

An ambient system exists under  weak conditions on the DLS, but in general it is not unique. For globally defined DLS from Section \ref{ssec:DLS} the ambient system  is the map $T:M\times M\to M\times M$. The problem of building of the ambient system is not important for our purposes because we usually start with a dynamical system $F:P\to P$ and the corresponding DLS $\calL$ is a
technical tool for studying dynamics of $F$. A standard example is
the Lagrangian representation of a symplectic map in the next
section.

\subsection{DLS generated by a symplectic map}
\label{ssec:DLS_map}

A DLS determines a local (in general) symplectic map. Conversely,
one can represent the dynamics of a symplectic map by a DLS, but
the Lagrangian will be only locally defined. Let
\begin{equation}
\label{sympl}
  F : P\to P, \quad
  F^*\omega = \omega,
\end{equation}
be a symplectic self-map of a symplectic $2m$-dimensional manifold
$(P,\omega)$ \cite{Arnold:matmex}.

Let $(x^-,y^-)$ and $(x^+,y^+)$ be symplectic coordinates in small
neighborhoods $D^-$ and $D^+$ of the points $z^-$ and $z^+ =
F(z^-)$:
$$
  \omega|_{D^-} = dy^- \wedge dx^-, \quad
  \omega|_{D^+} = dy^+ \wedge dx^+.
$$
If $F(D^-)\cap D^+\ne\emptyset$, the map $F$ can be represented as
\begin{equation}
\label{eq:can}
  y^+ = y^+(x^-,y^-), \quad
  x^+ = x^+(x^-,y^-), \qquad
  dy_+\wedge dx_+ = dy_-\wedge dx_-.
\end{equation}
The coordinates can be chosen so that $\det (\partial y^+ /
\partial x^-) \ne 0$. Then $F:D^-\to D^+$ is locally determined by
a generating function $L(x^-,x^+)$ on an open set $U^-\times U^+$:
$$
  y^- = \partial L/\partial x^-, \quad
  y^+ = - \partial L/\partial x^+.
$$

Since the construction is local, we obtain a collection
$\calL=\{L_\kappa\}_{\kappa\in J}$ of generating functions defined
on open sets $U_\kappa^-\times U_\kappa^+$. More precisely, take a
collection of  symplectic charts in $P$, i.e.\ open sets
$\{D_k\}_{k\in I}$ and symplectic maps
$$
  \phi_k:D_k\to U_k\times\mR^m, \quad
  U_k\subset\mR^m, \qquad
  \phi_k^*(dy\wedge dx) = \omega|_{D_k}.
$$
Let $J$ be the set of $\kappa=(\kappa_-,\kappa_+)\in I^2$ such
that $F(D_{\kappa_-})\cap D_{\kappa_+}\ne \emptyset$. For any
$\kappa\in J$, the map $F:D_{\kappa_-}\cap F^{-1}(D_{\kappa_+})\to
F(D_{\kappa_-})\cap D_{\kappa_+}$ satisfies
$$
F(x^-,y^-)=(x^+,y^+)\quad\Leftrightarrow\quad
y^+\,dx^+-y^-\,dx^-=dL_\kappa
$$
for some function $L_\kappa$. Changing coordinates if needed, we
can locally express $L_\kappa$ as a smooth function
$L_\kappa(x^-,x^+)$ on a subset in $U_{\kappa_-}\times
U_{\kappa_+}$. We obtain a DLS $\calL = \{L_\kappa\}_{\kappa\in
J}$. The corresponding graph has the set of vertices $J$, and there is an edge from $\kappa$ to $\kappa'\in J$ if $\kappa_+ = \kappa'_-$.

An orbit $\bdz=(z_i)$ of $F$ such that $z_i\in D_{k_{i}}$ for all $i$ corresponds to a trajectory $(\bdk,\bdx)$ of the DLS. Here $x_i\in U_{k_i}$ and $\kappa_i=(k_i,k_{i+1})$. Any trajectory of $\calL$ gives a unique trajectory of $F$.

If the charts $\{D_k\}_{k\in I}$ cover $P$ then any trajectory of $F$ corresponds to a trajectory (may be, not unique) of the DLS $\calL$.
However, in the
majority of examples below $\cup_{k\in I} D_k \ne P$.

The map $F:P\to P$ determines an ambient dynamical system for the
DLS $\calL$. A concrete example is given in Section
\ref{sec:Smale}.

\begin{rem} There are other  ways to represent trajectories of a
symplectic map $F$ by a DLS. If equations (\ref{eq:can}) can be
solved for $y^-(x^-,y^+)$, then $F:D^-\to D^+$ can be locally
represented by a generating function $S(x^-,y^+)$:
$$
  F(x^-,y^-)=(x^+,y^+)\quad\Leftrightarrow\quad
  d S (x^-,y^+) = y^-\,dx^-+x^+\,dy^+.
$$
We obtain a collection of functions $S_\kappa$ which generate
symplectic maps $F:D_\kappa^-\to D_\kappa^+$ of open sets in $P$.
Define the discrete Lagrangian by
$$
    L_{\kappa}(z^-,z^+)
  = \langle x^-,y^- \rangle - S_\kappa(x^-,y^+),\quad
    z^\pm
  = (x^\pm,y^\pm).
 $$
Then orbits $\bdz = (z_i)$ of $F$  are critical points of the
Poincar\'e discrete action functional
\begin{equation}
\label{eq:calA}
  \calA_\bdk(\bdz) = \sum L_{\kappa_i}(z_i,z_{i+1}) .
\end{equation}
Thus  $\bdz$ is a trajectory of a DLS with $2m$ degrees of
freedom. This representation is standard in symplectic topology,
see \cite{Arnold:matmex,McDuff}.
\end{rem}

\subsection{DLS generated by a Lagrangian flow}
\label{ssec:DLS_flow}

DLS naturally appear in Lagrangian systems with continuous time (CLS).
Consider a CLS with the Lagrangian $L(q,\dot q,t)$ on $TM\times
\mR$, where $M$ is a smooth manifold (the configurational space).
We assume that $L$ satisfies the Legendre condition: the Legendre
transform $\dot q\to p=\partial_{\dot q}L$ is a diffeomorphism.
Then $L$ defines a Lagrangian flow on $TM$ or Hamiltonian flow
$\phi^t$ with Hamiltonian $H(q,p,t)=\langle p,\dot q\rangle-L$ on
$T^*M$.

If the Lagrangian is $T$-periodic in time, then dynamics is described by the  monodromy (Poincar\'e) map $\phi^T:T^*M\to T^*M$. This map is symplectic and can be represented by a DLS. For Lagrangian systems this can be made more explicit.

Define Hamilton's action function as the action of a trajectory
$\gamma:[0,T]\to M$ joining two points $x_-,x_+\in M$:
$$
  S(x_-,x_+)=\int_0^T L(\gamma(t),\dot\gamma(t),t)\,dt.
$$
If the end points $x_\pm$ are nonconjugate along $\gamma$, then
$S$ is locally well defined and smooth. By Hamilton's first
variation formula, the initial and final momenta $p_-=p(0)$ and
$p_+=p(T)$ satisfy
\begin{equation}
\label{eq:dS} p_+\,dx_+-p_-\,dx_-=dS(x_-,x_+).
\end{equation}
Thus $S$ is the local generating function of the  monodromy map
$(q_-,p_-)\mapsto (q_+,p_+)$.

In general there can be several  trajectories joining $x_-$ and $x_+$, so the
function $S$ is multivalued. Thus we have a collection of generating functions $S_\kappa$ defined on open sets $U_\kappa^-\times U_\kappa^+\subset M^2$. Let $\gamma:\mR\to M$ be a trajectory of the CLS and  $x_i=\gamma(t_i)$, $t_i=Ti$. If the points $x_i$ and $x_{i+1}$ are nonconjugate along $\gamma$, then
$\bdx=(x_i)$ is a trajectory of the DLS $\calL=\{S_\kappa\}$.
Under the Legendre condition, also the converse is true: a
trajectory $\bdx$ of the DLS   gives a trajectory $\gamma$ of the
CLS with the momentum satisfying
$$
p(t_i+0)=
-\partial_{x_i}L_{\kappa_i}(x_i,x_{i+1})=\partial_{x_i}L_{\kappa_{i-1}}(x_{i-1},x_i)=p(t_i-0).
$$
Indeed, then
the momentum determines the velocity and so $\Delta p(t_i)=0$
implies $\Delta \dot\gamma(t_i)=0$. Hence $\gamma$ is smooth at
$t_i$ and so it is a trajectory of the CLS. Moreover the discrete
action functional corresponds to the action functional:
$\int_\gamma L\,dt=A_\bdk(\bdx)$.

\begin{rem}  This construction can be generalized. Take hypersurfaces
$\{\Sigma_k\}_{k\in I}$ in $M\times \mR$. For
$\kappa=(\kappa_-,\kappa_+)\in I^2$  define the discrete
Lagrangian $L_\kappa(x_-,x_+)$, $x_\pm=(q_\pm,t_\pm)\in
\Sigma_{\kappa_\pm}$ as the action of a
trajectory $\gamma:[t_-,t_+]\to M$ joining a  pair of
points $q_-,q_+$. In general there can be several such
trajectories so we obtain a collection of Lagrangians defined on
open sets  in $\Sigma_{\kappa_-}\times \Sigma_{\kappa_+}$. We use
this construction in Section  \ref{ssec:ail_cls_slow}.
\end{rem}

For autonomous systems, any $T>0$ is a period and can be used to define the  monodromy map $\phi^T$. Then $\phi^T$ has the energy integral $H$ and the
symmetry group defined by the  Hamiltonian flow $\phi^t$. The
corresponding DLS will have the energy integral and a (local)
symmetry group. Hence this DLS will never be anti-integrable.

We can avoid this difficulty by replacing the monodromy map
$\phi^T$ by the Poincar\'e map. Fix energy and consider a local cross-section $N$ in the energy level $H=E$. We can do this by choosing a hypersurface $\Sigma\subset M$ (in applications $\Sigma=\cup\Sigma_k$ may be a union of several hypersurfaces) and setting
\begin{equation} \label{eq:N}
N=\{(q,p):q\in\Sigma,
H(q,p)=E\}.
\end{equation}
 This is a symplectic manifold and the Poincar\'e map $F:U\to N$ is a symplectic map of an open set $U\subset N$.

We represent $F$ by a DLS with $m-1$ degrees of freedom as
follows. Define the discrete Lagrangian as the Maupertuis action
$$
  S(x_-,x_+)=\int_\gamma p\,dq,\qquad p=\partial_{\dot q}L,
$$
of a trajectory $\gamma:[0,T]\to M$ with energy $H=E$ joining given points $x_\pm\in M$. Here $T=T(\gamma)>0$ is arbitrary. The action $S$ is smooth provided $x_\pm$ are nonconjugate along $\gamma$ for the Maupertuis action functional. The action satisfies (\ref{eq:dS}) but the Hamilton--Jacobi equation $H(x_+,\partial_{x_+}S)=E$ implies that the twist condition always fails.

A trajectory $\gamma$ not always exists, and there can be several
of them, so we obtain a collection of discrete Lagrangians
$S_\kappa$ on open sets $U_\kappa^-\times U_\kappa^+\subset M^2$.
Take a hypersurface $\Sigma\subset M$ and set
$\Sigma_\kappa^\pm=\Sigma \cap U_\kappa^\pm$. Then
$L_\kappa=S_\kappa|_{\Sigma_\kappa^-\times \Sigma_\kappa^+}$ is a
local generating function of the Poincar\'e map of the cross
section (\ref{eq:N}). Thus $\calL=\{L_\kappa\}$ is a  DLS
describing trajectories of the continuous Lagrangian system on the
energy level $H=E$. If $\gamma:\mR\to M$ is a trajectory with
$H=E$ such that $\gamma(t_i)\in\Sigma$ and $x_i,x_{i+1}$ are
nonconjugate (for fixed energy) along $\gamma$, then $\bdx=(x_i)$
is a trajectory of the DLS.

Note that not all trajectories of the DLS correspond to trajectories
of the continuous Lagrangian system: some correspond to billiard
trajectories with elastic reflection from $\Sigma$. Indeed, for
given $x\in\Sigma$, the equations
$$
  H(x,p)=H(x,p'),\quad \Delta p=p'-p_-\perp T_x\Sigma ,
$$
do not imply $p'=p$. If $H$ is convex in $p$, there are 2
solutions for $p'$: one describing elastic reflection from
$\Sigma$ and the other ($p'=p$) a smooth trajectory.

However, in the local problems discussed in this paper this
difficulty does not appear: locally equations above have unique
solution.

\section{Anti-integrable limit in DLS}

\subsection{Main theorem}
\label{subsec:main}

The original idea of AI limit \cite{aubry0} was extended in
various directions. The case of multidimensional standard map was
considered in \cite{mackmei} and other papers, see the references in \cite{East-Meiss}. General equivariant case is
presented in \cite{TZ}. Billiard systems with small convex
scatterers are studied in \cite{Chen2004}. A particle with a small
mass in a potential force field as a continuous Lagrangian system
near AI limit is treated in \cite{Bol-Mac:anti}. Other examples
and references can be found in Section \ref{sec:ail_cls}. In this
section we present a general approach to AI limit in DLS which
covers essentially all known examples.

Consider a  DLS $\calL=\{L_\kappa\}_{\kappa\in J}$ as defined in Section \ref{ssec:DLS_multi}. In the AI limit the Lagrangians  are small perturbations of twistless Lagrangians $L_\kappa^0$:
$$
  L_\kappa(x,y) = L_\kappa^0(x,y)+u_\kappa(x,y), \qquad
  \kappa\in J,
$$
where the function $u_\kappa$ is small in the $C^2$ norm and
$L_\kappa^0$ has identically zero twist: $B_\kappa^0(x,y)\equiv
0$. Then
$$
  L_\kappa^0(x,y) = V_\kappa^-(x) + V_\kappa^+(y),
$$
where the functions $V_\kappa^\pm$  are defined on the sets
$U_\kappa^\pm$  respectively.

Define an oriented  graph $\Gamma$ with the set of vertices $J$
and the set of edges $E$ as follows. We join the vertices
$\kappa,\kappa'\in J$ with an edge $\gamma=(\kappa,\kappa')$ if
$W_{\gamma}=U_{\kappa}^+\cap U_{\kappa'}^-\ne\emptyset$ and the
function $\Psi_{\gamma}=V_{\kappa}^+ + V_{\kappa'}^-$ on
$W_\gamma$ has a nondegenerate critical point.  Then to each edge
$\gamma=(\kappa,\kappa')\in E$ joining vertices $\kappa,\kappa'\in
J$ there corresponds a nondegenerate critical point  $a_\gamma$ of
the function $\Psi_{\gamma}$.

\begin{rem}
In many applications for given $\kappa,\kappa'$ the critical point
$a_\gamma$ will be unique. If there are several nondegenerate
critical points we  join
$\kappa,\kappa'$ by several edges $\gamma$, and to each of them
there corresponds the nondegenerate critical point  $a_\gamma$. In
this case $\Gamma$ is not a simple graph: there  are several edges $\gamma$
joining two vertices $\kappa,\kappa'$. Theorem \ref{thm:anti} and its proof   still work, but a
path in the graph will mean a sequence of edges, not vertices.
\end{rem}

According to the traditional definition, a path in the graph $\Gamma$ is a sequence of vertices $\bdk=(\kappa_i)_{i\in\mZ}$. It defines a sequence of edges $(\gamma_i)_{i\in\mZ}$, where $\gamma_i$ joins the vertex $\kappa_{i-1}$ with $\kappa_i$. The corresponding sequence of critical points
$$
  \bda(\bdk) = (a_{\gamma_i})_{i\in\mZ}.
$$
is a trajectory of the twistless Lagrangian system
$\calL^0=\{L_\kappa^0\}$ i.e.\ a critical point of the uncoupled
action functional:
\begin{equation}
\label{eq:split}
  A_\bdk^0(\bdx)=\sum_{k\in\mZ}\Psi_{\gamma_i}(x_i).
\end{equation}

The set $\Pi\subset J^\mZ$ of all paths is invariant under the shift $\calT:\Pi\to\Pi$, $(\kappa_i)\to (\kappa_{i+1})$. For simplicity suppose that the graph $\Gamma$ is finite. The general case is discussed later in Section \ref{sec:equi}. Under this assumption the set $\Pi$ of all paths is a compact shift invariant subset in the Cantor set $J^\mZ$. The dynamical system  $\calT:\Pi\to\Pi$ is called a topological Markov chain or a subshift of finite type \cite{KH}. When the graph is complete, i.e. $\Pi=J^\mZ$, this is the Bernoulli shift with $\# J$ symbols.

\begin{theo}
\label{thm:anti} Suppose $\#\Gamma < \infty$. Then there exist
$C,\eps_0,\sigma>0$ such that if
\begin{equation}
\label{eps0}
  \eps=\max_{\kappa}\|u_\kappa\|_{C^2}<\eps_0,
\end{equation}
then the following is true.

{\bf (a)} For any path $\bdk=(\kappa_i)$ in $\Gamma$, there exists
a unique orbit $\bdx(\bdk)=(x_i)$ of the DLS $\calL$ which shadows
the sequence $\bda=\bda(\bdk)$:
$$
  \rho(\bdx,\bda) = \sup_{i}\, \dist(x_i,a_{\gamma_i}) < \sigma.
$$
Moreover $\rho(\bdx,\bda)< C\eps$.

{\bf (b)} The orbit $\bdx$ is hyperbolic with Lyapunov exponent
$\lambda\ge C|\log\eps|$.

{\bf (c)} The map $\bdk\mapsto\bdx$ is continuous, so
$\Lambda=\{(\bdk,\bdx(\bdk)):\bdk\in\Pi\}$ is a compact set in
$J^\mZ\times (\mR^m)^\mZ$.
\end{theo}

\begin{rem}
The assumption that $\Psi_\gamma$ has  a nondegenerate critical
point can be weakened. For example, for the existence of a
shadowing trajectory it is sufficient that $\Psi_\gamma$ has a
compact set of minimum points in $W_\gamma$. However, then the
shadowing orbit may be not unique or hyperbolic. Further generalizations
can be obtained by variational methods of Mather, see  \cite{Mather}.
\end{rem}

If the graph $\Gamma$ is branched,  i.e.\ there are 2 cycles
through the same vertex, then $\calT:\Lambda\to\Lambda$ has
positive topological entropy. If the DLS admits an ambient system
$F:P\to P$, we obtain a homeomorphism $\Pi\to\Lambda\to\calK$ onto
a compact hyperbolic set $\calK\subset P$, and $F:\calK\to\calK$
is conjugate to the topological Markov chain $\calT:\Pi\to\Pi$. In
this case we have the following

\begin{cor}
If the graph $\Gamma$ is branched, then  the set $\calK$ has
positive topological entropy satisfying (\ref{top_ent}).
\end{cor}

\begin{rem} For our purposes only the local existence of the ambient
system in a  neighborhood of the hyperbolic set is needed. Thus
only the local twist condition is required:
$$
  \det B_\kappa(a_\gamma,a_{\gamma'})\ne 0
$$
for any two edges such that $\gamma$ ends at $\kappa$ and
$\gamma'$ starts at $\kappa$.
\end{rem}

Theorem \ref{thm:anti} does not require the twist condition. Hence, we have to generalize the usual definition of hyperbolicity. Indeed, without the twist condition the present state of a trajectory does not determine the past and the future so the usual definition of the stable (contracting) and unstable (expanding) subspaces does not work. We reformulate the cone hyperbolicity criterion of V.M.~Alexeyev
\cite{KH}.

Linearization of equation (\ref{eq:crit}) at $\bdx$ yields the
variational equation of the trajectory $(\bdk,\bdx)$:
\begin{equation}
\label{linearization}
  G_{i-} u_{i-1} + G_i u_i + G_{i+} u_{i+1} = 0,
\end{equation}
where $G_i$ and $G_{i\pm}$ are linear operators. We define the
cones
\begin{eqnarray}
\label{coneH}
     H_i
 &=& \big\{ (u_i,u_{i+1})  :
              \|u_i\| \le \alpha_H \|u_{i+1}\| \big\}, \\
\label{coneV}
     V_i
 &=& \big\{ (u_i,u_{i+1})  :
              \|u_{i+1}\| \le \alpha_V \|u_i\| \big\} .
\end{eqnarray}
We say that $(\bdk,\bdx)$ is hyperbolic if there exists $\mu > 1$
such that for any $i\in\mZ$ and for any $u_{i-1},u_i,u_{i+1}$
satisfying (\ref{linearization})
\begin{eqnarray}
\label{hypH}
  (u_{i-1},u_i)\in H_{i-1} \quad
  \mbox{implies} \quad
  (u_i,u_{i+1})\in H_i \;
  \mbox{and} \;
  \|(u_i,u_{i+1})\| \ge \mu \|(u_{i-1},u_i)\|, \\
\label{hypV}
  (u_i,u_{i+1})\in V_i \quad
  \mbox{implies} \quad
  (u_{i-1},u_i)\in V_{i-1} \;
  \mbox{and} \;
  \|(u_{i-1},u_i)\| \ge \mu \|(u_i,u_{i+1})\|.
\end{eqnarray}
The cone definition has this form only for a good choice of the
metric. But this is not important for us because we use it only as
a sufficient condition for hyperbolicity.

If $(u_i,u_{i+1})\in H_i$, then $\|(u_j,u_{j+1})\|\ge \mu^{j-i}
\|(u_i,u_{i+1})\|$ for $j\ge i$, so the Lyapunov exponent is at
least $\log\mu$.
\medskip

We say that a compact $\calT$-invariant set  $\Lambda$ of
trajectories $(\bdk,\bdx)$  is hyperbolic if every trajectory
$(\bdk,\bdx)$ in $\Lambda$ satisfies conditions above with the
same constants. If a DLS admits an ambient system, then the
corresponding
compact invariant set $\calK\subset P$  is hyperbolic in the traditional \cite{KH} sense.

\begin{rem}   Another condition for hyperbolicity for $(\bdk,\bdx)$ is
that the Hessian operator $A''_\bdk(\bdx):l_\infty(\mR^m)\to
l_\infty(\mR^m)$ defined by the right hand side of the variational
equation  has bounded inverse in the $l_\infty$ norm. The
equivalence of these definitions was proved in
\cite{Aubry-MacKay}.
\end{rem}

\subsection{Proof of Theorem \protect\ref{thm:anti}}

The proof of Theorem \ref{thm:anti} is almost the same as for the
Standard map. By the implicit function theorem for any edge
$\gamma=(\kappa,\kappa') \in E$ joining $\kappa,\kappa'\in J$ the
map $g_\gamma=D\Psi_{\gamma}$ is a diffeomorphism from a
neighborhood $W_\gamma$ of $a_\gamma$ to a neighborhood of the
origin in $\mR^m$.   We may assume that
\begin{equation}
\label{Ugamma}
    W_\gamma
  = \{x\in \mR^m :
           |x - a_\gamma| < r \}.
\end{equation}
Here $r>0$ can be taken independent of $\gamma$ because the graph
$\Gamma$ is finite. Let $\phi_\gamma$ denote the inverse map:
$\phi_\gamma = g_\gamma^{-1}$. Then $\phi_\gamma(0) = a_\gamma$.

Let $\bdk=(\kappa_i)$ be a path in the graph $\Gamma$ and
$(\gamma_i)$ the corresponding sequence of edges.   Consider the
metric space $(X,\rho)$, where $X = \prod_{i\in\mZ} W_{\gamma_i}$
and for any $\bdx,\bdx'\in X$
$$
    \rho(\bdx,\bdx')
  = \sup_{i\in\mZ} |x_i - x'_i|.
$$
A trajectory $\bdx$ is a critical point of the functional
$A_{\bdk}$ in a neighborhood of $\bda = \bda(\bdk)$ in $X$ iff
\begin{equation}
\label{phi}
    x_i
  = \phi_{\gamma_i}
             \big( \partial_{x_i}
                 \big( u_{\kappa_{i-1}}(x_{i-1},x_{i})
                      + u_{\kappa_i}(x_i,x_{i+1}) \big) \big).
\end{equation}
It remains to apply contraction principle in a neighborhood of the
point $\bda$ in $X$. Let $\calB\subset X$ be the ball
$$
  \calB = \{\bdx\in X : \rho(\bdx,\bda) \le \sigma\},\qquad \sigma<r.
$$
Consider the operator $\Phi : \calB_\sigma\to X$,
$$
  \bdx\mapsto\bdy = \Phi(\bdx), \quad
    y_i
  = \phi_{\gamma_i}\big( \partial_{x_i}
                 \big( u_{\kappa_{i-1}}(x_{i-1},x_{i})
                      + u_{\kappa_i}(x_i,x_{i+1}) \big) \big).
$$
Any fixed point of $\Phi$ is a trajectory of the DLS $\calL$.

As in Section \ref{sec:standard} it is easy to show that if $\eps$
is sufficiently small, $\Phi(\calB_\sigma) \subset \calB_\sigma$
and $\Phi$ is contracting. Indeed,
\begin{eqnarray*}
     \rho(\bdy,\bda)
 &=& \sup_{i\in\mZ}
      \big| \phi_{\gamma_i}\big( \partial_{x_i}
                 \big( u_{\kappa_{i-1}}(x_{i-1},x_{i})
                      + u_{\kappa_i}(x_i,x_{i+1}) \big) \big)
           - \phi_{\gamma_i}(0) \big| \\
 &\le& \lambda \sup_{i\in\mZ}
      \big| \partial_{x_i} \big( u_{\kappa_{i-1}}(x_{i-1},x_{i})
                      + u_{\kappa_i}(x_i,x_{i+1}) \big) \big|
  \le 2\lambda\eps.
\end{eqnarray*}
Hence $\bdy\in\calB_\sigma$ if we set $\sigma = 2\lambda\eps$.

We put $\bdy = \Phi(\bdx)$, $\bdy' = \Phi(\bdx')$,
$\bdx,\bdx'\in\calB_\sigma$. Then
\begin{eqnarray*}
     \rho(\bdy,\bdy')
 &=& \sup_{i\in\mZ}
      \big| \phi_{\gamma_i}\big( \partial_{x_i}
                 \big( u_{\kappa_{i-1}}(x_{i-1},x_{i})
                      + u_{\kappa_i}(x_i,x_{i+1}) \big) \big) \\
 && \qquad
          - \phi_{\gamma_i}\big( \partial_{x'_i}
                 \big( u_{\kappa_{i-1}}(x'_{i-1},x'_{i})
                      + u_{\kappa_i}(x'_i,x'_{i+1}) \big) \big) \big| \\
 &\le& \lambda \sup_{i\in\mZ}
      \big| \partial_{x_i} \big( u_{\kappa_{i-1}}(x_{i-1},x_{i})
                      + u_{\kappa_i}(x_i,x_{i+1}) \big)
         - \partial_{x'_i} \big( u_{\kappa_{i-1}}(x'_{i-1},x'_{i})
                      + u_{\kappa_i}(x'_i,x'_{i+1}) \big) \big| \\
 &\le& 2\lambda\eps \rho(\bdx,\bdx').
\end{eqnarray*}
Hence, $\Phi$ is contracting for $\eps < 1 / (2\lambda)$.
\smallskip

Now we prove assertion {\bf (b)}. We take in
(\ref{coneH})--(\ref{coneV}) $\alpha_H=\alpha_V=1/2$. Let $\bdx$
be the orbit corresponding to $\bdk\in\Pi$ by assertion {\bf (a)}.
Differentiating (\ref{phi}), we have:
\begin{equation}
\label{eq:PQR}
  u_i = P_i u_{i-1} + Q_i u_i + R_i u_{i-1} .
\end{equation}
This equation is equivalent to (\ref{linearization}).

Norms of the linear operators $P_i,Q_i,R_i$ are small:
$$
  \|P_i\| \le \lambda\eps, \quad
  \|Q_i\| \le \lambda\eps, \quad
  \|R_i\| \le \lambda\eps .
$$
The equations
$$
  R_i u_{i+1} = (I-Q_i) u_i - P_i u_{i-1}, \quad
  (u_{i-1}, u_i) \in H_{i-1}
$$
imply
$$
      \frac{|u_i|}{|u_{i+1}|}
  \le \frac{|(I-Q_i) u_i - P_i u_{i-1}|}
           {(1 - \|Q_i\| - \|P_i\|/2)\cdot |u_{i+1}|}
  \le \frac{|R_i u_{i+1}|}{(1 - 3\eps\lambda/2)\cdot |u_{i+1}|}
  \le \frac{\eps\lambda}{ 1 - 3\eps\lambda/2 }.
$$
Therefore
$$
  (u_i, u_{i+1}) \in H_i \quad
  \mbox{provided}\quad
  \frac{\eps\lambda}{ 1 - 3\eps\lambda/2 } < \frac12.
$$
This implies the first assertion (\ref{hypH}). We also have:
$$
      \|(u_{i-1}, u_i)\|^2
  \le |u_{i-1}|^2 + \frac{|(I-Q_i) u_i - P_i u_{i-1}|^2}
                         {(1 - \|Q\| - \|P\|/2)^2}
  \le \frac{|u_i|^2}4 + \frac{\eps^2\lambda^2 |u_{i+1}|^2}
                             {(1 - 3\eps\lambda/2)^2}
  \le \frac14 \|(u_i, u_{i+1})\|^2
$$
provided $\frac{\eps\lambda}{ 1 - 3\eps\lambda } < \frac12$. Hence
we have the second assertion (\ref{hypH}) with $\mu = 2$.
Assertions (\ref{hypV}) can be checked analogously.

Now we prove {\bf (c)}. If two codes $\bdk$ and $\bdk'$ are close
in the product topology in $J^\mZ$, they have a long identical
segment: $\kappa_i=\kappa'_i$ for $|i|\le n$. Then the corresponding trajectories $\bdx$ and $\bdx'$ satisfy
$$
  |x_i-x'_i|\le C\alpha^{|i|-n}, \qquad
  0<\alpha<1, \quad
  |i|\le n.
$$
This follows from (\ref{phi}) as in the proof of Lemma \ref{lem:a'=a''}. If $n$ is large, $\bdx$ is close to $\bdx'$ in the product topology. Hence
the map $\bdk\mapsto\bdx$ is continuous. \qed
\smallskip

Another  way to check hyperbolicity is to show that the Hessian
$A_\bdk''(\bdx)$ has bounded inverse in the $l^\infty(\mR^m)$
norm. Since it is 3-diagonal with invertible diagonal and small
off diagonal terms, this is almost evident. We can also describe
the stable and unstable subspaces of the hyperbolic trajectory.

\begin{prop}\label{prop:hyp}
For any $u\in\mR^m$ there exists a unique trajectory
$\bdu=(u_j)_{j\ge 0}$ of the variational system such that $u_0=u$
and $\|u_j\|$ is bounded as $j\to\infty$. The trajectory $\bdu$
exponentially tends to zero: $\|u_j\|\le \mu^{-j}\|u\|$, $\mu>1$. Thus
$\bdu$ belongs to the stable subspace of the trajectory
$(\bdk,\bdx)$. Similarly for the unstable subspace.
\end{prop}

The proof of Proposition \ref{prop:hyp} is obtained  by the same
contracting mapping argument applied to the map
$$
  (u_i)_{i>  0} \mapsto (v_i)_{i>  0}, \qquad
  v_i = P_i u_{i-1} + Q_i u_i + R_i u_{i-1},\quad u_0=u,
$$
(see (\ref{eq:PQR})) on the space of bounded sequences
$(u_j)_{j>0}$.

\subsection{$G$-equivariant DLS}
\label{sec:equi}

In many applications we use anti-integrable limit to construct
trajectories going to infinity. Then we have to consider infinite
graphs $\Gamma$. In the case of an infinite graph $\Gamma$ we need
a certain uniformity.
\medskip

\noindent {\bf U: Uniform anti-integrability}. There exist positive constants $r$ and $\lambda$ such that for any $\gamma \in E$,

(a) $g_\gamma:W_\gamma\to B_\gamma\subset\mR^m$ is a
diffeomorphism where $B_\gamma$ is a neighborhood of the origin
and $W_\gamma$ is determined by (\ref{Ugamma}) with $r>0$
independent of $\gamma$.

(b) the map $\phi_{\gamma} = {g_\gamma}^{-1}:B_\gamma\to W_\gamma$
is Lipschits with Lipschits constant $\lambda$.

(c) $\eps=\sup_{\kappa\in J}\|u_\kappa\|_{C^2}$ is finite and
sufficiently small.

\medskip

Obviously, condition {\bf U} holds if $\Gamma$ is finite.

\begin{theo}
\label{thm:anti_inf} Theorem \ref{thm:anti} remains true for an
infinite graph $\Gamma$ if the DLS is uniformly anti-integrable.
\end{theo}

The proof coincides with the proof of Theorem \ref{thm:anti}. \qed

Condition {\bf U} is a restrictive assumption. In this section we
present a class of discrete Lagrangian systems with $\# E =
\infty$ where {\bf U} requires a finitely many conditions to hold.

Suppose that a discrete group $G$ acts on $M$.\footnote {i.e.\
there is a homomorphism of the group $G$ to the group of
diffeomorphisms of $M$.} We assume that the action is discrete:
any point $x\in M$ has a neighborhood $U$ such that the sets
$g(U)$, $g\in G$, do not intersect: $g'(U) \cap g''(U) =
\emptyset$ for $g' \ne g''$.   In this case the quotient space
$\widetilde M=M/G$ is a smooth manifold and $\pi:M\to \widetilde
M$ is a covering. The action of $G$ on $M$ generates the diagonal
action of $G$ on the product $M\times M$: for any pair $x,y\in M$
and $g\in G$ we have $g(x,y) = (g(x),g(y))$.

Let $L: M\times M\to\mR$ be invariant with respect to the action
of $G$:
$$
  L(x,y) = L(g(x,y))\qquad
  \mbox{for all} \quad  x,y\in M \quad
  \mbox{and } g\in G.
$$
The corresponding DLS is called $G$-equivariant. The ambient
system here may be defined as a map $F:P\to P$ of $P=(M\times
M)/G$. In the case of the multidimensional standard map with the
Lagrangian (\ref{eq:standard}) we have $M=\mR^m$ and $G\cong\mZ^m$
is the group of shifts preserving the potential $V$.

As another example suppose that $M$ is a Riemannian manifold and
the group $G$ acts on $M$ by isometries. Let
$\mbox{dist}(\cdot\,,\cdot)$ be the distance induced by the
Riemannian metric. Then the Lagrangian $L(x,y) = \mbox{dist}^2
(x,y)$ is a smooth function for any pair of sufficiently close
points $x,y$. It is invariant with respect to the diagonal action
of $G$.

The function $z= z(x,y)$ (see (\ref{-0+})) in this example has a
simple geometric meaning. Suppose there exists a unique shortest
geodesic $\gamma:[0,1]\to M$ joining $x$ with $z$ (for example,
the points $x,z$ are close to one another). Then $y=\gamma(1/2)$
is situated on the same geodesic, $y$ lies between $x$ and $z$ and
$$
    \mbox{dist}\,(x,y)
  = \mbox{dist}\,(y,z)
  = \frac12 \mbox{dist}\,(x,z).
$$

For smooth $G$-invariant functions $V : M\to\mR$ and $u: M\times
M\to \mR$ consider the DLS, determined by the discrete Lagrangian
$L$
$$
    L(x,y) = u(x,y) + V(y).
$$
We define the AI limit in such a system as the limit $\|u\|_{C^2}\to 0$.

Suppose that the configurational space $\widetilde M=M/G$ is
compact. We fix  a Riemannian metric on $\widetilde M$ which is
pulled back up to a $G$-invariant  metric on $M$. Let $\dist$ be
the corresponding distance on $M$.

Take    a finite set of nondegenerate
critical points of $ V$ on $\tilde M$, and  let $I$ the corresponding
$G$-invariant set of critical points of $V$ on $M$. The latter is a
finite union of the orbits
$$
 \{g(k):g\in G\}, \qquad k\in I,
$$
of the group action.

Taking a large constant $N$ we define the graph $\Gamma$ with the set of vertices $J$ and the set of edges $E$:
$$
    J
  = \big\{ \kappa = (\kappa_-,\kappa_+)\in I^2 :
             \dist(\kappa_-,\kappa_+) < N \big\}, \quad
    E
  = \big\{\gamma = (\kappa,\kappa')\in J^2 :
                     \kappa_+ = \kappa'_- \big\}.
$$

To make notations consistent with  the notation  of Section \ref{subsec:main}, we introduce the DLS $\calL$ as follows.
Take small $G$-invariant $\rho$-neighborhood of $U=\cup_{k\in I}U_k$ of the set $I$ in $M$. Then
for $\kappa\in J$ set
$$
L_\kappa=L|_{U_{\kappa_-}\times U_{\kappa_+}},\quad  V_\kappa^- =
0, \quad
   V_\kappa^+ = V|_{U_{\kappa_+}}, \quad
   u_\kappa = u|_{U_{\kappa_-}\times U_{\kappa_+}}
$$

If the group $G$ is infinite, the  graph $\Gamma$ is infinite, but
due to $G$-invariance of $L$, condition   ${\bf U}$ follows from
finiteness of the set $\Cr$ and boundedness of $N$. One can notice
that $G$ acts on $\Gamma$ as well and $\tilde\Gamma=\Gamma/G$ is a
finite graph. Then we obtain:

\begin{prop}
The map $F$ has a hyperbolic set $K$ such that $F|_K$ is conjugate to the
topological Markov chain determined by the graph $\Gamma$.
\end{prop}

In fact we need only the local existence of the ambient system
near the hyperbolic set. It will exist provided $\det
B(\kappa_-,\kappa_+)\ne 0$ for all $(\kappa_-,\kappa_+)\in J$.

Note that if $M$ is not compact, $M^2 / G$ is also not compact,
but the set $\calK$ is compact because all its points are located
on a distance less than $N+\sigma$ from the compact set
$\{(x,x)\in M^2\} / G$.

\section{Examples: systems with discrete time}
\label{sec:exa_discr}

\subsection{Light particle and periodic kiñks}

Consider a particle with a small mass $\eps^2$ which moves in the
space $\mR^m$ in the force field with potential
\begin{equation}
\label{potentsing}
  \calV(x,t) = \frac{1}{2\pi} V(x) \delta(t),
\end{equation}
where the function $V$ is smooth on $\mR^m$ and $\delta$ is the periodic
$\delta$-function:
$$
    \delta(t)
  = \left\{\begin{array}{cl} \infty, &\quad t\in 2\pi\mZ,\\
                                0  , &\quad t\in \mR\setminus 2\pi\mZ
           \end{array} \right. , \qquad
   \int_{2\pi k-\sigma}^{2\pi k + \sigma} \delta(t) \, dt = 1\quad
   \mbox{for any $k\in\mZ$, $\sigma\in (0,\pi)$}.
$$

The Hamiltonian of the system has the form $H = \frac1{2\eps^2}
|p|^2 + \frac 1{2\pi} V(x) \delta(t)$, where $p =
(p_1,\ldots,p_m)$ is the momentum canonically conjugate to the
coordinates $x=(x_1,\dots,x_m)$. The Hamiltonian equations read
$$
  \dot p = - \frac 1{2\pi} \frac{\partial V}{\partial x}(x)
\delta(t),\quad
  \dot x = \frac{p}{\eps^2}.
$$
Therefore $p$ gets increments $-(2\pi)^{-1} \partial V / \partial
x$ at the time moments $2\pi l$, $l\in\mZ$. During the remaining
time the particle is free.

For any integer $l$ we put $x(2\pi l - 0) = x_l$, $p(2\pi l - 0) =
p_l$. Then
\begin{eqnarray*}
  \left(\begin{array}{l} x_l \\ p_l \end{array}\right)
  &\mapsto &
  \left(\begin{array}{l} x(2\pi l + 0) \\ p(2\pi l + 0) \end{array}
  \right)
 =
  \left(\begin{array}{l} x_l \\ p_l - \frac{1}{2\pi}
                                     \frac{\partial V}{\partial x} (x_l)
  \end{array}\right) \\
  &\mapsto &
  \left(\begin{array}{l} x_{l+1} \\ p_{l+1} \end{array}
  \right)
 =
  \left(\begin{array}{l} x_l + 2\pi\eps^{-2} p_{l+1}  \\
                         p_l - \frac{1}{2\pi}
                               \frac{\partial V}{\partial x} (x_l)
  \end{array}\right).
\end{eqnarray*}

The quantities $x_{l-1},x_l,x_{l+1}$ satisfy the equation $
x_{l+1} - 2x_l + x_{l-1}
  = \eps^{-2} \frac{\partial V}{\partial x}(x_l)$. We have a DLS
  with the Lagrangian which has the form (\ref{eq:standard}), where
$B=\eps^2 I$.
\begin{equation}
\label{light}
  L(x,y) = \eps^2\frac{|x-y|^2}{2}
              -   V (y).
\end{equation}
The corresponding discrete dynamical system is the
multidimensional standard map. By Theorem \ref{thm:anti} if $V$
has two nondegenerate critical points, then for small $\eps$ the
system has a chaotic hyperbolic set. If $V$ is $\mZ$-periodic on
$\mR^m$, then we have an $\mZ$-equivariant DLS corresponding to a
symplectic twist self-map of $\mT^m\times\mR^m$.

Note that for small $\eps$ the system remains close to the AI limit if the potential $V$ is a smooth periodic function close (as a distribution) to~(\ref{potentsing}).

If a light particle travels in a potential force field, where the
potential $V(x,t)$ does not satisfy (\ref{potentsing}), the theory
of the AI limit becomes technically more complicated.
We discuss these results in Section \ref{ssec:ail_cls_slow}.

\subsection{Billiard in a wide strip}

\begin{figure}[ht]
\begin{center}
\includegraphics{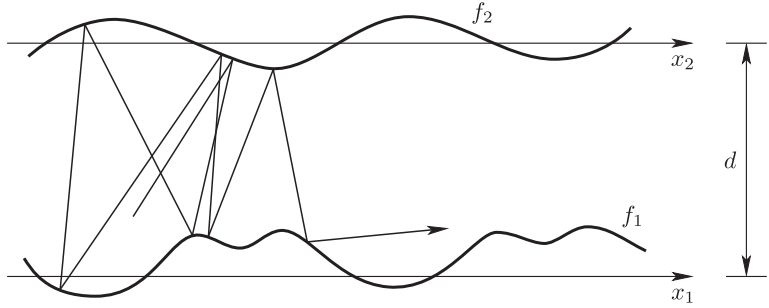}
\caption{ Billiard in a wide strip.} \label{fig:bill}
\end{center}
\end{figure}

Consider a plane billiard system in a wide strip bounded by graphs
of two $1$-periodic functions. In other words, we assume that a
particle moves in the domain
$$
  D = \{(x,y)\in \mR^2 : f_1(x) \le y \le f_2(x) + d \},
$$
where $f_1,f_2$ are $1$-periodic functions and the parameter $d$ is large (Fig. \ref{fig:bill}). The motion of the particle inside the domain
is assumed to be free. Reflections from the boundary are elastic.

This is a billiard system (see Section \ref{ssec:DLS}). Let $L$ be
the length of the line segment between two subsequent points of
the impact with the boundary. We will consider motions such that
the particle collides alternately with the upper and lower walls.

Let $x_1$ be the coordinate on the lower boundary and $x_2$ on the
upper one. The length of the corresponding line segment is
$$
  L(x_1,x_2) = \sqrt{ (x_2 - x_1)^2 + (d + f_2(x_2) - f_1(x_1))^2}.
$$
The Lagrangian is $\mZ$-invariant and for large $d$,
$$
    L(x_1,x_2)
  = d + f_2(x_2) - f_1(x_1) + \frac \eps{2} (x_2 - x_1)^2
      + O(\eps^2),\qquad
       d = \eps^{-1},
$$
it has the form (\ref{eq:standard}).
The action functional has the form
$$
     A(\bdx)
  = \sum_{i\in\mZ} \big( L(x_{2i-1},x_{2i}) + L(x_{2i+1},x_{2i}) \big).
$$

In the notation of  section \ref{ssec:DLS_multi} we have a DLS $\calL=\{L_1,L_2\}$ with the configuration space $M=\mR_1\cup\mR_2$ a union of two copies of $\mR$ and the Lagrangians $L_1(x_1,x_2)=L(x_1,x_2)$
and $L_2(x_2,x_1)=L(x_1,x_2)$. The action functional
$$
     A(\bdx)
  = \sum_{i\in\mZ} L_{\kappa_i}(x_i,x_{i+1})=  \sum_{i\in\mZ}(d + 2(-1)^i f_{\kappa_i}(x_i)+O(\eps)),  \qquad \kappa_i=i\bmod 2,
$$
has the uncoupled form (\ref{eq:split}). Thus the DLS is anti-integrable provided  both functions $f_1$ and $f_2$ have nondegenerate critical points.
Let $\Cr_1\subset\mR_1$  and $\Cr_2\subset\mR_2$ denote the $\mZ$-invariant sets of nondegenerate critical points of $f_1$ and $f_2$. Suppose that the sets $\Cr_1/\mZ$ and $\Cr_2/\mZ$ are finite and nonempty.

The corresponding graph $\Gamma$ has vertices of two types: $J = J_1\cup J_2$
\begin{eqnarray*}
     J_1
 &=& \big\{ (\kappa,x,y) \in \{1\}\times\Cr_1\times\Cr_2 : |x-y| < N
     \big\}, \\
     J_2
 &=& \big\{ (\kappa,x,y) \in \{2\}\times\Cr_2\times\Cr_1 : |x-y| < N
     \big\}.
\end{eqnarray*}
There is an edge from the vertex $(\kappa,x,y)$ and $(\kappa',x',y')$ if and only if $\kappa\ne\kappa'$ and $y = x'$.

As in the standard equivariant situation (Section \ref{sec:equi}) we obtain the existence of a hyperbolic set carrying dynamics conjugated to the dynamics in a topological Markov chain. The graph $\Gamma$ is infinite, but the condition of  uniform anti-integrability holds.

Multidimensional analog of this system is straightforward: it is
sufficient to say that $x\in\mR^m$ and the functions $f_1,f_2$ are
periodic in all components of the vector variable $x$.

\subsection{Billiard systems with small scatterers}
\label{ssec:small_scat}

Let $D\subset\mR^m$ be a domain with smooth boundary
$\Sigma=\partial D$. Following \cite{Chen2004}, consider a
billiard system in the domain $\Omega_\eps = D\setminus
(\cup_{j=1}^N D_j)$, where $D_j \subset D$  are small  subdomains
which play the role of scatterers. They are small in the following
sense. Each domain $D_j$ is associated with some point $a_j\in D$.
The boundary $\partial D_j$ is
$$
     a_j+\eps\Sigma_j
  =  \{ q\in\mR^m : q = a_j + \eps \phi_j(x), \; x\in S^{m-1} \},
$$
where the vector-functions $\phi_j:S^{m-1}\to\mR^m$ are smooth
embeddings which define smooth submanifolds
$\Sigma_j=\phi_j(S^{m-1})$. If $\eps$ is sufficiently small, then
$D_j\subset D$. When $\eps\to 0$, $\Omega_\eps$ degenerates to
$\Omega_0=D\setminus A$, where $A = \{a_1,\dots, a_N\}$.

\begin{figure}[ht]
\begin{center}
\includegraphics{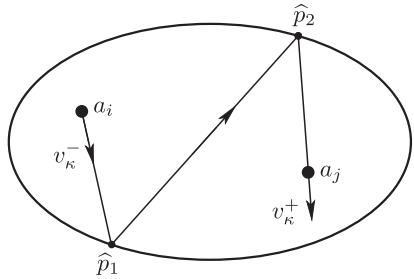}
\caption{An orbit of the billiard  in $\Omega_0$.}
\label{fig:bill_ss1}
\end{center}
\end{figure}

We plan to explain that for small $\eps$ this billiard system
generates a DLS with the discrete Lagrangian
$\calL=\{L_\kappa\}$, where the index $\kappa$ corresponds to
a nondegenerate trajectory of a billiard in $D$ starting and
ending in $A$, see Fig. \ref{fig:bill_ss1}.

Consider a billiard trajectory $\kappa$ in $D$ starting and ending
at the points $a_i,a_j\in A$. It may be a segment
$\kappa=(a_i,a_j)$ joining $a_i$ and $a_j$, or a broken line
\begin{equation}
\label{piece}
  \kappa = (a_i,p_1,\ldots,p_n,a_j),
  p_1,\ldots,p_n \in \Sigma,
\end{equation}
joining the points $a_i, p_1,\ldots,p_n, a_j$.
Then $(p_1,\dots,p_n)$ is a critical point  of the length function
$$
    l(a_i,p_1,\dots,p_n,a_j)
  = |a_i -  p_1| + |p_1 - p_2| + \ldots + |p_{n-1} - p_n|
         + |p_n - a_j|
$$
on $\Sigma^n$. The trajectory $\kappa$ is called nondegenerate if
$n = 0$ or the critical point $(p_1,\ldots,p_n)$ is nondegenerate.
We call the broken line (\ref{piece})  a quasi-trajectory  if it
is nondegenerate and, except the end points $a_i$ and $a_j$,
contains no point from the set $A$.

We denote by $a_\kappa^-=a_i$ and $a_\kappa^+=a_j$ the initial and
final points of the quasitrajectory, and by $v_\kappa^-$ and
$v_\kappa^+$ its initial and final velocity vectors issuing from
the points $a_i$ and $a_j$ respectively:
$$
    v_\kappa^-
  = \frac{p_1 - a_i}{|p_1 - a_i|}, \quad
    v_\kappa^+
  = \frac{a_j - p_n}{|a_j - p_n|}.
$$
In the case $n=0$ we should take
$v_\kappa^-=v_\kappa^+=(a_j-a_i)/|a_j-a_i|$.

If $D$ is convex, there always exist nondegenerate
quasitrajectories with $n=0$ joining $a_i$ and $a_j$. There also
always exist quasitrajectories with $n=1$ corresponding to the
minimum of $l$ on $\Sigma$. Generically, they will be
nondegenerate.

For $m=2$, by using the same argument as in the Birkhoff theorem
\cite{birk} on periodic trajectories of a convex billiard  (see
also \cite{Koz-Tre}), it is easy to show that, if $D$ is convex,
then for any $n\ge 1$ there are at least $2n$ orbits
(\ref{piece}).

\begin{figure}[ht]
\begin{center}
\includegraphics{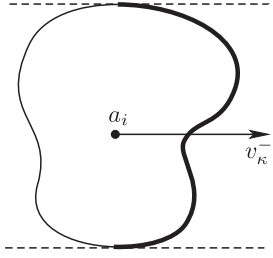}
\caption{The bold part of the boundary $\partial D_i$ is
 $a_i + \eps \phi_i(U_\kappa^-)$.}
\label{fig:bil_ss2}
\end{center}
\end{figure}

Any nondegenerate trajectory (\ref{piece}) is smoothly deformed
when we slightly change its starting and end points. Hence for
small $\eps\ge 0$ there exist smooth functions
$$
 p_1(x,y),\ldots,p_n(x,y),\qquad x,y\in S^{m-1},
$$
such that
$$
    B_\eps(x,y)
  = \big( a_i+\eps \phi_i(x), p_1(x,y),\ldots, p_n(x,y),
           a_j+\eps \phi_j(y) \big)
$$
is a trajectory of the billiard  in $D$ which coincides with
$\kappa$ when $\eps=0$. It will be contained in $\Omega_\eps$
provided the only common points of $B_\eps(x,y)$ with
$\partial\Omega_\eps$ are the end points $a_i+\eps q_i(x)$ and
$a_j+\eps q_j(y)$. This holds if $\eps$ is sufficiently small and
$x\in U_\kappa^-$, $y\in U_\kappa^+$, where $U_\kappa^\pm\subset
S^{m-1}$ are some  open sets, see Fig. \ref{fig:bil_ss2}. For
example, $U_\kappa^-$ is the set of $x$ such that the ray starting
at $\phi_i(x)$ in the direction of $v_{\kappa}^-$  does not cross
$\Sigma_i$.

The discrete Lagrangian corresponding to $\kappa$ is
\begin{eqnarray*}
       L_\kappa(x,y)
 &=&   \eps^{-1}l(B_\eps) \\
 &=& \eps^{-1} \big| a_i - p_1(x,y) + \eps \phi_i(x) \big|
    + \eps^{-1} \big| p_1(x,y) - p_2(x,y) \big|
    + \dots \\
 && +\, \eps^{-1} \big| p_{n-1}(x,y)
                        - p_n(x,y) \big|
    + \eps^{-1} \big| p_n(x,y) - a_j - \eps \phi_j(y) \big|.
\end{eqnarray*}

Let us show that
$$
    L_\kappa(x,y)
  = \eps^{-1} l(\kappa) + \langle v_\kappa^+,\phi_i(x)\rangle
                       - \langle v_\kappa^-,\phi_j(y)\rangle
                       + O(\eps),\qquad (x,y)\in S^{m-1}\times S^{m-1}.
$$
Indeed, by Hamilton's first variation formula,
$$
     d l(B_\eps)
  = \nu_\kappa^+ \,d q_+ - \nu_\kappa^-\, d q_-, \qquad
    d q_- = \phi_i(x)\,d\eps,\quad
    d q_+ = \phi_j(y)\,d\eps.
$$
Therefore
$$
\frac{\partial}{\partial\eps}\bigg|_{\eps=0}l(B_\eps)
  = \langle\nu_\kappa^+, \phi_j(y)\rangle
    - \langle\nu_\kappa^-, \phi_i(x)\rangle.
$$

Let the set of vertices of the graph $\Gamma$ be a finite collection of quasitrajectories $J = \{\kappa\}$.  We connect vertices $\kappa$ and $\kappa'$ with an edge  provided
\begin{enumerate}
\item end of $\kappa$ is the start of $\kappa'$: $a_{\kappa}^+=
a_{\kappa'}^- = a_k\in A$,

\item direction change: $v_{\kappa}^+\ne v_{\kappa'}^-$.

If $\Sigma_k$ is strictly convex, no more conditions are needed. If not, we need one more. Let $v\in\mR^m$ be the unit vector
$$
  v=\frac{ v_{\kappa'}^--v_{\kappa}^+}{| v_{\kappa'}^--v_{\kappa}^+|},
$$
see Fig. \ref{fig:nununu}. There always exists $s\in S^{m-1}$ such
that $v$ is the outer normal to the tangent plane
$T_{\phi_k(s)}\Sigma_k$ and $s\in U_{\kappa}^+\cap U_{\kappa'}^-$.

\begin{figure}[ht]
\begin{center}
\includegraphics{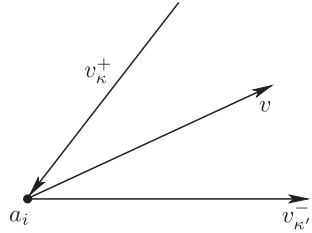}
\caption{The vectors $\nu_\kappa^+$, $\nu_{\kappa'}^-$ and $\nu$.}
\label{fig:nununu}
\end{center}
\end{figure}

\item  We assume that the second fundamental form $\langle
v,d^2\phi_k(s)\rangle$ of $\Sigma_k$ at $\phi_k(s)$ is
nondegenerate.
\end{enumerate}

Condition (3) always holds if  $\Sigma_k$ is strictly convex, then
there exists unique $s$ obtained by maximizing $\langle
\phi_k(s),v\rangle$. If $\Sigma_k$ is not convex there can be
several such points $s$, and then we will have several edges
$\gamma$ joining the vertices $\kappa,\kappa'$.

Thus we defined the graph with the set of vertices $J=\{\kappa\}$
and the set of edges $E=\{\gamma\}$. The Lagrangian $L_\kappa$ has
the anti-integrable form with $V_\kappa^+=\langle
\nu_\kappa^+,\phi_j\rangle$ and $V_\kappa^-=-\langle
\nu_\kappa^-,\phi_i\rangle$. Any path in the graph defines a code.
Now Theorem \ref{thm:anti} implies:

\begin{theo}
\label{theo:chen} (compare with \cite{Chen2004}). Suppose that
$\eps$ is sufficiently small. Then for any code $\bdk$ there is a
billiard trajectory $\bdx$ in $\Omega_\eps$ shadowing the chain of
quasitrajectories $\bdk=(\kappa_i)$. The orbit $\bdx$ is
hyperbolic.
\end{theo}

The billiard system in $\Omega_\eps$ is the ambient system for the
DLS we have constructed. Lower estimates for topological entropy
of the ambient system can be obtained with the help of the
topological Markov chain we obtain. A more detailed study of the
problem of topological entropy for a billiard with small
scatterers is presented in \cite{Chen2010}.

Note that if the scatterers are convex, then under certain conditions
the billiard will be hyperbolic, and then much stronger results hold,
in particular the metric entropy will be positive
\cite{Sin_bill}.

\subsection{Birkhoff-Smale-Shilnikov theorem  }

\label{sec:Smale}

As a general example of an application of the AI limit in
symplectic maps, we prove the existence of a chaotic hyperbolic
invariant set for a discrete dynamical system with a transverse
homoclinic orbit. This theorem goes back to Birkhoff \cite{birk} and Smale
\cite{Smale} and in the final form was proved by Shilnikov
\cite{Shil}, see \cite{KH,SSTC}.

We will prove the Birkhoff--Smale--Shilnikov theorem  for
symplectic maps. However, a general (nonsympectic)  map can be
reduced to a symplectic one by doubling the dimension, see the
remark in the end of the section.

Let $F:P\to P$  be a symplectic map of a $2m$-dimensional
symplectic manifold and $O$ a hyperbolic fixed point. Suppose that
there exist transverse homoclinic orbits to $O$. Take a finite
collection $\{\gamma_k\}_{k\in I}$ of them. We will show that
trajectories of $F$ which stay in a neighborhood of the homoclinic
set $\cup_{k\in I}\gamma_k\cup O$ are described by an
anti-integrable DLS. This provides a symbolic representation of
the trajectories in vicinity of the homoclinic set.

It is convenient to denote by $W^+$ and $W^-$ the stable and
unstable manifolds\footnote{The standard notation is $W^s$ and
$W^u$.} of $O$:
$$
W^\pm=\{x:F^n(x)\to O\;\mbox{as}\; n\to\pm\infty\}.
$$
 There exist symplectic coordinates $q,p$ in a
neighborhood $D$ of $O$ such that $O=(0,0)$ and the local stable
and unstable manifolds $W_\loc^\pm\subset D$  are Lagrangian
graphs over a ball $ U\subset\mR^m$:
\begin{equation}
\label{eq:Wloc} W^+_\loc=\{(q,p): q\in U, \; p=- \nabla
S_+(q)\},\quad W^-_\loc=\{(q,p): q\in U, \; p= \nabla S_-(q)\}.
\end{equation}

Let
$$
  \lambda_j,\lambda_j^{-1},\qquad 0<|\lambda_j|<1,\quad
  j=1,\dots,m,
$$
be the eigenvalues of $DF(O)$. Fix
$$
  \alpha\in (\max\{|\lambda_j|\},1).
$$
The local dynamics near $O$ is described by the next statement
which can be deduced from the Shilnikov lemma \cite{Shil} or strong the
$\lambda$-lemma \cite{Deng}.

\begin{lem}\label{lem:Shil}
Let $N$ be sufficiently large. Then for any $n\ge N$ and $q^\pm\in
U$ there exist $p^\pm$ such that $F^n(q^+,p^+)=(q^-,p^-)$ and
$F^i(q^+,p^+)\in D$ for $0\le i\le n$. The map $(q^+,p^+)\mapsto
(q^-,p^-)$ is a symplectic diffeomorphism $F^n:D_n^+\to D_n^-$ of
open sets in $D$ given by the generating function (defined up to a
constant)
\begin{equation}
\label{eq:un}
  S_n(q^+,q^-)=S_+(q^+)+S_-(q^-)+u_n(q^+,q^-),
\end{equation}
where $\|u_n\|_{C^2(U\times U)}\le C\alpha^n$.
\end{lem}

Thus
$$
  F^n(q^+,p^+)=(q^-,p^-)\quad
  \Leftrightarrow\quad
  p^+\,dq^+ - p^-\,dq^- = dS_n(q^+,q^-).
$$
For large $n$ the generating function $S_n$ of the local map
$F^n:D_n^+\to D_n^-$ has an anti-integrable form. To obtain
chaotic dynamics we need to come back to $O$ by using a global
return map along homoclinic orbits.

Let $\gamma_k$  be a transverse homoclinic orbit to the fixed
point $O$. Then there are points $z_k^\pm = (q_k^\pm,p_k^\pm)\in
\gamma_k\cap W_\loc^\pm$ and integers $m_k$ such that
$F^{m_k}(z_k^-)=z_k^+$.
 The   map $F^{m_k}$
from a neighborhood $G_k^-$ of $z_k^-$ to a neighborhood $G_k^+$
of $z_k^+$ is called the global map. Since the  map  is
symplectic, there exist functions $\Phi_k$ such that
\begin{equation}
\label{eq:glob}
  F^{m_k}(q^-,p^-) = (q^+,p^+)\quad
  \Leftrightarrow\quad
  p^+\,dq^+ - p^-\,dq^- = d\Phi_k.
\end{equation}
Perturbing the coordinates if needed, we  can assume
$\det(\partial q^+/\partial p^-)\ne 0$. Then we can locally
express $\Phi_k$ as a function $\Phi_k(q^-,q^+)$ on $U_k^-\times
U_k^+$, where $U_k^\pm\subset U$ is a neighborhood of $q_k^\pm$.
Then $\Phi_k$ is the generating function of the global map
$F^{m_k}:G_k^-\to G_k^+$.

Equations (\ref{eq:Wloc}) and (\ref{eq:glob}) imply that the
function
$$
  R_k(q^-,q^+) = S_-(q^-) + \Phi_k(q^-,q^+) + S_+(q^+)
$$
has a critical point $(q_k^-,q_k^+)$ which corresponds to the
homoclinic orbit $\gamma_k$. Since $\gamma_k$ is transverse, the
critical point is nondegenerate.

Trajectories of $F$ which stay in a neighborhood of the homoclinic
set correspond to trajectories of the DLS with Lagrangians
$\{\Phi_k,S_n\}_{k\in I,n\ge N}$. However, formally this system is
not anti-integrable. Application of Theorem \ref{thm:anti} is
cleaner if we consider another  DLS with $2m$ degrees of freedom.

Let $J=\{\kappa=(k,n):k\in I,\; n\ge N\}$. We define the DLS
$\{L_\kappa\}_{\kappa\in J}$ with the discrete Lagrangian
$L_\kappa$ representing the composition $F^n\circ F^{m_k}$ of the
global and local map. Set
$$
  L_\kappa(x,y)=\Phi_k(x)+S_n(x^+,y^-),\qquad
  x=(x^-,x^+)\in
  U_k^-\times U_k^+,\quad y=(y^-,y^+)\in U^2.
$$

\begin{rem}The Lagrangian $L_\kappa$ has $2m$ degrees
of freedom. We can replace it by the reduced Lagrangian
$$
    \tilde L_\kappa(x^-,y^-)
  = \mathrm{Crit}_{x^+}(\Phi_k(x^-,x^+))+S_n(x^+,y^-))
$$
with $m$ degrees of freedom. This requires an extra nondegeneracy
condition which can be satisfied by perturbing the coordinates.
\end{rem}

Orbits of $F$ shadowing the homoclinic chain $(\gamma_{k_i})$
correspond to critical points of the action functional
$$
     A_\bdk(\bdx)
 = \sum L_{\kappa_i}(x_i,x_{i+1}),\qquad
   x_i=(x_i^-,x_i^+)\in U_{k_i}^-\times U_{k_i}^+,\quad
   \kappa_i=(k_i,n_i).
 $$

 To obtain an anti-integrable  DLS, we replace
 $L_\kappa$  with a gauge equivalent Lagrangian (with the same
 action functional)
 $$
 \hat L_\kappa(x,y)=L_\kappa(x,y)+S_-(x^-)-S_-(y^-)=R_{k}(x) +
 O(\alpha^n),
 $$
 where  $R_{k}$ has a nondegenerate critical point $(q_k^-,q_k^+)$.
In the notation of  Theorem \ref{thm:anti} we have an
anti-integrable DLS defined by a graph $\Gamma$ with the set of
vertices $J$, and any two vertices are joined by an edge. If there
is no restriction on $n_i$ from above, then the graph will be
infinite, but it is easy to see that Theorem \ref{thm:anti_inf}
works. We obtain the Birkhof-Smale-Shilnikov theorem:

\begin{theo}
\label{thm:Smale} For any code $\kappa_i=(k_i,n_i)\in J$ the
corresponding homoclinic chain  $(\gamma_{k_i})$ is shadowed by a
unique hyperbolic trajectory of $F$ which follows $\gamma_{k_i}$,
then stays for $n_i$ iterations in a neighborhood of $O$, then
shadows $\gamma_{k_{i+1}}$ and so on. Thus $F$ has a chaotic
hyperbolic invariant set.
\end{theo}

 The standard graph which gives symbolic representation
of dynamics in a neighborhood of a transverse homoclinic orbit is
a bit different \cite{KH}: there is one vertex for every
homoclinic orbit $\gamma_k$ and one more for the fixed point $O$.
After $\gamma_{k_i}$, a path in the graph may stay at $O$ for
several steps $n_i$ before following $\gamma_{k_{i+1}}$. Thus a
path corresponds to a sequence $(k_i,n_i)$ as described above, so
dynamics is equivalent.

\begin{rem}
The Birkhof-Smale-Shilnikov theorem holds for a general
(nonsymplectic) map $f$. Indeed,  the map $q^+=f(q^-)$ can be
reduced to a symplectic map $F$ with the generating function
$S(q^-,p^+)=\langle f(q^-),p^+\rangle$ by introducing the
conjugate momentum:
$$
  p^- =\partial_{q^-} S,\quad
  q^+ = \partial_{p^+}S = f(q^-).
$$
If $f$ has a hyperbolic fixed point possessing  a transverse
homoclinic orbits, then so does $F$. Thus we proved Theorem
\ref{thm:Smale} also for nonsymplectic maps.
\end{rem}

\subsection{Shadowing a chain of invariant tori}

The following application of the AI limit is based on  \cite{Bol:homtori}.
Let $F:P\to P$ be a smooth symplectic diffeomorphism which has a
$d$-dimensional hyperbolic invariant torus $\Gamma$. Then $\Gamma$
is the image of a smooth embedding $h:\mT^d\to P$, and $F|_\Gamma$
is a translation with the rotation vector $\rho\in\mR^{d}$:
$$
  F(h(x))=F(x+\rho).
$$
If the rotation vector is Diophantine:
$$
 |\langle\rho,j\rangle - j_0| \ge \alpha |j|^{-\beta},\qquad
 \alpha,\beta>0,
$$
for any $j\in \mZ^{d}\setminus\{0\}$ and $j_0\in\mZ$, the torus is
said to be {\it Diophantine.} We assume the torus to be {\it
isotropic}, i.e.\ $\omega|_{\Gamma}=0$. If the torus is
Diophantine and the symplectic structure exact, this is always so.

\begin{df}\label{defin:hyp}
The torus $\Gamma$ is called hyperbolic if there exist two smooth
$(m-d)$-dimensional subbundles $E^{\pm}$ of the bundle $T_\Gamma
P$ such that
\begin{itemize}
\item $E^{\pm}$ are invariant for the linearized map $DF$, i.e.\
$DF(x) E^\pm_x = E^\pm_{F(x)}$ for all $x\in \Gamma$. \item The
linearized map is contracting on $E^+$  and expanding on $E^-$,
i.e.\ for some $c>0,\lambda>1$ and all $x\in \Gamma$,
$$
\|DF^k(x)\big|_{E_x^+} \|\le c\lambda^{-k}, \quad
\|DF^{-k}(x)\big|_{E_x^-} \| \le c\lambda^{-k}, \qquad k\in\mN.
$$
\end{itemize}
\end{df}

We fixed some Riemannian metric, and $\|\cdot\|$ is the operator
norm defined by this metric. Since $\Gamma$ is compact, the
definition is independent of the metric.

\begin{df}
\label{defin:nondeg} A   torus $\Gamma$ is {\it nondegenerate} if
all bounded trajectories of the linearized map are tangent to
$\Gamma$. Thus if $x\in \Gamma$, $v\in T_x M$, and
$\|DF^k(x)v\|\le c$ for all $k\in\mZ$, then $v\in T_x \Gamma$.
\end{df}

We can rewrite the definition of a hyperbolic torus in the
coordinate form.

\begin{df}
\label{defin:coord} An invariant torus $\Gamma$ is called
hyperbolic if in its tubular neighborhood $D$ there exist
symplectic coordinates $ x\in\mT^d$, $y\in \mR^d$, $z_\pm\in
\mR^{m-d}$  such that:
\begin{itemize}
\item $\omega|_D=dy\wedge dx+dz_+\wedge dz_-$; \item $\Gamma $ is
given by the equations $y=0,z_-=z_+=0$. \item The map $F|_D$ has
the form
\begin{equation}
\label{eq:normal} \left(
\begin{array}{c}
x\\
y\\
z_-\\
z_+
\end{array}
\right) \mapsto \left(
\begin{array}{l}
  x+\rho+Ay\\ y\\ B^*(x)^{-1} z_-\\ B(x)z_+
\end{array}
\right)+O_2(y,z_-,z_+).
\end{equation}
\item The Lyapunov exponents for the skew product map
$(x,z)\mapsto(x+\rho,B(x)z)$ are negative: there exists a norm
such that $\|B(x)\|\le \alpha<1$. \item The symmetric $d\times d$
matrix $A$ is constant.
\end{itemize}
\end{df}

For a Diophantine torus Definitions \ref{defin:hyp} and
\ref{defin:coord} are equivalent, see \cite{Bol-Tre:hyp}. If we
use Definition \ref{defin:coord}, then for the sequel  it is not
necessary to assume that the torus is Diophantine. However,
hyperbolic invariant tori arising in applications are usually
Diophantine. A hyperbolic torus is nondegenerate iff $\det A\ne
0$. By KAM-theory, hyperbolic nondegenerate Diophantine invariant
tori survive small  smooth exact symplectic perturbations of the
map.

\medskip

The hyperbolic torus $\Gamma$ has $m$-dimensional stable $W^+$ and
unstable $W^-$ manifolds  consisting of trajectories $F^n(z)$ that
are asymptotic to $\Gamma$ as $n\to +\infty$ and $n\to -\infty$
respectively.

Let $\{\Gamma_k\}_{k\in I}$ be a finite set of nondegenerate
hyperbolic tori of dimension $0<d_k<m$, and let $W_k^\pm$ be their
stable and unstable manifolds. A point $z\in W_j^-\cap W_k^+$ is
called a heteroclinic point, and its orbit
$\gamma=(F^n(z))_{n\in\mZ}$ is called  a heteroclinic orbit
connecting $\Gamma_j$ with $\Gamma_k$. It is called transverse if
$T_zW_j^-\cap T_zW_k^+=\{0\}$.

 Let $\{\gamma_\kappa\}_{\kappa\in J}$ be a finite
set of transverse heteroclinic\footnote{In the sequel heteroclinic
always means that homoclinic are included.} orbits connecting
pairs of tori from the set $\{\Gamma_k\}_{k\in I}$. Let $G$ be the
oriented graph with vertices $k\in I$ and edges $\kappa\in J$
corresponding to heteroclinic orbits.

\begin{theo}
\label{thm:shadow} Let $N>0$ be sufficiently large. Take any sequence of integers $(n_i)_{i\in\mZ}$, $n_i\ge N$.
Let $\kappa=(\kappa_i)_{i\in\mZ}$ be a path on the graph $G$
corresponding to the sequence $(\gamma_{\kappa_i})_{i\in\mZ}$ of
heteroclinics connecting $\Gamma_{k_{i-1}}$ with $\Gamma_{k_i}$.
Then there exists a (nonunique) trajectory of $F$ shadowing the
heteroclinic chain $(\gamma_{\kappa_i})_{i\in\mZ}$ and staying for
$n_i$ iterations in a neighborhood of the torus $\Gamma_{k_i}$
after shadowing $\gamma_{\kappa_i}$.
\end{theo}

Such results can be proved by Easton's  window method \cite{GR},
but the AI limit proof seems simpler and more constructive.
An alternative to AI limit approach is given by variational methods which do not require transverality
of heteroclinics (see e.g.\ \cite{Kalo-Levi,Bessi,Cheng-diff1,Cheng-diff2,Kal-Zhang}),
but these methods are limited to positive definite  Lagrangian systems.

Theorem~\ref{thm:shadow} is analogous to the shadowing lemma in
the general theory of hyperbolic sets of dynamical systems
\cite{KH}. The difference is that the set $\bigcup_{k\in
I}\Gamma_k\cup\bigcup_{\kappa\in J}\gamma_\kappa$ is not
hyperbolic, so the standard results of the hyperbolic theory do
not apply. Theorem~\ref{thm:shadow} can be used to construct
diffusion orbits in the problem of Arnold diffusion for an
$\eps$-small perturbation of an priori unstable system (see
section \ref{sec:sm_mult}). However, this is possible only away
from strong resonances and moreover the diffusion speed given by
Theorem~\ref{thm:shadow} will be of order $O(\eps^2)$ which is
much weaker than the diffusion speed $O(\eps/|\log\eps|)$ obtained
by using the AI limit approach for the separatrix map, see section
\ref{sec:sm_mult}.

 When the hyperbolic tori  are  hyperbolic fixed points, we
obtain a weaker version of Theorem \ref{thm:Smale}. In general the
proof is similar. We will replace Lemma  \ref{lem:Shil} with the
following Lemma \ref{lem:local}.

Let $\Gamma$ be a hyperbolic torus. Choose symplectic coordinates
$(q,p)$ in a small tubular neighborhood $D$  of $\Gamma$ so that
the  Lagrangian local stable and unstable manifolds are given by
generating functions:
$$
  W_\loc^+=\{(q,p)\in D:q\in U,\;p=-\nabla S^+(q)\},\quad
  W_\loc^-=\{(q,p)\in D:q\in U,\; p=\nabla S^-(q)\}.
$$
The projection of $\Gamma=W_\loc^+\cap W_\loc^-$ to $U$ is the
torus $T=\{q:\nabla (S^++S^-)(q)=0\}$.

Since $U$ is a tubular neighborhood of the torus $T$, we have a
(noncanonical) projection $U\to\ \mT^d$. Let $\widetilde U$ be the
universal coverings of $U $. Then we have a map $\phi:\widetilde
U\to \mR^d$.   The map $F$ is lifted to a map $F:\widetilde
D\to\widetilde D$.

\begin{lem}
\label{lem:local} Let   $r>0$.  There exists $N>0$ such that:
\begin{itemize}
\item For  all $n\ge N$ and all $(q_+,q_-)$ in the set
\begin{equation}
\label{eq:Y} Y^{n}=\{(q_+,q_-)\in\widetilde U\times \widetilde U:
|\phi(q_-)-\phi(q_+)-n\rho|\le r\}
\end{equation}
there exist unique $z_+=(q_+,p_+)\in \widetilde D$ and
$z_-=(q_-,p_-)\in \widetilde D$ such that $F^n(z_+)=z_-$ and
$F^j(z_+)\in\widetilde D$ for $j=0,1,\dots, n$. \item The  map
$F^n:(q_+,p_+)\to (q_-,p_-)$ is given by the  generating function
$Q^n$:
$$
F^n(q_+,p_+)=(q_-,p_-)\quad\Leftrightarrow\quad p_-=\partial_{q_-}
Q^n(q_+,q_-),\quad p_+=-\partial_{q_+} Q^n(q_+,q_-).
$$
\item  $Q^n$  is a smooth function on $Y^{n}$ and has the form
\begin{equation}
\label{eq:Shil2}
Q^n(q_+,q_-)=S^+(q_+)+S^-(q_-)+n^{-1}v_n(q_+,q_-),
\end{equation}
where $\|v_n\|_{C^2(Y^{n})}\le Cr^2$. The constant $C>0$ is
independent of $n$.
\end{itemize}
\end{lem}

We do not need the rotation vector $\rho$ to be Diophantine or
nonresonant, but it is essential that $A$ is nondegenerate.
Lemma~\ref{lem:local} is proved in \cite{Bol:homtori} for $d=m-1$.
For any $0<d<m$ the proof is similar. For $d=0$, Lemma
\ref{lem:Shil} gives a stronger estimate.

For every torus $\{\Gamma_k\}_{k\in I}$ we define the tubular
neighborhoods $U_k$ and $D_k$ and generating functions $S_k^\pm$
and $Q_k^n$ as above. Lemma~\ref{lem:local} provides the discrete
Lagrangian describing the local map near the hyperbolic tori
$\Gamma_k$. Next we define the discrete Lagrangian describing the
global map along a transverse heteroclinic $\gamma_\kappa$.

Take two  tori $\Gamma_{j}$, $\Gamma_k$ joined by a heteroclinic
(homoclinic if $k=j$) orbit $\gamma_\kappa\subset W_j^-\cap
W_k^+$.
 There exist  points
$$
z_\kappa^-=(q_{\kappa}^-,p_\kappa^-)\in\gamma_\kappa\cap
D_{j},\quad z_\kappa^+=(q_\kappa^+,p_\kappa^+)\in\gamma_\kappa\cap
D_k
$$
and
$m_\kappa\in\mN$ such that
$F^{m_\kappa}(z_{\kappa}^-)=z_\kappa^+$.

The  symplectic coordinates  can be chosen  in such a way that the
global map $F^{m_\kappa}$ from a neighborhood of $z_\kappa^-$ to a
neighborhood of $z_\kappa^+$ is given by a generating function
$\Phi_\kappa$ on $X_{\kappa}^-\times X_\kappa^+$, where
$X_\kappa^\pm$ is a small neighborhood of $q_\kappa^\pm$:
\begin{equation}
\label{eq:gen2} F^{m_\kappa}(q_-,p_-)
=(q_+,p_+)\quad\Leftrightarrow\quad p_+=\partial_{q_+}
\Phi_\kappa(q_-,q_+),\quad p_-=-\partial_{q_-}
\Phi_\kappa(q_-,q_+).
\end{equation}
Since the heteroclinic $\gamma_\kappa$ is transverse, exactly as
in the previous section we conclude that $(q_\kappa^-,q_\kappa^+)$
is a nondegenerate critical point of the function
$$
  R_\kappa(q_-,q_+)=S_{j}^-(q_-)+\Phi_\kappa(q_-,q_+)+S_k^+(q_+),\qquad
  (q_-,q_+)\in X_{\kappa}^-\times X_\kappa^+.
$$

 Let $\pi:\widetilde U_k\to U_k$ be the universal covering with the transformation group
$$
\tau_v:\widetilde U_k\to\widetilde U_k,\qquad v\in\mZ^{d_k},\quad
d_k=\dim\Gamma_k.
$$
Fix connected components $Z_\kappa^\pm$ of
$\widetilde X_\kappa^\pm=\pi^{-1}(X_\kappa^\pm)$. Then $\widetilde
X_\kappa^\pm=\cup_{v\in\mZ^{d_k}}\tau_v Z_\kappa^\pm$. To decrease
the nonuniqueness, we assume that $Z_\kappa^\pm$ intersect with
the same fundamental domain $K$ of the group action. Then
\begin{equation}
\label{eq:diam}
  \max\{d(x,y):x\in Y_\kappa^+,y\in Y_\kappa^-\}\le 2\,\mathrm{diam}\, K.
\end{equation}

Let $\rho_k\in\mR^{d_k}$ be the rotation vector of the torus
$\Gamma_k$. Take $R>\sqrt{d_k}$. For any $n\in\mN$, there exists
$v\in\mZ^{d_k}$ such that
$$
\label{eq:m} |v-n\rho_k |\le R .
$$
Then
$$
|\phi(y)-\phi(x)-n\rho_k|\le r=R+2\,\mathrm{diam}\,K
$$
for all $x\in Z_{\kappa}^+$ and $y\in \tau_{v}Z_\kappa^-$. Let
$N>0$ be sufficiently large and $n\ge N$. Then
$$
Z_{\kappa}^+\times \tau_{v}Z_\kappa^-\subset Y_k^{n},
$$
where  $Y_k^{n}$ is the set (\ref{eq:Y}) corresponding to the
torus $\Gamma_k$.

Let $\sigma=(\kappa,n,v)$, where $n\ge N$ and $v$ is chosen as
above. We define the discrete Lagrangian $L_\sigma$ by
$$
L_\sigma(x,y)=\Phi_\kappa(x)+Q_{k}^{n}(x_+,\tau_v y_-),\qquad
x=(x_-,x_+),\quad y=(y_-,y_+),
$$
where $Q_{k}^n$ is the function (\ref{eq:gen2}) corresponding to
the torus $\Gamma_k$. The Lagrangian $L_\sigma$ represents
trajectories which shadow the heteroclinic orbit $\gamma_\kappa$
and then travel in a neighborhood of $\Gamma_k$ for $n$ iteration
of the map $F$. As in the previous section, $L_\sigma$ has $2m$
degrees of freedom, but, if desired, we can replace it by a
Lagrangian with $m$ degrees of freedom.

To obtain an anti-integrable Lagrangian, we make a gauge
transformation:
$$
\hat L_\sigma(x,y)=L_\sigma(x,y)+S_j^-(x_-)-S_k^-(y_-).
$$
By (\ref{eq:Shil2}), $\hat L_\sigma=R_\kappa(x)+O(n^{-1})$, and
$R_\kappa$ has a nondegenerate critical point
$(q_\kappa^-,q_\kappa^+)$.

Suppose we are given a chain of heteroclinics
$(\gamma_{\kappa_i})_{i\in\mZ}$, where $\gamma_{\kappa_i}$ joins
$\gamma_{k_{i-1}}$ and $\gamma_{k_i}$.   Consider the infinite
product
$$
Z= \prod_{i\in\mZ} Z_{\kappa_i}^-\times Z_{\kappa_i}^+
$$
which is the set of sequences
$$
\bdx=(x_i),\qquad x_i=(x_i^-,x_i^+)\in Z_{\kappa_i}^-\times
Z_{\kappa_i}^+.
$$

For a sequence $(n_i,v_i)$ such that $n_i\ge N$ and
$|v_i-n_i\rho_{k_i}|\le R$, set $\sigma_i=(\kappa_i,n_i,v_i)$ and
define on $Z$ a formal functional
$$
A_\sigma(\bdx)=\sum_{i\in\mZ}
L_{\sigma_i}(x_i,x_{i+1})=\sum_{i\in\mZ} \hat
L_{\sigma_i}(x_i,x_{i+1}).
$$
Its critical points correspond to shadowing orbits. Since $\hat
L_\sigma$ is anti-integrable, Theorem~\ref{thm:shadow} follows
from Theorem \ref{thm:anti}. \qed

\section{AI limit in continuous Lagrangian systems}
\label{sec:ail_cls}

In the previous sections we gave examples of application of AI
limit in DLS. Now we  discuss applications of AIL to continuous
Lagrangian systems, two examples for autonomous systems and one
for time dependent.

\subsection{Turaev--Shilnikov theorem for Hamiltonian systems}
\label{sec:Turayev}

As an example of the anti-integrable limit in autonomous
Lagrangian or Hamiltonian systems, we prove  the
Turayev--Shilnikov theorem \cite{Tur-Shil} for a Hamiltonian
system with a hyperbolic equilibrium. This situation is more
delicate than for hyperbolic fixed points of a symplectic map:
system with transverse homoclinics may be integrable \cite{Devaney}. The
proof we give  follows the approach in \cite{Bol-Rab:revers}.

Consider a Hamiltonian system with Hamiltonian $H$ on a symplectic
manifold $P$. Let $\phi^t$ be the flow and $O$ a hyperbolic
equilibrium point. We may assume that $H(O)=0$. Let
 $$
 \pm\lambda_j,\qquad  0<\Re\lambda_1\le \dots\le \Re\lambda_m,
 $$
 be the eigenvalues. Suppose that the eigenvalue $\lambda_1$ with the smallest real part is
 real and
 $$
0<\lambda_1<\Re\lambda_2.
$$
The case of complex $\lambda_1$ is somewhat simpler, it is studied
in \cite{Lerman} and \cite{Buf-Sere}.

Let $v_+$ be an eigenvector corresponding to $-\lambda_1$ and
$v_-$ an eigenvector corresponding to $\lambda_1$. We fix a metric
and assume $|v_\pm|=1$.  To decrease nonuniqueness of the
eigenvectors we may assume that $\omega(v_-,v_+)>0$, where
$\omega$ is the symplectic 2-form. Then the eigenvectors $v_\pm$
are uniquely defined up to the change $(v_+,v_-)\mapsto
(-v_+,-v_-)$.

Let $W^{\pm}$ be the stable and unstable manifolds of the
equilibrium $O$. They contain   strong stable and unstable
manifolds $W^{++}\subset W^+$ and $W^{--}\subset W^{-}$
corresponding to the eigenvalues $\pm\lambda_j$ with $j>1$. Any
trajectory in $W^+\setminus W^{++}$ is tangent to $v_+$ as
$t\to+\infty$, and any trajectory in $W^-\setminus W^{--}$ is
tangent to $v_-$ as $t\to -\infty$.

Let $\gamma(t)=\phi^t(z)$ be a homoclinic trajectory:
$\lim_{t\to\pm\infty}\gamma(t)=O$.   It is called transverse  if
$T_{\gamma(0)}W^+\cap T_{\gamma(0)}W^-=\mR \dot\gamma(0)$. Then
intersection of $W^+$ and $W^-$ along $\gamma$ is transverse in
the energy level $H=0$.

Suppose  that there exist several transverse homoclinic orbits
$\{\gamma_\kappa\}_{\kappa\in J}$  which do not belong to the
strong stable and unstable manifolds. Then $\gamma_\kappa$ will be
tangent to $v_\pm$ as $t\to\pm\infty$ respectively:
$$
\lim_{t\to\infty}\frac{\dot\gamma_\kappa(t)}{|\dot\gamma_\kappa(t)|}=s_\kappa^\pm
v_\pm,\qquad s_\kappa^\pm= +\;\mbox{or}\; -.
$$

Let us define a graph $\Gamma_-$ as follows. The vertices
$\kappa\in J$ correspond to transverse homoclinics
$\gamma_\kappa$. We join the vertices $\kappa$ and $\kappa'$ by an
edge if $s_\kappa^+=s_{\kappa'}^-$. The graph $\Gamma_+$ is
defined in the same way but we join $\kappa,\kappa'$ by an edge if
$s_\kappa^+=-s_{\kappa'}^-$.

\begin{theo}\label{thm:tur-shil}
There exists $\eps_0>0$ such that for any $\eps\in (0,\eps_0)$ and
any path $\bdk=(\kappa_i)_{i\in\mZ}$ in the graph $\Gamma_+$ there
exists a hyperbolic  trajectory with energy $H=\eps$ shadowing the
homoclinic chain $(\gamma_{\kappa_i})$. Similarly, for any path
$\bdk=(\kappa_i)_{i\in\mZ}$ in the graph $\Gamma_-$ the chain
$(\gamma_{\kappa_i})$ is shadowed by a hyperbolic trajectory with
energy $H=-\eps$.
\end{theo}

To prove this result we construct an anti-integrable  DLS which
describes shadowing trajectories. Then it remains to use Theorem
\ref{thm:anti}.

There exist symplectic coordinates $q,p$ in a neighborhood $D$ of
$O$ such that $O=(0,0)$ and the local stable and unstable
manifolds $W_\loc^{\pm}$ are Lagrangian graphs over a small ball
$U$ in $\mR^m$:
\begin{equation}
\label{eq:Wpm}
  W^+_\loc=\{(q,p): q\in U,\; p=-\nabla S_+(q)\},\quad
  W^-_\loc=\{(q,p): q\in U,\; p=\nabla S_-(q)\}.
\end{equation}
Set $S_\pm(0)=0$. We may assume that the coordinates are chosen in
such a way that the strong local stable and unstable manifolds are
given by
$$
  W_\loc^{++}=\{(q,p)\in W_\loc^+:q_1=0\},\quad
  W_\loc^{--}=\{(q,p)\in W_\loc^-:q_1=0\}.
$$
Fix small $\delta>0$ and for $s=+$ or $s=-$ set
$$
  U_s=\{q\in U:sq_1>\delta |q|\},\qquad
  W_{s}^{\pm}=\{(q,p)\in W_\loc^{\pm}: q\in U_s\}.
$$
We choose the signs in such a way that $v_\pm$ points towards
$W^\pm_+$. Then the homoclinic orbit $\gamma_\kappa$ satisfies
$\gamma_\kappa(t)\in W_{s_\kappa^\pm}^\pm$ for $t\to\pm\infty$.

\begin{figure}[ht]
\begin{center}
\includegraphics{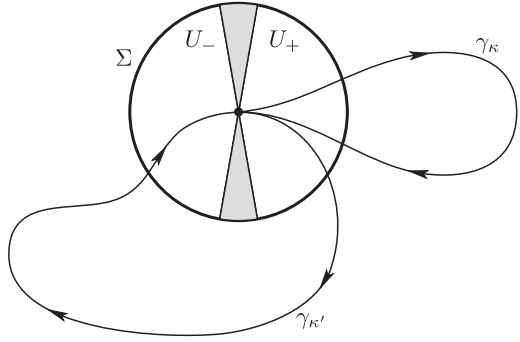}
\caption{Homoclinics $\gamma_\kappa$ and $\gamma_{\kappa'}$ with
$s_\kappa^-=s_\kappa^+=s_{k'}^-=+1$ and $s_{\kappa'}^+=-1$.}
\label{fig:Turayev}
\end{center}
\end{figure}

Take any $\alpha\in (0,1)$.

\begin{lem}\label{lem:connect}
Let $\eps_0>0$ be sufficiently small and  $\eps\in (0,\eps_0)$.
Let $s=+$ or $s=-$. For any $q_\pm\in U_{s}$ there exist $p_\pm$ and
$\tau>0$ such that $H(q_\pm,p_\pm)=- \eps$ and
$\phi^\tau(q_+,p_+)=(q_-,p_-)$. The trajectory
$\gamma(t)=\phi^t(q_+,p_+)$, $0\le t\le\tau$, stays in $D$ for
$0\le t\le\tau$. The Maupertuis action of the trajectory $\gamma $
has the form
$$
S(q_+,q_-,\eps)=\int_\gamma p\,dq=S_+(q_+)+S_-(q_-)+
w(q_+,q_-,\eps),
$$
 where $\|w\|_{C^2}\le  C\eps^\alpha$. If $q_+\in U_s$ and $q_-\in
 U_{-s}$, then the same is true but the trajectory will have
 energy $H=\eps$.
\end{lem}

Hence for small $\eps>0$ the local flow near $O$ is described by
an anti-integrable discrete Lagrangian. Lemma \ref{lem:connect} is
proved in \cite{Bol-Rab:revers}.

\medskip
\noindent{\it Proof of Theorem \ref{thm:tur-shil}.} For
definiteness we will consider shadowing trajectories with small
negative energy $H=-\eps$.

Let $\Sigma=\partial U$. The homoclinic orbit $\gamma_\kappa$
exits $W_\loc^-$ at the point $(q_\kappa^+,p_\kappa^+)$ and enters
$W_\loc^+$ at the point $(q_\kappa^+,p_\kappa^+)$, where
$q_\kappa^\pm \in \Sigma$. Then there are time moments $t_\kappa$
such that
$\phi^{t_\kappa}(q_\kappa^-,p_\kappa^-)=(q_\kappa^+,p_\kappa^+)$.
We can assume that $\det(\partial q_+/\partial p_-)\ne 0$. Then
there are neighborhoods $X_{\kappa}^\pm$ of $q_\kappa^\pm$ in
$\Sigma$ such that for $(q_-,q_+)\in X_\kappa^-\times X_\kappa^+$
 and any $\eps\in (0,\eps_0)$ there
exist $p_\pm$ such that $H(q_\pm,p_\pm)=-\eps$ and a trajectory $
\phi^t(q_-,p_-)$, $0\le t\le \tau=t_\kappa+O(\eps)$, with energy
$-\eps$ such that $\phi^\tau(q_-,p_-)=(q_+,p_+)$.  The map
$(q_-,p_-)\mapsto (q_+,p_+)$ has a generating function
$\Phi_\kappa(q_-,q_+,\eps)$ such that
$d\Phi_\kappa=p_+\,dq_+-p_-\,dq_-$.

For $\eps=0$, (\ref{eq:Wpm}) implies that the function
$$
    R_\kappa(q_-,q_+)
  = S_-(q_-) + \Phi_\kappa(q_-,q_+,0) + S_+(q_+), \qquad
    (q_-,q_+)\in X_\kappa^-\times X_\kappa^+,
$$
has a critical point $(q_\kappa^-,q_\kappa^+)$ corresponding to
the homoclinic orbit $\gamma_\kappa$. Since $\gamma_\kappa$ is
transverse, the critical point is nondegenerate.

For small $\eps>0$, trajectories with $H=-\eps$ shadowing
heteroclinic orbits $\{\gamma_\kappa\}$ correspond to trajectories
of a DLS $\{\Phi_\kappa,S\}$ with Lagrangians $\Phi_\kappa$ and
$S$ defined on open sets in $\Sigma$. It has $m-1$ degrees of
freedom. To apply Theorem \ref{thm:anti} it is convenient to
introduce a DLS with $2m-2$ degrees of freedom.

We may assume that $X_\kappa^\pm\subset U_{s_\kappa^\pm}$. Define
the discrete Lagrangian $L_\kappa$ on $X_\kappa^-\times
X_\kappa^+\times (U_{s_\kappa^+}\cap\Sigma)$ by
$$
L_\kappa(x,y,\eps)=\Phi_\kappa(x,\eps)+S(x_+,y_-,\eps),\qquad
x=(x_-,x_+),\; y=(y_-.y_+).
$$
This Lagrangian has $2m-2$ degrees of freedom.

For a given path $\bdk=(\kappa_i)$ in the graph $\Gamma_-$, we
obtain the discrete action functional
$$
A_\bdk(\bdx)=\sum  L_{\kappa_i}(x_i,x_{i+1},\eps),\qquad
x_i=(x_i^-,x_i^+).
$$
As explained in section  \ref{ssec:DLS_flow}, trajectories with
$H=-\eps$ shadowing the homoclinic chain $(\gamma_{\kappa_i})$
correspond to critical points of $A_\bdk$, i.e.\   trajectories of
the DLS $\calL=\{L_\kappa\}_{\kappa\in J}$.

By a gauge transformation, we can replace  $L_\kappa$ with an
anti-integrable Lagrangian
$$
\hat L_{\kappa}(x,y,\eps)=L_{\kappa}(x,y)+S_-(x_-)-S_-(y_-) =
R_{\kappa}(x)+O(\eps^\alpha),
$$
where the function $R_{\kappa}$ has a nondegenerate critical point
$(q_\kappa^-,q_\kappa^+)$. Now Theorem \ref{thm:tur-shil} follows
from Theorem \ref{thm:anti}. \qed

\medskip

 Theorem \ref{thm:tur-shil}  was formulated in \cite{Tur-Shil} (in a different form),
 and proved in  \cite{Bol-Rab:revers} (for positive definite Lagrangian systems)  and in \cite{Turayev:RCD}.
  In  \cite{Bol-Rab:revers}   variational
methods were used, so transversality of the homoclinics wasn't
assumed.  Similar results hold for systems with several hyperbolic
equilibria on the same energy level.

Consider for example a Lagrangian system on a compact manifold $M$
with $L(q,\dot q)=\frac12\|\dot q\|^2-V(q)$, where $V$ attains its
maximum on the set $A=\{a_1,\dots, a_n\}$ and each maximum point
is nondegenerate. There exist many minimal heteroclinics joining
the points in $A$, for example any pair of points can be joined by
a chain of heteroclinics.  Suppose that minimal heteroclinics
$\gamma_\kappa$ do not belong to the strong stable or unstable
manifolds, and define the numbers $s_\kappa^\pm$ as above. Define
the graph $\Gamma_-$ by joining $\kappa$ and $\kappa'$ by an edge
if $s_\kappa^+=s_\kappa^-$. Then for small $\eps>0$ and any path
in the graph $\Gamma_-$ there exists a trajectory with energy
$-\eps$ shadowing the corresponding chain of heteroclinics. See
\cite{Baesens} for  interesting concrete examples.

\subsection{n-center problem with small masses}
\label{sec:centers}

The problem we consider in this section  is somewhat similar to
the one discussed in Section \ref{ssec:small_scat}: instead of
small scatterers we have small singularities of the potential.

Let $M$ be a smooth manifold and $A=\{a_1,\dots,a_n\}$ a finite
set in $M$. Consider a Lagrangian system with the configuration
space $M\setminus A$ and the Lagrangian
\begin{equation}
\label{eq:Leps} L(q,\dot q,\eps)=L_0(q,\dot q)-\eps V(q) .
\end{equation}
We assume that $L_0$  is  smooth on $TM$ and   quadratic in the
velocity:
\begin{equation}
\label{eq:L} L_0(q,\dot q)=\frac12\|\dot q\|^2+\langle w(q),\dot
q\rangle-W(q),
\end{equation}
where $\|\;\|$  is a Riemannian metric on $M$.

The potential $V$ is a smooth function on $M\setminus A$  with
Newtonian singularities:  in a small ball $U_k$ around $a_k$,
\begin{equation}
\label{eq:sing} V(q)=-\frac {\phi_k(q)}{\dist(q,a_k)},
\end{equation}
where $\phi_k$ is a smooth positive function on $U_k$. The
distance is defined by means of the Riemannian metric $\|\;\|$. We
call the system with the Lagrangian (\ref{eq:Leps}) the $n$-center
problem. For $\eps=0$ the limit system with the Lagrangian
(\ref{eq:L})  has no singularities.

Let
\begin{equation}
\label{eq:energy} H(q,\dot q,\eps)=H_0+\eps V,\qquad H_0(q,\dot
q)=\frac12\|\dot q\|^2+W(q)
\end{equation}
be the energy integral. We fix $E$ such that $D=\{W<E\}$ contains
the set $A$ and study the system   on the energy level $\{H=E\}$.

We call a trajectory $\gamma:[a,b]\to D$  of the limit system
with the Lagrangian (\ref{eq:L}) a nondegenerate collision orbit
if $\gamma(a),\gamma(b)\in A$, $\gamma(t)\notin A$ for $a<t<b$,
and the endpoints are nonconjugate along $\gamma$ on the energy
level $H=E$, i.e.\ for the Maupertis action functional.

Suppose there are several nondegenerate collision orbits
$\{\gamma_\kappa\}_{\kappa\in J}$ connecting points
$a_\kappa^-,a_\kappa^+\in A$. We denote by $v_{\kappa}^-$ and
$v_\kappa^+$ the initial and final velocity of $\gamma_\kappa$.
Consider a graph $\Gamma$ with vertices $J$. We join $\kappa$,
$\kappa'$ by an edge if $a_\kappa^+=a_{\kappa'}^-$ and
$v_{\kappa}^+\ne \pm v_{\kappa'}^-$.

The next result is proved in \cite{Bol-Mac:centers}.

\begin{theo}
\label{thm:centers} There exists $\eps_0>0$ such that for all
$\eps\in(0,\eps_0]$ and any   path $\kappa=(\kappa_i)_{i\in\mZ}$
in the graph $\Gamma$ there exists a unique (up to a time shift)
trajectory of energy $E$ shadowing the chain
$(\gamma_{\kappa_i})_{i\in\mZ}$ of collision orbits.
\end{theo}

Hence there is an invariant subset in $\{H=E\}$ on which the
system is  conjugate to a suspension
 of a  topological Markov chain.
The topological entropy is positive provided that the graph
$\Gamma$ has a connected branched subgraph.

Conditions of Theorem~\ref{thm:centers} are formulated in terms of
the Lagrangian $L_0$ and the set $A$ only, not involving the
potential $V$ provided it has Newtonian singularities at $P$. We
can also add to $L$ a smooth   $O(\eps)$-small perturbation.

\begin{cor}
Suppose $M$ is a closed manifold and
$$
E> \min_{q\in M}\left(\frac12 \|w(q)\|^2+W(q)\right).
$$
Then for any $n\ge 2$,  almost all points $a_1,\dots,a_n\in M$,
and small $\eps>0$, there exists a chaotic hyperbolic invariant
set of trajectories of energy $E$ close to chains of collision
orbits.
\end{cor}

Indeed, the assumption implies that the Jacobi metric
$$
ds_E=\sqrt{2(E-W(q))}\|dq\|+\langle w(q),dq\rangle
$$
is a positive definite Finsler metric on $M$, and by Morse theory
any two generic points can be connected by an infinite number of
nondegenerate geodesics, i.e.\ trajectories of energy $E$. Thus
any two points in a generic finite set $A\subset M$ can be
connected by an infinite number of nondegenerate collision
trajectories of energy $E$, and any other points in $A$ do not lie
on these trajectories.

If $n$ is large enough, then chaotic trajectories exist for purely
topological reasons \cite{Bol-Neg:reg}, so genericity of $A$ and
smallness of $\eps$ is not needed.

\begin{rem} There is another  corollary of Theorem \ref{thm:centers}
for systems with no small parameter:
$$
  L(q,\dot q) = \frac12\|\dot q\|^2 - V(q),
$$
but with large energy $E=\eps^{-1}$. After a time change $s= \eps^{-1/2}t$ we obtain the Lagrangian (\ref{eq:Leps}) with $L_0=\|\dot q\|^2/2$. Hence the conclusion of Theorem \ref{thm:centers} holds.
See \cite{Knauf} for the classical $n$ center problem.
\end{rem}

\medskip

Next we prove   Theorem~\ref{thm:centers}. First consider the
limit system with $\eps=0$. Let $\Sigma_k=\partial U_k$. For any
$x\in \Sigma_k$, there is a unique trajectory $\gamma_x^+$ of
energy $E$ connecting $x$ with $a_k$ in $U_k$ and a unique trajectory $\gamma_x^-$ of energy $E$ connecting $a_k$
with $x$. Denote
\begin{eqnarray}
S_k^\pm(x)= \int_{\gamma_x^\pm}p\,dq=\int_{\gamma_x^\pm}ds_E.
\label{eq:Sp}
\end{eqnarray}
Then $S_k^\pm$ are smooth functions on $\Sigma_k$.

Denote by $u^+(x)$ and $u^-(x) $ the velocity vectors  of
$\gamma_x^\pm$ at the point $a_k$. Fix arbitrary small $\delta>0$
and let
\begin{equation}
\label{eq:X} X_k=\{(x,y)\in\Sigma_k^2\mid
\|u_+(x)-u_-(y)\|\ge\delta\}.
\end{equation}

\begin{lem}
\label{lem:connect2} Suppose $\eps_0>0$ is sufficiently small and
let $\eps\in (0,\eps_0]$.
\begin{itemize}
\item For any $(x,y)\in X_k$, there exists a unique trajectory
$\gamma=\gamma_{x,y}^\eps $  of energy $E$
 connecting $x$ with $y$ in $U_k$.
\item $\gamma$ smoothly depends on $x,y$. \item The Maupertuis
action
\begin{equation}
\label{eq:S} S_k(x,y,\eps)= \int_\gamma p\,dq
\end{equation}
is a smooth function on $X_k\times(0,\eps_0]$ and, modulo a
constant $\eps\log\eps$,
$$
  Q_k(x,y,\eps)=S_k^-(x)+S_k^+(y)+\eps u_k(x,y,\eps),
$$
where $u_k$ is uniformly $C^2$ bounded on $X_k$ as $\eps\to 0$.
\end{itemize}
\end{lem}

The proof is  based on regularization of singularities and a
version of Lemma \ref{lem:Shil}. In fact
 all positive eigenvalues for the regularized system will be  equal, so
Lemma \ref{lem:Shil} needs  modification, see
\cite{Bol-Mac:centers}.

For any $\kappa\in J$, let $x_\kappa\in\Sigma_{\kappa^-}$ and
$y_\kappa\in\Sigma_{\kappa^+}$ be the intersection points of the
collision orbit $\gamma_\kappa$ with $\Sigma_{\kappa^-}$ and
$\Sigma_{\kappa^+}$ respectively.
 If the spheres $\Sigma_{k}$ are small enough,
the points $x_\kappa$ and $y_\kappa$ are nonconjugate  along
$\gamma_\kappa$.

Let $U_\kappa^-\subset\Sigma_{\kappa^-}$ be a small neighborhood
of $x_\kappa$, and  $U_\kappa^+\subset\Sigma_{\kappa^+}$ a small
neighborhood of $y_\kappa$. Taking $\delta>0$ small enough, it can
be assumed that for any edge $(\kappa,\kappa')$ we have
$U_{\kappa}^+\times U_{\kappa'}^-\subset X_k$, where
$a_k=a_\kappa^+=a_{\kappa'}^-$. If the neighborhoods
$U_\kappa^\pm$ are small enough and $\eps\in (0,\eps_0)$, any
points $x\in U_\kappa^-$ and $y\in U_\kappa^+$ are joined by  a
unique trajectory $\beta_\eps$ of energy $H=E$ which is close to
$\gamma_\kappa$. Let
$$
\Phi_\kappa(x,y,\eps)=\int_{\beta_\eps} p\,dq
$$
be its action. Then $\Phi_\kappa$ is a  smooth function on
$U_\kappa^-\times U_\kappa^+$.

\begin{lem}
\label{lem:nondeg} The function
$R_\kappa(x,y)=\Phi_\kappa(x,y,0)+S_{\kappa^-}^-(x)+S_{\kappa^+}^+(y)$
on $U_\kappa^-\times U_\kappa^+$ has a nondegenerate critical
point $(x_\kappa,y_\kappa)$.
\end{lem}

Lemma~\ref{lem:nondeg} follows from the assumption that
$\gamma_\kappa$ is a nondegenerate critical point of the action
functional. Indeed, $R_k(x,y)$ is the Maupertuis action of the
piecewise smooth trajectory of the limit system obtained by gluing
together the trajectories $\gamma_x^-$, $\beta_0$ and
$\gamma_y^+$. Hence $R_k$ is a restriction of the action
functional to a finite-dimensional submanifold consisting of broken
trajectories (with break points $x,y$) connecting $a_{\kappa^-}$
with $a_{\kappa^+}$. \qed

\medskip

Define the discrete Lagrangian with $2m-2$ degrees of freedom by
$$
L_\kappa(z_-,z_+,\eps)=\Phi_\kappa(z_-)+Q_{\kappa^+}(y_-,x_+,\eps),\qquad
z_-=(x_-,y_-),\quad z_+=(x_+,y_+).
$$
For any path $\bdk=(\kappa_i)$ in the graph $\Gamma$, critical
points of the functional
$$
A_\bdk(\bdz)=\sum L_{\kappa_i}(z_i,z_{i+1},\eps),\qquad z_i \in
U_{\kappa_i}^-\times U_{\kappa_i}^+.
$$
correspond to trajectories with energy $H=E$ which shadow the
collision chain $(\gamma_{\kappa_i})$.

 As in Section \ref{sec:Turayev},  we replace $L_\kappa$ by a gauge
equivalent anti-integrable Lagrangian
$$
\hat
L_\kappa(z_-,z_+,\eps)=L(z_-,z_+,\eps)+S_{\kappa^-}^-(x_-)-S_{\kappa^+}^-(x_+)=R_\kappa(z_-)+O(\eps).
$$
By Lemma \ref{lem:nondeg}, $R_\kappa$ has a nondegenerate critical
point. Now Theorem \ref{thm:centers} follows from Theorem
\ref{thm:anti}. \qed

\medskip

For a concrete example, consider the spatial circular restricted 3
body problem (Sun, Jupiter, and Asteroid) and suppose that the
mass $\eps$ of Jupiter is small with respect to the mass $1-\eps$
of the Sun. The center of mass is stationary and the first two
masses move in circular orbits about it, having separation and
angular frequency both normalized to $1$.

Consider the motion of the Asteroid in the frame $Oxyz$ rotating
anti-clockwise about the $z$-axis through the Sun  at $O=(0,0,0)$.
Jupiter can be chosen at $P=(1,0,0)$. The motion of the Asteroid
$q=(x,y,z)$ is described by a Lagrangian system   of the form
(\ref{eq:Leps}), where
\begin{eqnarray*}
&  L_0(q,\dot q)
 = \frac12|\dot q|^2 + x\dot y - y\dot x +
 \frac12|q|^2+\frac1{|q|},\qquad V(q)
 = \frac1{|q|}-\frac1{|q-P|}+x.
\end{eqnarray*}
We have $M=\mR^3\setminus\{O\}$ and the singular set consists of
one point $P$.  The energy integral in the rotating coordinate
frame
$$
H=  \frac 12 |\dot q|^2-\frac 12 |q|^2
-\frac{1-\eps}{|q|}-\frac{\eps}{|q-p|}+\eps x
$$
is called the Jacobi integral, and $C=-2H$  is called the {\em
Jacobi constant}.

For $\eps=0$, the limit system is the Kepler problem of
Sun--Asteroid. Its bounded orbits
 are transformations to the rotating frame of
ellipses with parameters $a,e,\iota$, where $a$ is the semi-major
axis, $e$ is the eccentricity, and $\iota$ the inclination of the
orbit to the plane of the orbit of the Sun and Jupiter. They have
angular frequency $\Omega = a^{-3/2}$ and Jacobi constant $C =
a^{-1} + 2\sqrt{a(1-e^2)} \cos{\iota}$.  Collision orbits are
rotating Kepler arcs starting and ending at $P$.

Given  $C \in (-2,+3)$ we define the set $A_C$ of  allowed
frequencies of Kepler ellipses to be
\begin{itemize}
\item $(0,1)$ if $C \in [-1,+2]$, \item $(0,(2+C)^{3/2})$ if $C
\in (-2,-1)$, \item $((3-C)^{3/2},1)$ if $C\in (2,3)$.
\end{itemize}
The next result was proved in \cite{Bol-Mac:spatial}

\begin{theo}\label{th:3body}
For any $C \in (-2,+3)$ there exists a  subset $S\subset
A_C$  such that for any $\Omega\in S$ there
is a nondegenerate collision orbit $\gamma_\Omega$ with frequency
$\Omega$ and inclination $\iota=\arccos{
C/2-\Omega^{2/3}}$.
\end{theo}

Then Theorem \ref{thm:centers} implies:

\begin{cor}
For any finite set $\Lambda\subset S$ there exists $\eps_0
> 0$ such that for any sequence $(\Omega_n)_{n \in \mZ}$ in
$\Lambda$ and $\eps \in (0, \eps_0)$ there is a trajectory of the
spatial circular restricted 3 body problem with Jacobi constant
$C$ which   $O(\eps)$-shadows
a concatenation of collision orbits formed from Kepler arcs
$\gamma_{\Omega_n}$.  The resulting invariant set is uniformly
hyperbolic.
\end{cor}

Trajectories of the 3 body problem which shadow chains of
collision orbits of the Kepler problem were named by Poincar\'e
 second species solutions.

In \cite{Bol:Nonlin} and \cite{Bol-Neg:cel} similar results were obtained for the elliptic restricted and
nonrestricted plane 3 body problem with 2 masses small. The proof
is also based on a reduction to a DLS, but the Lagrangian is only partly
anti-integrable. The
corresponding generalizations of  Theorem \ref{thm:anti} are proved
in \cite{Bol:DCDS,Bol-Neg:RCD}.

\subsection{Lagrangian systems with slow time dependence}
\label{ssec:ail_cls_slow}

A straightforward  analog of an anti-integrable  DLS is a
continuous Lagrangian system with slow  time dependence:
\begin{equation}
\label{eq:Lf} L=L(q,\eps\dot q,t,\eps),\qquad q\in M,
\end{equation}
where $\eps>0$ is a small parameter. For example
\begin{equation}
\label{eq:La}
  L = \frac{\eps^2}2\|\dot q\|^2+U(q,t),
\end{equation}
where $\|\ \|$ is a Riemannian metric on $M$, possibly depending
on $t$. It looks similar to the anti-integrable discrete
Lagrangian (\ref{light}). Thus we may expect that the limit
$\eps\to0$  is similar to the anti-integrable limit in  DLS.

Introducing the fast time $s=t/\eps$, we obtain a system depending
slowly on $s$:
\begin{equation}
\label{eq:Ls} L=L(q,q',t,\eps),\qquad  q'=\frac{dq}{ds}.
\end{equation}
If $L$ satisfies the Legendre condition, the Lagrangian system can
be represented as a Hamiltonian system:
\begin{equation}
\label{eq:Ham} q'= \partial_p H,\quad p'=-\partial_qH,\quad
t'=\eps,
\end{equation}
where \begin{equation} \label{eq:p}
p=\partial_{q'}L(q,q',t,\eps),\quad H(q,p,t,\eps)=\langle
p,q'\rangle -L.
\end{equation}
 Denoting $z=(q,p)$, we obtain the
differential equation of  the form
\begin{equation}
z'=v(z,t,\eps),\quad t'=\eps, \label{eq:v}
\end{equation}
This is the standard form of a singularly perturbed differential
equation \cite{Vas-But}. Next we present a simplified version of
some results from \cite{Bol-Mac:anti}. References to earlier
classical results can be found in this paper. Similar anti-integrable approach methods was used in \cite{Pif}
for the Mather acceleration problem. Interesting results
on the chaotic energy growth for systems of type (\ref{eq:Lf})
were obtained recently in \cite{Gel-Tur} by a very different approach.

For $\eps=0$ the frozen Lagrangian (\ref{eq:Ls}) takes the form
\begin{equation}
\label{eq:L0} L_t(q,q')=L(q,q',t,0),
\end{equation}
where  the time $t$ is now frozen. Hence the energy
$H_t(q,p)=H(q,p,t,0)$ is a first integral of the frozen system.

Suppose that for each $t$ the frozen system has a hyperbolic
equilibrium point $O_t$, and there exists a homoclinic trajectory
$\gamma_t:\mR\to M$ with $\gamma_t(\pm\infty)=O_t$. Without loss
of generality we may assume that $H_t=0$ at the equilibrium $O_t$.

\begin{rem}  For a natural system
 (\ref{eq:La}) on a compact manifold the
homoclinic trajectory always exists if $O_t$ is a point of strict
nondegenerate minimum of $U(q,t)$, see  \cite{Bol:libr}. If
$M$ is not simply connected, there are at least as many of them as
the number of generators of the semigroup $\pi_1(M)$.
\end{rem}

Let $I$ be the open set of $t$ such that the homoclinic $\gamma_t$
is transverse: the intersection of the stable and unstable
manifolds $W_t^{\pm}$ along $\gamma_t$ is transversal in the
energy level $\{H_t=0\}$. The Maupertuis action
$$
  f(t)=J(\gamma_{t})=\int_{\gamma_{t}}p\,dq
$$
is a smooth function on $I$.
 Hamilton's principle implies that
$$
f'(t)=-\int_{-\infty}^\infty\left(\partial_tH_t\right)
\big|_{\gamma_{t}(s)}\,ds,
$$
where the homoclinic orbit $\gamma_{t}$ is parameterized by the
fast time $s$. Hence  $f$ is a version of the Poincar\'e--Melnikov
function.

For simplicity assume that $L$ is 1-periodic in time. Take a
finite set of  nondegenerate critical points of $f$ in
$\mT=\mR/\mZ$. The corresponding points in $\mR$ form a discrete
periodic set $K\subset\mR$. If $\rho>0$ is small enough, then
$I_k=(k-\rho,k+\rho)$, $k\in K$, are nonintersecting intervals.

\begin{theo}
\label{thm:cont} Fix small $\delta>0$ and any $0<c_1<c_2$. Suppose
that $\eps>0$ is sufficiently small. Then for any increasing
sequence $k_i\in K$ there exists a unique trajectory $q(t)$ and
sequences $t_i\in I_{k_i}$, $c_1<T_i<c_2$, such that
\begin{itemize}
\item $d(q(t),\gamma_{k_i}(\mR)\le \delta$ for $t\in[t_i,t_i+\eps
T_i]$; \item $d(q(t),O_t )\le\delta$  for $t\in [t_i+ \eps T_i,
t_{i+1}]$. \item Moreover
 $|t_i-k_i|\le C\eps$ and
$d(q(t),\gamma_{k_i}(\mR))\le C\eps$ for $t\in[t_i,t_i+\eps T_i]$.
\end{itemize}
\end{theo}

\begin{figure}[ht]
\begin{center}
\includegraphics{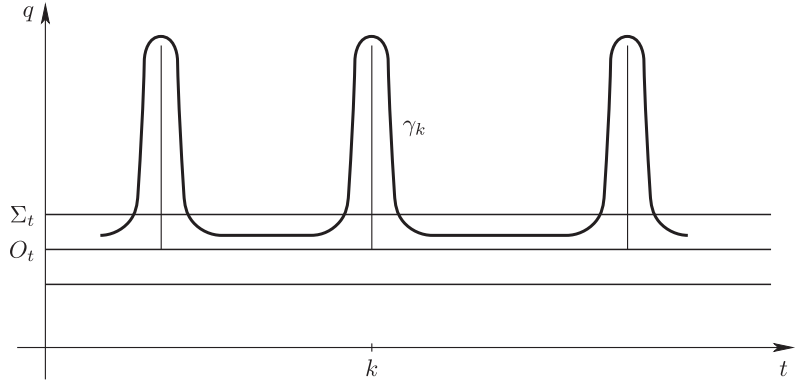}
\caption{Multibump trajectory shadowing the homoclinic
$\gamma_k$.} \label{fig:multibump}
\end{center}
\end{figure}

We call $q(t)$ a multibump trajectory, see Figure
\ref{fig:multibump}.
 Theorem \ref{thm:cont} holds also for
nonperiodic $L$ if certain uniformity is satisfied.

In \cite{Bol-Mac:anti} a similar result was proved without the
transversality assumption. Consider system (\ref{eq:La}) on a non
simply connected compact manifold and let $f(t)$ be the minimum of
actions of  noncontactible homoclinics  for the frozen system.
Then it is enough to assume that $f$ is nonconstant.

\medskip

\noindent{\it Proof.} First let $\eps=0$. Since $O_t $ is a
hyperbolic equilibrium of the limit system, it has local stable
and unstable Lagrangian manifolds $W_t^+$ and $W_t^-$ in the phase
space $T^*M$.  We assume for simplicity that the projection
$W_t^{\pm}\to M$ is nondegenerate at $O_t$. This holds if the
limit system is natural as in (\ref{eq:La}); in general we can
achieve this by changing symplectic coordinates near $O_t $. Hence
\begin{eqnarray*}
\label{eq:stable} W_t^+=\{(q,p):q\in D_t,\;p=-\nabla
S_t^+(q)\},\quad W_t^-=\{(q,p):q\in D_t,\;p=\nabla  S_t^-(q)\},
\end{eqnarray*}
where $D_t$ is a small $\delta$-neighborhood of $O_t $.
 Since $S_t^\pm$ is defined up to a function
of time, without loss of generality we may assume that
$S_t^\pm(O_t )=0$. Then $-S_t^+(q)$ is the action of the
trajectory of the frozen system starting at $q$ and asymptotic to
$O_t$  as the fast time $s\to+\infty$, while $S_t^-$ is the action
of the trajectory asymptotic to $O_t$  as $s\to-\infty$ and ending
at $q$.

Take small $\eps>0$. In the next lemma the slow time $t$ is used.

\begin{lem}
\label{lem:bound} Fix $0<c_1<c_2$. For sufficiently small
$\eps>0$, any $a<b$ such that $b-a\in (c_1,c_2)$ and any $x\in
D_{a}$, $y\in D_{b}$, there exists a unique trajectory $q(t)\in
D_t$, $a\le t\le b$,   such that $q(a)=x$ and $q(b)=y$. Moreover:
\begin{itemize}
\item the function $q(t)=q(t,a,b,x,y,\eps)$ is  smooth for
$\eps>0$. \item The trajectory $q(t)$ is $C^0$ close to the
concatenation of the asymptotic trajectories of the frozen system.
\item Let
$$
      S(a,b,x,y,\eps)
  =  \eps^{-1}\int_{a}^{b}L(q(t),\eps\dot q(t),t,\eps)\,dt
$$
 be the action of the trajectory $q(t)$  normalized to the fast time. Then
\begin{equation}
\label{eq:action}
 S =S^+(x)+S^-(y)+\eps v(a,b,x,y,\eps),
\end{equation}
where $v$ is a smooth function for $\eps>0$ and
  $\|v\|_{C^2(D_a\times D_b)}\le C$ for some  constant $C$ independent of $\eps$.
\end{itemize}
\end{lem}

Note that $q(t)$ behaves badly as $\eps\to 0$: its time derivative
 is unbounded of order $\eps^{-1}$, while the generating
function $S$ is regular as $\eps\to 0$. If  the Lagrangian
(\ref{eq:Lf})  is independent of $t$ and $\eps$, then
Lemma~\ref{lem:bound} follows from Shilnikov's lemma \cite{Shil},
or the strong $\lambda$-lemma \cite{Deng}.

In the general case, the existence of the solution of the given
boundary value problem can be obtained from the results of the
theory of singularly perturbed differential equations
\cite{Vas-But}.

\medskip

Let $p(t)=\partial_{q'}L(q(t),\eps \dot q(t),t,\eps)$ be the
momentum (\ref{eq:p}) of the trajectory $q(t)$. By the first
variation formula,
\begin{equation}
\label{eq:var1}
\begin{array}{c}
\partial_x S=-p(a),\quad \partial_yS=p(b),\\
\partial_{a}S=H(x,p(a),a,\eps),\quad \partial_{b}
S=-H(y,p(b),b,\eps).
\end{array}
\end{equation}

Lemma \ref{lem:bound} describes trajectories which stay near the
equilibrium $O_t$ during a finite interval of the slow time $t$.
Next we describe trajectories which travel near the homoclinic
orbit $\gamma_k$ in a short interval of slow time of order
$\eps$.  Hence the fast time $s=\eps^{-1}t$ will be used.

To simplify notation, we can assume that   $D_t=D_{k}$ is
independent of $t$ for $t\in I_k$. Let $\Sigma_k=\partial D_{k}$.
Let $\gamma_{k}(0),\gamma_{k}(\tau_k)$ be the intersection points
of the transverse homoclinic $\gamma_{k}:\mR\to M$ with
$\Sigma_k$. Changing $\delta$ if needed, we may assume that  they
are nonconjugate along $\gamma_k$. Then for $z=(x,y,t,T)$ close to
$z_k^0=(\gamma_{k}(0),\gamma_{k}(\tau_k),k,\tau_k)$ and
sufficiently small $\eps>0$, there exists a trajectory
$\beta_\eps(s)$, $0\le s\le T$, for system with the Lagrangian
$L(q,q',t+\eps s,\eps)$  with boundary conditions
$\beta_\eps(0)=x$, $\beta_\eps(T)=y$, which is close to the
homoclinic orbit $\gamma_{k}$. Let
$$
    \Phi_k(z,\eps)
  = \int_0^T L(\beta_\eps(s),\beta'_\eps(s),t+\eps s,\eps)\,ds
  $$
be its action. The derivatives of $\Phi_k$ satisfy the equations
similar to (\ref{eq:var1}):
\begin{equation}
\label{eq:var2}
\begin{array}{c}
\partial_x \Phi_k=-p(0),\quad \partial_yF_k=p(T),\\
\partial_{t}\Phi_k=H(x,p(0),t,\eps),\quad \partial_{T}
\Phi_k=-H(y,p(T),t+\eps T,\eps),
\end{array}
\end{equation}
where $p(s)$ is the momentum of $\beta_\eps(s)$.

\begin{lem}The function
$$
R_k(z)=S_t^-(x)+\Phi_k(x,y,t,T,0)+S_t^+(y),\qquad
    x,y\in \Sigma_k,
$$
has a nondegenerate critical point $z_k^0$.
\end{lem}

Indeed,  $R_k$ is the action of the concatenation of an asymptotic
trajectory of the frozen system starting at $O_t$ and ending at
$x$, trajectory $\beta_0$  joining $x$ and $y$ in a neighborhood
of $\gamma_k$, and a trajectory starting at $y$ and asymptotic to
$O_t$. The equation $\partial_TR_k=0$ implies that $H_t=0$ along
$q(s)$. Then equations $\partial_xR_k=0$ and $\partial_yR_k=0$
imply that the concatenation is   a smooth trajectory of the
frozen system homoclinic to $O_t$. Thus it coincides with
$\gamma_t$. Finally $\partial_tR_k=0$ means that $t$ is a critical
point of $f(t)$. \qed

\medskip

Next we define the DLS describing multibump trajectories. Let
$J=\{\kappa=(\kappa_-,\kappa_+)\in K^2:\kappa_-<\kappa_+\}$ be the
set of vertices of the graph. We join $\kappa,\kappa'\in J$ with
an edge if $\kappa_+=\kappa'_-$. Define the discrete Lagrangian
with $2m$ degrees of freedom by
$$
L_\kappa(z_-,z_+,\eps)=\Phi_k(z_-,\eps)+S(t_-+\eps
T_-,t_+,y_-,x_+,\eps),\qquad z_\pm=(x_\pm,y_\pm,t_\pm,T_\pm).
$$
Thus $L_\kappa(z_-,z_+,\eps)$ is the action of the broken
trajectory which starts at $x_-\in\Sigma_{\kappa_-}$ at time
$t_-$, travels near the homoclinic trajectory $\gamma_{\kappa_-}$
till it reaches $\Sigma_{\kappa_+}$ at the point $y_-$ at time
$t+\eps T_-$, then travels close to the hyperbolic equilibrium
$O_t$ for $t_-+\eps T_-\le t\le t_+$ and ends at
$x_+\in\Sigma_{\kappa_+}$ at time $t_+$.

By (\ref{eq:var1})--(\ref{eq:var2}), critical points of  the
functional
\begin{eqnarray*}
      A_\bdk(z)
  = \sum L_{\kappa_i}(z_i,z_{i+1},\eps) ,\qquad
  \kappa_i=(k_i,k_{i+1}),
\end{eqnarray*}
where
$$
z_i=(x_i,y_i,t_i,T_i),\qquad x_i,y_i\in\Sigma_{k_i},\quad
  t_i\in I_{k_i}, \quad T_i>0,
$$
correspond to trajectories shadowing the chain of homoclinics
$(\gamma_{k_i})$.

We make a gauge transformation replacing $L_\kappa$ by the
Lagrangian
$$
\hat L_\kappa(z_-,z_+,\eps)=
L_\kappa(z_-,z_+,\eps)+S_{\kappa^-}(x_-)-S_{\kappa^+}(x_+)
=R_{\kappa^-}(z_-)+O(\eps),
$$
where $R_{\kappa^-}$ has a nondegenerate critical point
$z_{\kappa^-}^0$. Then $\hat L_\kappa$ has an anti-integrable
form.

Note that the graph $\Gamma$ describing the anti-integrable system
is infinite. However, the Lagrangian is invariant with respect to
$\mZ$-action on $\mR$, so uniform anti-integrability  in Theorem
\ref{thm:anti_inf} holds.
Theorem \ref{thm:cont} is proved.

\section{Separatrix map}

\subsection{AI limit in Zaslavsky separatrix map}
\label{sec:sm}

Consider an integrable area-preserving map $F_0$ having a hyperbolic
fixed point with two homoclinic separatrix loops. Let $F_\eps$ be
a perturbed map also assumed to be area-preserving. In an
$\eps$-neighborhood of the unperturbed separatrix loops the
dynamics of $F_\eps$ is determined by the separatrix map.

\begin{figure}[ht]
\begin{center}
\includegraphics{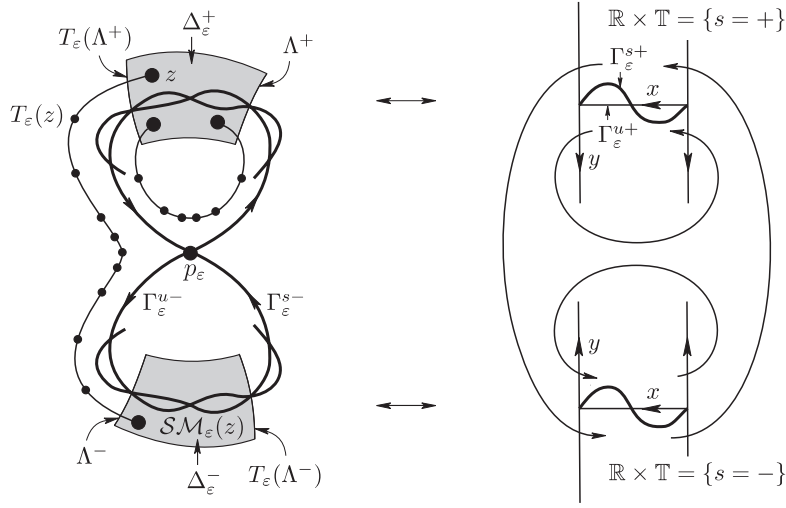}
\caption{The separatrix map.} \label{fig:sep}
\end{center}
\end{figure}

The construction is presented in Fig \ref{fig:sep}. The left-hand
side of the figure presents the phase space of the map $F_\eps$.
We see the hyperbolic fixed point $p_\eps$ and its asymptotic
curves (separatrices), splitted for $\eps\ne 0$. We also see two
grey domains $\Delta_\eps^\pm$ on which the separatrix map will be
defined. Boundaries of these domains are curvilinear quadrangles.
The ``horizontal'' sides can be regarded as lying on invariant KAM
curves (this is convenient, but not necessary) while the ``vertical''
sides for each quadrangle are images of each other under the maps
$F_\eps$ and $F_\eps^{-1}$. For any point $z\in\Delta_\eps =
\Delta_\eps^+\cup\Delta_\eps^-$ its image $F_\eps(z)$ lies outside
$\Delta_\eps$. By definition the image of $z$ under the separatrix
map is $F_\eps^n(z)$, where $n = n(\eps)$ is the minimal natural
number such that $F_\eps^n(z)\in\Delta_\eps$. In Figure
\ref{fig:sep} two such points $z$ and their images are presented.

In some convenient coordinates the separatrix map can be computed in the form: an explicit part plus small error terms,
\begin{equation}
\label{sm_model}
    \left(\begin{array}{c} y \\ x \\ \sigma \end{array}\right)
  \mapsto
    \left(\begin{array}{c} y_+ \\ x_+ \\ \sigma_+ \end{array}\right),
  \qquad
  \begin{array}{rcl}
     y_+ &=& y + \lambda \partial V_\sigma/\partial x
               + O(\eps), \\
     x_+ &=& x + \frac {1+O(\eps)}\lambda (\omega_\sigma + \log|y_+|) , \\
  \sigma_+ &=& \sigma\, \mbox{sign\,} y_+, \qquad
  \sigma,\sigma_+\in\{-1,1\}.
  \end{array}
\end{equation}
Here the variable $x\bmod 1$ changes along the separatrix, $y$
changes across the separatrix, and the discrete variable
$\sigma=\pm1$ indicates near which loop the orbit is located at
the given time moment. The functions $V_{\sigma}(x)$ (the
Poincar\'e-Melnikov potentials) are periodic with period 1.
Details can be found for example, in \cite{PT,TZ}).

The map (\ref{sm_model}) has the form
\begin{equation}
\label{sm_H}
       y = \partial W / \partial x,\quad
     x_+ = \partial W / \partial y_+,\quad
  \sigma_+ = \sigma\,\mbox{sign\,} y_+,
\end{equation}
where the generating function
$$
  W = W(y_+,x,\sigma)
    = xy_+ - \lambda V_\sigma(x)
           + \frac{1+O(\eps)}{\lambda}
                \big( \omega_\sigma + \log |y_+| - 1 \big)\, y_+ .
$$
The first two equations (\ref{sm_H}) can be written as follows:
$$
  x_+ dy_+  +  y dx = dW(y_+,x,\sigma).
$$

To present (\ref{sm_model}) in the Lagrangian form, consider
another generating function (the Legendre transform of $W$)
\begin{equation}
\label{form_Lagr}
  \lambda L = x_+ y_+ - W, \qquad
  y_+ dx_+ - ydx = \lambda dL.
\end{equation}
Here we use the fact that a Lagrangian is defined up to a nonzero
constant multiplier. It is easy to obtain an explicit formula for
$L$:
\begin{eqnarray}
\label{L_sm} &      L(x,x_+,\sigma,\vartheta_+)
    =  (1 + O(\eps)) \vartheta_+
                e^{\lambda(x_+ - x - \hat\omega_\sigma)}
      + V_\sigma (x), & \\
\nonumber &   \vartheta_+ = \mbox{sign\,} y_+, \quad
    \hat\omega_\sigma = \omega_\sigma + \lambda^{-1} \log\lambda^2.
\end{eqnarray}
Thus the separatrix map is a $\mZ$-equivariant DLS, where $\mZ$
acts on $M = \mR\times\mZ_2\times\mZ_2$ by the shifts
$$
    M\ni    (x,\sigma,\thet)
   \mapsto k(x,\sigma,\thet)
      =    (x + k,\sigma,\thet).
$$

The Lagrangian form of the separatrix map is as follows:
\begin{equation}
\label{eq:sm_Lagr}
  \left(\begin{array}{c} x_- \\ x \\ \sigma_- \\ \vartheta
  \end{array}\right)
  \mapsto
  \left(\begin{array}{c} x \\ x_+ \\ \sigma \\ \vartheta_+
  \end{array}\right) ,\quad
  \begin{array}{l}
  \sigma = \sigma_- \vartheta, \\
  \vartheta_+ = \sigma\sigma_+, \\
  \displaystyle{
  \frac{\partial}{\partial x} (L_-  +  L) = 0, }
  \end{array}
\end{equation}
where
$$
   L_- = L(x_-,x,\sigma_-,\vartheta),\quad
   L  = L(x,x_+,\sigma,\vartheta_+).
$$
The first two equations (\ref{eq:sm_Lagr}) are generated by the
definition of $\vartheta$ $(\vartheta = \mbox{sign}\,I)$ and by
the last equation (\ref{sm_model}). The third equation
(\ref{eq:sm_Lagr}) follows from (\ref{form_Lagr}) since, according
to (\ref{form_Lagr}) and to the analogous equation $y dx - y_-
dx_- = dL_-$, we have
$$
  y = - \partial L / \partial x
    =   \partial L_- / \partial x.
$$
It is easy to check that the quantities $x_+$, $\sigma$ and
$\vartheta_+$ are computed uniquely from (\ref{eq:sm_Lagr}) in
terms of $x_-,x,\sigma_-$, and $\vartheta$.

We obtain a DLS with extra discrete variables $\sigma$ and
$\thet$. Any sequence
$$
   \bdz = \{z_j\},\quad
    z_j = \left(\begin{array}{c} x_j\\ \sigma_j\\ \vartheta_j
                 \end{array}\right),\qquad
   \sigma_{j+1} = \sigma_j \vartheta_{j+1}
$$
is said to be a path. Let $\Sigma$ be the set of all paths.

In general the index $j$ takes all integer values. However, it is
possible to consider also semifinite and finite paths. Paths
finite from the left begin with a triple $f_j$, where $x_j =
+\infty$. Paths finite from the right end with $f_j$, where $x_j =
-\infty$. Paths finite from the left and from the right are called
finite.

The action $A$ is defined as the formal sum
$$
  A = A(f) = \sum_j L(x_j,x_{j+1},\sigma_j,\thet_{j+1}).
$$
The path $f^0$ is said to be an extremal (or a trajectory) iff $\partial A /
\partial x_j |_{f = f^0} = 0$ for any $j$.

Note that semifinite trajectories belong to separatrices. Finite
ones belong to both stable and unstable separatrices, and
therefore they are homoclinic trajectories.

We define the distance $\rho$ on $\Sigma$ in the following way.
Let $\bdz'$ and $\bdz''$ be paths, where
$$
 \bdz'_j = \left(\begin{array}{c} x'_j\\ \sigma'_j\\ \vartheta'_j
                  \end{array}\right),\quad
 \bdz''_j = \left(\begin{array}{c} x''_j\\ \sigma''_j\\ \vartheta''_j
                  \end{array}\right).
$$
We put $\rho(\bdz',\bdz'') = \infty$ if the sequences
$\sigma'_j,\thet'_j$ do not coincide with $\sigma''_j,\thet''_j$
or if for some $j$ only one of the triples is defined. Otherwise
we put
$$
  \rho(\bdz',\bdz'') = \sup_j |x'_j - x''_j|.
$$
Here we put $|-\infty - (-\infty)| = |+\infty - (+\infty)| = 0$.

Let $\mbox{Cr}(\sigma)$ denote a finite set of nondegenerate
critical points of the function $V_\sigma$. The set $\Sigma$
contains the subset $\Pi$ of simple paths (codes). By the
definition a path $\bdz$ is simple if $x_j\in\mbox{Cr}(\sigma_j)$
for any $j$.

\begin{theo}
\label{theo:sep+anti} Suppose that the the constants $c_1$, $c_2 =
c_2(c_1)$ are sufficiently large. Then, for any simple path
$\bdz^*$ such that $x^*_j - x^*_{j+1} > c_1$ for all $j$, there
exists a unique trajectory $\tilde\bdz$ in the
$c_2^{-1}$-neighborhood of $\bdz^*$. The trajectory $\tilde\bdz$
is hyperbolic.
\end{theo}

Theorem \ref{theo:sep+anti} establishes a symbolic dynamics in a
neighborhood of separatrices of an area-preserving map. Theorem
\ref{theo:sep+anti}  can be deduced from Theorem \ref{thm:anti} by
introducing the DLS with the Lagrangians
$$
L_\kappa(x_-,x_+)=L(x_-+k,x_+,\sigma,\theta),\qquad
\kappa=(k,\sigma,\theta),\quad x_\pm\in (0,2\pi),
$$
where $k>c_1$ is large enough.  The corresponding graph has vertices $\kappa$, and vertices $\kappa$, $\kappa'$ are joined by an edge if $\sigma'=\sigma\theta$.  There is an edge $\gamma$ for every nondegenerate  critical point of $V_\sigma$ in $(0,2\pi)$. For large $k$, the DLS is anti-integrable. The graph $\Gamma$ is infinite but condition {\bf U} of uniform anti-integrability evidently holds.
\medskip

The traditional approach to symbolic dynamics near separatrices is
presented in \cite{birk,Alex,Nitec}. Another version of the
separatrix map is constructed by Shilnikov and Afraimovich,
\cite{AS}, see also \cite{SSTC}. This construction is also given
in Section \ref{sec:Turayev}.

\subsection{Separatrix map and Arnold diffusion}
\label{sec:sm_mult}

Ideas of AI limit can be applied in the problem of
Arnold diffusion. Here we only discuss the {\it a priori} unstable
case where unlike the original ({\it a priori stable} situation,
\cite{Arn64}) the unperturbed integrable system contains normally
hyperbolic manifold $N$. In the perturbed system chaos is
essentially concentrated near $N$ and its asymptotic manifolds.

We consider a nonautonomous near-integrable Hamiltonian system on
the phase space $\mT^n_x\times\overline\calD\times D\times\mT_t$,
where $\calD\subset\mR^n_y$, is an open domain with compact
closure $\overline\calD$, $D\subset \mR^2_{(v,u)}$ is an open
domain. The Hamiltonian function and the symplectic structure are
as follows:
\begin{eqnarray*}
\nonumber &   H(y,x,v,u,t,\eps)
  = H_0(y,v,u) + \eps H_1(y,x,v,u,t) + \eps^2 H_2(y,x,v,u,t,\eps), & \\
\nonumber &   \omega
  = dy\wedge dx + dv\wedge du. &
\end{eqnarray*}
As usual $\eps\ge 0$ is a small parameter. Hamiltonian equations
have the form
\begin{equation}
\label{eq:ham}
  \dot y = - \partial H / \partial x,\quad
  \dot x =  \partial H / \partial y,\quad
  \dot v = - \partial H / \partial u,\quad
  \dot u =  \partial H / \partial v.
\end{equation}

We assume that in the unperturbed Hamiltonian the variables $y$
are separated from $u$ and $v$ i.e., $H_0(y,v,u) = F(y,f(v,u))$.
The function $f$ has a nondegenerate saddle point $(v,u) = (0,0)$,
a unique critical point on a compact connected component of the
set
$$
  \gamma = \{(v,u)\in D : f(v,u) = f(0,0)\} .
$$
In dynamical terminology $(0,0)$ is a hyperbolic equilibrium of
the Hamiltonian system $(D,dv\wedge du,f)$ with one degree of
freedom, and the corresponding separatrices $\gamma$ are doubled.
Topologically these separatrices form a figure-eight: two loops,
$\gamma^\pm$, issuing from one point, $\gamma = \gamma^+ \cup
\gamma^-$.

Hence the unperturbed normally hyperbolic manifold is
$$
  N = \mT^n_x\times \overline\calD\times (0,0)\times\mT_t.
$$
It is foliated in tori
$$
  N_y = \{x,y,u,v,t) : u=v=0, \; y = \mbox{const}\}
$$
carrying quasi-periodic dynamics with frequencies
$$
  \Big(\begin{array}{c}\nu(y)\\ 1\end{array}\Big), \qquad
  \nu(y) = \frac{\partial H_0}{\partial y}(y,0,0).
$$
We are interested in the perturbed dynamics near asymptotic
manifolds of $N$.

The problem of Arnold diffusion in a priori unstable case has
three aspects, contained in the following
\smallskip

\begin{conj} \cite{AKN,Tre2012}.
\label{conj:diffu}

{\bf A} (genericity). Diffusion exists for an open dense set of
$C^r$-perturbations, where $r\in\mN\cup\{\infty,\omega\}$ is sufficiently large.

{\bf B} (freedom). Projection of a diffusion trajectory to the
$y$-space can go in a small neighborhood of any smooth curve
$\chi\subset\calD$.

{\bf C} (velocity). There are ``fast'' diffusion trajectories
whose average velocity along $\chi$ is of order $\eps/\log\eps$.
\end{conj}

There are several different approaches to the problem. Traditional
approach is based on the construction of transition chains of
hyperbolic tori \cite{Arn64,CG,DLS2,Fontich,GGM}, later it was
effectively supplemented with the idea of the scattering map
\cite{DLS_scat,Delsh-Hug1,Delsh-Hug2} and symbolic dynamics in
polysystems \cite{Bounemoura}. Variational approach was developed
in \cite{Bessi,Berti,BeBo,Cheng-diff1,Cheng-diff2,GR,Kalo-Levi,Kal-Zhang}.

Conjecture \ref{conj:diffu} is proven only in the case $n=1$ (2
and a half degrees of freedom), \cite{Tre2004}. Here we explain
briefly ideas and methods of \cite{Tre2004}. The system
(\ref{eq:ham})$|_{\eps=0}$ has the $n$-parametric family of
(partially) hyperbolic $(n+1)$-tori $N_y$ which foliate the normally hyperbolic manifold $N$. Asymptotic manifolds of $N$ form two components
$\{y\}\times\mT_x^n\times\gamma^\pm\times\mT_t$. Hence, after
perturbation the situation reminds the one discussed in the
previous section. In the phase space of the time-one map a picture
analogous to that presented in Figure \ref{fig:sep} appears.
However unlike the case considered in Section \ref{sec:sm} now we
have to deal with a family of hyperbolic tori and their asymptotic
manifolds. The separatrix map still can be defined
\cite{Tre_sm_res} and explicit formulas for it can regarded as a
multidimensional generalization of (\ref{sm_model}).

The separatrix map can be presented in the form
\begin{eqnarray}
\nonumber &        (\zeta,\rho,\tau_-,\tau,\sigma_-,\theta)
 \mapsto (\zeta_+,\rho_+,\tau,\tau_+,\sigma,\theta_+), & \\
\label{sm_mult} &    \rho = \frac{\partial\calR}{\partial\zeta},
\quad
  \zeta_+ = \frac{\partial\calR}{\partial\rho_+}, \quad
  \frac{\partial}{\partial\tau} (\calR_- + \calR) = 0, \quad
  \sigma = \sigma_- \thet, \quad
  \thet_+ = \sigma\sigma_+ , & \\
\nonumber & \calR =
R(\zeta,\rho_+,\tau,\tau_+,\sigma,\theta_+,t_+,\eps), \quad
  \calR_- = R(\zeta_-,\rho,\tau_-,\tau,\sigma_-,\theta,t,\eps).
&
\end{eqnarray}
(compare with (\ref{eq:sm_Lagr})). Here up to small error terms
$\eps\rho=y$, $\zeta - x$ is a function of $y,u,v$. The integer
vatiable $t$ has a meaning of time the trajectory of the time-one
map travels outside analogs of the domains $\Delta_\eps^\pm$, see
Fig. \ref{fig:sep}. The variables $\sigma$ and $\thet$ are
analogous to that from Section \ref{sec:sm}, and
\begin{eqnarray*}
     \calR
 &=& \langle \rho_+,\zeta+\nu t_+\rangle
     - (\tau_+ - \tau - t_+) \bdH(\eps\rho_+,\zeta) \\
 &&\qquad\qquad\qquad
     - \thet_+ e^{\lambda (\tau_+ - \tau - \omega^\sigma_+)}
     + \hat\Theta^\sigma(\eps\rho_+,\zeta,\tau), \\
 &&\!\!\!\!\!\!\!\!\!\!\!\!\!\!\!\!\!\!\!\!
    \omega_+^\sigma = t_+ + \lambda^{-1}\log\eps + f(\eps\rho_+), \quad
    \nu = \nu(\eps\rho_+), \quad
    \lambda = \lambda(\eps\rho_+).
\end{eqnarray*}
If $n=0$, we do not have variables $\rho$ and $\zeta$, and only the last two terms in $\calR$ remain. In this case up to a constant multiplier $\calR$ converts into $L$, see (\ref{L_sm}).

We are interested in the case when the quantities
$t_+ + \lambda^{-1}\log\eps$ are greater than a big positive constant
$K_0$. So, $e^{K_0}$ plays the role of a large parameter in the AI
limit. The function $\bdH$ is essential only in small
neighborhoods of strong (low order) resonances
$$
  \{ \rho : \langle \nu(\rho),k \rangle + k_0 = 0 \}, \qquad
  k\in\mZ^n,\; k_0\in\mZ,\quad
  |k| + |k_0| < C.
$$

We see that the dynamical system (\ref{sm_mult}) is partially
Hamiltonian (w.r.t. $\rho,\zeta$) and partially Lagrangian (w.r.t.
$\tau$). We can represent dynamics of the separatrix map by a DLS
as follows.

 Because of the presence of new Hamiltonian variables we have to use a certain generalization of the method of AI limit.
 Unfortunately, in this generalization the symbolic dynamics is not as clean and standard as the one given in Sections
 \ref{subsec:main} or \ref{sec:sm}. The construction is as follows, \cite{Tre_traj}.
 Having a finite piece of a trajectory of the separatrix map and the corresponding quasitrajectory,
 a piece of the same length (the code), we present a rule according to which the code can be extended by adding a new point.
 Then according to the main result of \cite{Tre_traj} the trajectory can be slightly deformed and extended staying
 close to the longer code. So we have again the pair: a piece of trajectory with a code.
 By using a certain freedom in the rule for the extension of a code,
 one can hope to push the trajectory in the desired direction in the $y$-space.

Here another difficulty appears. This possibility to push the
orbit in the proper direction is relatively simple in the
nonresonant zone i.e., where the frequency vector
$(\nu(y),1)\in\mR^{n+1}$ (1 is the time frequency) does not allow
resonances of low orders. In a near-resonant domain construction
of a reasonable code is a separate delicate problem, solved
completely only in the case of two-and-a-half degrees of freedom
\cite{Tre2004}. Diffusion in domains free of low order resonances
in the case of an arbitrary dimension is established in
\cite{Tre2012}.

Finally we mention another application of the separatrix map.
According to \cite{CK,KZ,KZZ} in (so far special cases of) {\it a
priori} unstable systems with 2 and a half degrees of freedom a
large set of trajectories is constructed whose projections to the
$y$-axis for small $\eps>0$ behave like trajectories of the
Brownian motion. This shows that proposed by Chirikov term
``diffusion'' is quite adequate for the phenomenon discussed
above.


\begin{thebibliography}{100}

\bibitem{AS}    Afraimovich, V.S. Shilnikov, L.P.
              Small periodic perturbations of autonomous systems. (Russian) Dokl. Akad. Nauk SSSR  214  (1974), 739--742.

\bibitem{Alex}      Alekseev~V.~M.,
              Quasi-random dynamical systems I,II,III.
              {\it Mat. Sb.} 1968, V.~76, P. 72--134,
              1968, V.~77, P. 545--601, 1969, V.~78, P. 3--50.
              (Russian) Engl. transl. in {\it Math. USSR-Sb.}
              1968, V. 5, P. 73--128,
              1968, V. 6, 1969, V. 7, P. 1--43.
              \bibitem{Arn64}     Arnold V.I.,
              Instability of dynamical systems with several degrees of freedom.
              Sov. Math. Dokl. {\bf 5}, 581--585, 1964.
\bibitem{Arnold:matmex}       Arnold~V.~I.,
              Mathematical methods in classical mechanics.
              Springer-Verlag, New-York, Heidelberg,
              Berlin 1982.
\bibitem{AKN}       Arnold~V.~I., Kozlov~V.~V., Neishtadt~A.~I.,
              Mathematical aspects of classical and celestial
              mechanics. Encyclopedia of Math. Sciences,
              Vol.~3, Springer-Verlag 1989.
\bibitem{aubry0}    Aubry~S., Abramovici~G.,
              Chaotic trajectories in the standard map: the concept
              of anti-integrability.
              {\it Physica} 43 D, 1990, 199--219.
\bibitem{Aubry-MacKay}               Aubry S., MacKay R. S. and Baesens C.,
Equivalence of uniform hyperbolicity for symplectic twist maps and
phonon gap for Frenkel-Kontorova models Physica D {\bf 56},
123-134 (1992).
\bibitem{Baesens} Baesens, C., Chen, Y.-Ch., and MacKay, R. S.,
              Abrupt Bifurcations in Chaotic Scattering: View from the Anti-Integrable Limit,
              Nonlinearity, 2013, vol.26, no. 9, pp. 2703--2730.

\bibitem{Berti}    Berti~M. and Bolle~P.
              A functional analysis approach to Arnold Diffusion, Annales Institute Henri Poincare\'e, analyse non-lineaire,
              19, 4, 2002, pp. 395-450.
\bibitem{BeBo}    Berti~M., Biasco L., and Bolle P.,
              Drift in phase space: A new variational mechanism with
              optimal diffusion time.
              J. Math. Pures Appl. (9) 82 (2003), no. 6, 613--664.
              \bibitem{Bessi}     Bessi~U.,
              An approach to Arnold's diffusion through the calculus
              of variations.
              {\it Nonlin.~Anal.~TMA} {\bf 20}, 1303--1318 (1996)
\bibitem{Bia}       Bialy M.,
              Maximizing orbits for higher-dimensional convex billiards.
              {\it J. of Modern Dynamics}, Vol. 3, No. 1, 2009, 51-59.
\bibitem{birk}      Birkhoff~G.~D.,
            Dynamical systems.
            Amer. Math. Soc. Colloquium Publications V.~9, New York, 1927.
\bibitem{Bol:libr}
Bolotin~S.V., Libration motions of natural dynamical systems.
Vestnik Moskov. Univ. Ser. I Mat. Mekh. 1978, no. 6, 72--77.
\bibitem{Bol:homtori} Bolotin~S., Symbolic dynamics near minimal hyperbolic
invariant tori of Lagrangian systems, Nonlinearity, 14:5 (2001),
1123--1140.
 \bibitem{Bol:DCDS} Bolotin~S.,
            Shadowing chains of collision orbits, Discrete Contin.
            Dyn. Syst., 14:2 (2006), 235--260.
\bibitem{Bol:Nonlin}
S.~Bolotin, Symbolic dynamics of almost collision orbits and skew
products of symplectic maps. {\it Nonlinearity}, {\bf  19} (2006),
2041--2063.

\bibitem{Bol-Mac:anti}   Bolotin~S.~V. and MacKay~R.,
            Multibump orbits near the anti-integrable limit
            for Lagrangian systems
            {\it Nonlinearity}, V.10, No 5, 1997,  1015.

\bibitem{Bol-Mac:centers}
S.~Bolotin and R.S.~MacKay, Periodic and chaotic trajectories of
the second species for the $n$-centre problem. {\it Celest.\
Mech.\ \& Dynam.\ Astron.}, {\bf 77} (2000), 49--75.


\bibitem{Bol-Mac:spatial} Bolotin~S.~V. and MacKay~R.S.,
            Nonplanar second species periodic and chaotic trajectories for the circular restricted three-body problem,
            Celestial Mech. Dynam. Astronom., 94:4 (2006), 433-449.
            \bibitem{Bol-Neg:reg} Bolotin S., Negrini P., Global regularization for the
n-center problem on a manifold, Discrete Contin. Dyn. Syst., 8:4
(2002), 873–892.
\bibitem{Bol-Neg:cel}  S. Bolotin, P. Negrini,
     Variational approach to second species periodic solutions of Poincar\'e of the three-body problem,
     Discrete Contin. Dyn. Syst., 33:3 (2013), 1009-1032.


\bibitem{Bol-Neg:RCD} Bolotin S., Negrini P., Shilnikov lemma for a nondegenerate critical manifold of a Hamiltonian system,
Regul. Chaotic Dyn., 18:6 (2013), 774--800.
\bibitem{Bol-Rab:revers}
Bolotin S. and  Rabinowitz P.H., A variational construction of
chaotic trajectories for a reversible Hamiltonian system. {\it J.
Differ. Equat.,} {\bf 48} (1998), 365--387.
              \bibitem{Bol-Tre:hyp}   Bolotin~S. and Treschev~D., Remarks on the definition of hyperbolic tori of Hamiltonian systems, Regul. Chaotic Dyn., 5:4 (2000), 401-412

\bibitem{Bounemoura} Bounemoura A. and Pennamen E.
     Instability for a priori unstable Hamiltonian systems: a dynamical
     approach. (English summary)
     Discrete Contin. Dyn. Syst.  32  (2012),  no. 3, 753--793.
     \bibitem{Buf-Sere}
     Buffoni, B. and Sere, E.,
     A Global Condition for Quasi-Random Behavior in a Class of Conservative Systems,
     Comm. Pure Appl. Math., 1996, vol.49, no.3, pp. 285--305.
\bibitem{Sin_bill}  Bunimovich L.A., Sinai Ya.G., Chernov N.I.
              Statistical properties of two-dimensional hyperbolic billiards.
              (Russian) Uspekhi Mat. Nauk 46 (1991), no. 4(280), 43--92, 192;
              translation in Russian Math. Surveys 46 (1991), no. 4, 47--106
\bibitem{CK} O. Castejon, V. Kaloshin,
    Random Iteration of Maps on a Cylinder and diffusive behavior
    Preprint. June 13, 2015. http://www.math.umd.edu/~vkaloshi/papers/random-iteration.pdf
\bibitem{Chen2004}  Y.-C. Chen,
    Anti-integrability in scattering billiards,
    Dyn. Syst. 19 (2004) 145-159.
\bibitem{Chen2010} Chen, Yi-Chiuan,
              On topological entropy of billiard tables with small inner scatterers.
              Adv. Math.  224  (2010),  no. 2, 432--460.
\bibitem{Cheng-diff1} Cheng Ch.-Q. and Yan J.,
              Existence of diffusion orbits in a priori unstable Hamiltonian systems.
              J. Differential Geom. 67 (2004), no. 3, 457--517.
\bibitem{Cheng-diff2} Cheng Ch.-Q. and Yan J.,
              Arnold diffusion in Hamiltonian systems: a priori unstable case
              J. Differential Geom. Volume 82, Number 2 (2009), 229-277.
\bibitem{CG}        Chierchia L. and Gallavotti G.,
              Drift and diffusion in phase space.
              Ann. Inst. Henri Poincar\'e, {\bf 60} (1), 1994 1--144.



\bibitem{DLS2}      Delshams A., de la Llave R. and Seara T.,
              A geometric mechanism for diffusion in Hamiltonian systems overcoming the large gap problem:
              heuristics and rigorous verification on a model.
              Mem. Amer. Math. Soc. 179 (2006), no. 844, viii+141 pp.
\bibitem{DLS_scat} A. Delshams, R. de la Llave, T.M. Seara, Geometric
              properties of the scattering map of a normally hyperbolic invariant manifold, Adv. Math. 217 (3) (2008) 1096--1153.
\bibitem{Delsh-Hug1} Delshams A. and Huguet G.
              Geography of resonances and Arnold diffusion in a priori unstable Hamiltonian systems.
              Nonlinearity  22  (2009),  no. 8, 1997--2077.
\bibitem{Delsh-Hug2} Delshams A. and Huguet G.
              A geometric mechanism of diffusion: rigorous verification in a priori unstable Hamiltonian systems. (English summary)
              J. Differential Equations  250  (2011),  no. 5, 2601--2623.
 \bibitem{Deng}             B.~Deng, The Shilnikov problem, exponential expansion, strong
$\lambda$-lemma, $C^1$-linearization and homoclinic bifurcation.
{\it J.~Differ.~Equat.}, {\bf 79} (1989), 189--231.
 \bibitem{Devaney}              Devaney, R. L., Homoclinic Orbits in Hamiltonian Systems, J.
Differential Equations, 1976, vol.21, no.2, pp. 431--438.

\bibitem{Dua}    Duarte P.,
              Plenty of elliptic islands for the standard family
              of area preserving maps.
              {\it Ann. Inst. H. Poincare Anal. Non Lineaire}
              11 (1994), no. 4, 359--409.
\bibitem{East-Meiss}  R. W. Easton, J.D. Meiss, G. Roberts, Drift by
Coupling to an Anti-Integrable Limit, Physica D 156, 201--218.
(2001).
\bibitem{Fontich}   Fontich E. and Martin P.
              Arnold diffusion in perturbations of analytic integrable Hamiltonian systems.
              Discrete Contin. Dynam. Systems 7 (2001), no. 1, 61--84.
\bibitem{GGM}       Gallavotti G., Gentile G., Mastropietro V.
              Hamilton-Jacobi equation and existence of heteroclinic chains in three time scales systems.
              Nonlinearity 13 (2000) 323-340

\bibitem{Gel-Tur}
V.~Gelfreich and D.~Turaev, Unbounded energy growth in Hamiltonian
systems with a slowly varying parameter. {\it Comm. Math. Phys.}

\bibitem{GR}  Gidea M. and Robinson C. Diffusion along transition chains of
              invariant tori and Aubry-Mather sets. Ergodic Theory Dynam. Systems  33  (2013),  no. 5, 1401--1449.

\bibitem{Gole}
Gole C., Symplectic Twist Maps. Advanced Series in Nonlinear
Dynamics vol 18, Singapore: World Scientific (2001)
{\bf 283} (2008),  769--794.
\bibitem{KZ} M. Guardia, V. Kaloshin, J. Zhang
      A second order expansion of the separatrix map
      for trigonometric perturbations of a priori unstable systems
      Preprint. September 30, 2015, http://arxiv.org/pdf/1503.08301v2.pdf

\bibitem{Kalo-Levi} Kaloshin V. and Levi M.
              Geometry of Arnold diffusion.
              SIAM Rev.  50  (2008),  no. 4, 702--720.
\bibitem{Kal-Zhang}
      Kaloshin, V. and Zhang, K., Partial Averaging and Dynamics of
      the Dominant Hamiltonian, with Applications to Arnold Diffusion,
      arXiv:1410-1844 (2014).

\bibitem{KZZ} V. Kaloshin, J. Zhang, K. Zhang
      Normally Hyperbolic Invariant Laminations and diffusive behaviour for the generalized Arnold example away from resonances.
      Preprint, April 6, 2015
      http://www.math.umd.edu/~vkaloshi/papers/nhil-rand-model.pdf
\bibitem{KH}  Katok A., Hasselblatt B.
      Introduction to the Modern Theory of Dynamical Systems.
      Cambridge University Press, 1997.
      \bibitem{Knauf}     Klein, M., Knauf, A.,         Classical Planar Scattering by Coulombic
Potentials. Springer (1992)

\bibitem{Knill} Knill, O.,
      Topological entropy of standard type monotone twist maps.
      Trans. Amer. Math. Soc.  348  (1996),  no. 8, 2999--3013.
\bibitem{Koz-Tre}   Kozlov V. and Treschev D.,
         Billiards: a genetic introduction to the dynamics of sysetms with impacts.
         Translations of Mathematical Monographs, vol. 89, AMS, 1991.
\bibitem{Lerman}  Lerman, L. M.,
         Complex Dynamics and Bifurcations in a Hamiltonian System Having a Transversal Homoclinic Orbit to a Saddle Focus,
         Chaos, 1991, vol.1, no. 2, pp. 174--180.
\bibitem{mackmei}   MacKay~R.~S., Meiss~J.~D.,
         Cantori for symplectic maps near the anti-integrable limit.
         {\it Nonlinearity} {\bf 5}, V.~149, 1992, 1--12.
         \bibitem{Mather}
Mather, J.\ N., Existence of quasi-periodic orbits for twist
homeomorphisms of the annulus, {\it Topology}, {\bf 21}, (1982),
457--467.
         \bibitem{McDuff}
D.~McDuff and D.~Salamon, {\it Introduction to Symplectic
Topology.} Oxford Mathematical Monographs. The Clarendon Press,
Oxford University Press, New York, 1998.
\bibitem{Nitec}     Nitecki~Z.
         Differentiable Dynamics,
         The MIT Press, Cambridge, Massachusetts and London, England, 1971.
\bibitem{Pif}    Piftankin G.N. 2006
         Diffusion speed in the Mather problem
         {\it Nonlinearity} {\bf 19},  2617-2644
\bibitem{PT}  Piftankin G.N. and Treschev D.V.
         Separatrix map in Hamiltonian systems. (Russian)
         Uspekhi Mat. Nauk {\bf 62}, 2007, No. 2(374), 3--108;
         translation in Russian Math. Surveys {\bf 62}, 2007, No. 2, 219--322
\bibitem{Shil}  Shilnikov, L. P.,
         On a Poincar\'e-Birkhoff Problem,
         Math. USSR-Sb., 1967, vol.3, no. 3, pp. 353--371.
\bibitem{SSTC} Shilnikov L.P., Shilnikov A. L., Turaev D.V. and Chua L.O.
         Methods of qualitative theory in nonlinear dynamics. Part~II,  2001, Singapore: World Scientific.
         \bibitem{Smale} Smale S., Diffeomorphisms with many
         periodic points. Princeton University Press, 1965.
\bibitem{Tre_sm_res} Treschev~D.,
         Multidimensional separatrix maps
         J. Nonlin. Sci. {\bf 12}, No. 1, 2002, 27--58.
\bibitem{Tre_traj} Treschev D.,
         Trajectories in a neighborhood of asymptotic surfaces
         of a priori unstable Hamiltonian systems.
         Nonlinearity {\bf 15}, 2002, 2033--2052.
\bibitem{Tre2004} Treschev, D.
         Evolution of slow variables in a priori unstable Hamiltonian systems.
         Nonlinearity  17  (2004),  no. 5, 1803--1841.

\bibitem{Tre2012} Treschev D.,
         Arnold diffusion far from strong resonances in
         multidimensional a priori unstable Hamiltonian systems.
         Nonlinearity {\bf 25} (2012) 2717--2757.
         \bibitem{TZ}    Treschev D. and Zubelevich O.
         Introduction to the perturbation theory of Hamiltonian systems.
         Springer, 2009.



\bibitem{Turayev:RCD} Turaev, D. V.,   Hyperbolic Sets near Homoclinic Loops to a
 Saddle for Systems with a First Integral, Regul. Chaotic Dyn., 2014, 19
 (6), pp. 681-693.
 \bibitem{Tur-Shil}          Turaev, D. V. and Shilnikov, L. P., Hamiltonian Systems with
Homoclinic Saddle Curves, Soviet Math. Dokl., 1989, vol.39, no.1,
pp. 165--168.
\bibitem{Vas-But}
Vasilyeva, A.B., and V.F. Butuzov, Asymptotic expansion of the
singularly perturbed equations solutions. Nauka, Moscow  (1973).

\bibitem{Veselov}   Veselov A.P., Integrable mappings.
             Russian Math. Surveys {\bf 46}, 1991, No. 5, 1--51.
\bibitem{Walters} Walters P., An introduction to Ergodic Theory, p. 178,
             Springer Verlag: New-York, Heidelberg, Berlin, 1982.
\end{thebibliography}
\end{document}